\documentclass[11pt]{amsart}

\pdfoutput=1

\usepackage[a4paper]{geometry}
\usepackage{amsmath,amssymb,amsthm,amsfonts,color}
\usepackage{comment}
\usepackage{hyperref}
\usepackage{stmaryrd}
\usepackage{tikz}
\usepackage{bm}
\usepackage{accents}

\usepackage[shortlabels]{enumitem}
\usepackage{refcount}

\usepackage{threeparttable}
\usepackage{caption}
\theoremstyle{definition}
\newtheorem{Def}{Definition}

\newtheorem{Thm}[Def]{Theorem}
\newtheorem{Lem}[Def]{Lemma}
\newtheorem{Sublem}[Def]{Sublemma}

\newtheorem{Prop}[Def]{Proposition}

\newtheorem{Cor}[Def]{Corollary}
\newtheorem{Claim}{Claim}
\newtheorem{Subclaim}{Subclaim}
\newtheorem{Q}{Question}

\newtheorem{Conj}{Conjecture}
\newtheorem{Thm?}[Def]{Theorem?}

\newcommand{\ZFC}{\ensuremath{\mathsf{ZFC}}}
\newcommand{\ZF}{\ensuremath{\mathsf{ZF}}}
\newcommand{\AC}{\ensuremath{\mathsf{AC}}}
\newcommand{\DC}{\ensuremath{\mathsf{DC}}}
\newcommand{\AD}{\ensuremath{\mathsf{AD}}}
\newcommand{\BlAD}{\ensuremath{\mathsf{Bl}\text{-}\mathsf{AD}}}
\newcommand{\ADR}{\ensuremath{\mathsf{AD}_{\mathbb{R}}}}
\newcommand{\BlADR}{\ensuremath{\mathsf{Bl}\text{-}\mathsf{AD}_{\mathbb{R}}}}

\newcommand{\R}{\mathbb{R}}
\newcommand{\LL}{\mathrm{L}}

\renewcommand{\L}{\text{L}}

\begin{document}




\title[$\ADR$ and $\BlADR$]{The Axiom of Real Determinacy and the Axiom of Real Blackwell Determinacy}
\date{\today}

\author[D.\ Ikegami]{Daisuke Ikegami}
\address[D.\ Ikegami]{Institute of Logic and Cognition, Department of Philosophy, Sun Yat-sen University, Xichang Hall 602, 135 Xingang west street, Guangzhou, 510275 CHINA}

\email[D.\ Ikegami]{\href{mailto:ikegami@mail.sysu.edu.cn}{ikegami@mail.sysu.edu.cn}}

\author[W.\ H.\ Woodin]{W.\ Hugh Woodin}
\address[W.\ H.\ Woodin]{Department of Mathematics, Harvard University, Science Center 513, 1 Oxford Street, Cambridge, MA 02138 USA}

\email[W.\ H.\ Woodin]{\href{mailto:wwoodin@g.harvard.edu}{wwoodin@g.harvard.edu}}

\thanks{The research in this paper was partially conducted when the authors stayed in the Mittag-Leffler Institute for the semester program \lq\lq Mathematical Logic: Set Theory and Model Theory" in fall 2009. The authors are grateful to the Mittag-Leffler Institute for their hospitality and support.  This paper contains results from the first author's Ph.D. thesis~\cite{thesis}. The first author would like to thank Benedikt L\"{o}we for his supervision. The first author's research is partially supported by the Guangdong Philosophy and Social Science Foundation with grant number GD25YZX01. The second author acknowledges support from the US National Science Foundation with grant number DMS-0856201. 
The authors would like to thank the referee for their helpful comments which improved the presentation of this paper. }

\begin{abstract}
We show that the Axiom of Real Determinacy $\ADR$ and the Axiom of Real Blackwell Determinacy $\BlADR$ are equivalent in $\ZF$+$\DC$. This answers the question of L\"{o}we~\cite[Question~53]{MR2374765}. While we do not know whether they are equivalent in $\ZF$+$\AC_{\omega} (\mathbb{R})$, we show that the theories $\ZF$+$\ADR$ and $\ZF$+$\AC_{\omega} (\mathbb{R})$+$\BlADR$ are equiconsistent. 
\end{abstract}

\maketitle

\section{Introduction}


In this paper, we compare the stronger versions of determinacy of Gale-Stewart games and Blackwell games, i.e., the Axiom of Real Determinacy $\ADR$ and the Axiom of Real Blackwell Determinacy $\BlADR$. 

In 1953, Gale and Stewart~\cite{MR0054922} developed the general theory of infinite games, so-called {\em Gale-Stewart games}, which are two-player zero-sum infinite games with perfect information. The theory of Gale-Stewart games has been investigated by many logicians and now it is one of the main topics in set theory and it has connections with other topics in set theory as well as model theory and computer science. 

In 1928, John von Neumann proved his famous {\em minimax theorem} which is about finite games with imperfect information. Infinite versions of von Neumann's games were introduced by David Blackwell~\cite{Blackwell} where he proved the analogue of von Neumann's theorem for $\mathrm{G}_{\delta}$ sets of reals (i.e., $\undertilde{\mathbf{\Pi}}^0_2$ sets of reals). The games he introduced are called {\em Blackwell games} and they were called by him \lq \lq games with slightly imperfect information" in his paper~\cite{Blackwell2}. 

In 1998, Martin~\cite{MR1665779} proved that in most cases, Blackwell determinacy axioms follow from the corresponding determinacy axiom. Martin conjectured that they are equivalent, and many instances of equivalence have been shown 
(e.g., \cite{AD_Bl-AD} and Martin's proof of $\mathbf{\Pi}^1_1$ determinacy presented in \cite[Corollary 3.9]{MR2042937}). However, the general question, and in particular the most intriguing instance, namely whether $\AD$ and the axiom of Blackwell determinacy $\BlAD$ are equivalent, remains open.

De Kloet, L\"{o}we, and the first author~\cite{Bl_AD_R_sharp} turned to the other mentioned determinacy axiom, the stronger $\text{AD}_{\mathbb{R}}$ and its Blackwell analogue. They introduced the Axiom of Real Blackwell Determinacy $\text{Bl-AD}_{\mathbb{R}}$ and then proved that $\text{Bl-AD}_{\mathbb{R}}$ implies the existence of $\R^{\#}$, and that the consistency of $\ZF$+$\AC_{\omega}(\R)$+$\BlADR$ is strictly stronger than that of $\ZF$+$\AD$. 

In this paper, we investigate the relationship between $\ADR$ and $\BlADR$. 
Our main result is the following:
\begin{Thm}[Theorem~\ref{almost theorem}]\label{main theorem}
The axioms $\ADR$ and $\BlADR$ are equivalent in $\ZF$+$\DC$. 
\end{Thm}

Solovay~\cite{Solovay_AD_R} proved that the consistency of $\ZF$+$\DC$+$\ADR$ is strictly stronger than that of $\ZF$+$\ADR$. Hence assuming $\DC$ is not optimal when considering the equivalence between $\ADR$ and $\BlADR$. 
While we do not know whether they are equivalent in $\ZF$+$\AC_{\omega} (\mathbb{R})$, we show that they are equiconsistent: 
\begin{Thm}[Theorem~\ref{equicon}]\label{equiconsistency}
The theories $\ZF$+$\ADR$ and $\ZF$+$\AC_{\omega} (\R)$+$\BlADR$ are equiconsistent.
\end{Thm}

In \S\,\ref{sec:basics}, we introduce Blackwell determinacy as well as other notions, definitions, theorems, and lemmas we need throughout the paper. 
In \S\,\ref{sec:Rsharp}, we give a simple proof of the existence of a fine normal measure on $\wp_{\omega_1}(\mathbb{R})$, which was originally proved by de Kloet, L\"{o}we, and the first author~\cite{Bl_AD_R_sharp}. The key lemma for the existence of such a measure is the determinacy of range-invariant sets from $\BlADR$ (Lemma~\ref{range-invariant}), which will also be used in later sections. 
In \S\,\ref{sec:infBorel}, we show that $\BlADR$ implies that every set of reals has the Baire property. 
In \S\,\ref{sec:infBorelcodes}, we show that $\BlADR$ implies that every set of reals is $\infty$-Borel, and that $\BlADR$ and $\DC$ imply that every set of reals is strongly $\infty$-Borel, a key step towards the equivalence between $\ADR$ and $\BlADR$ in $\ZF$+$\DC$. 
In \S\,\ref{sec:equiv}, we prove Theorem~\ref{main theorem}, i.e., the equivalence between $\ADR$ and $\BlADR$ in $\ZF$+$\DC$. 
In \S\,\ref{sec:Con}, we prove Theorem~\ref{equiconsistency}, i.e., the equiconsistency between $\ADR$ and $\BlADR$. 
In the last section, we raise some open problems.  


\section{Preliminaries}\label{sec:basics}


Throughout this paper, we work in $\ZF$+$\AC_{\omega}(\mathbb{R})$, where $\AC_{\omega}(\mathbb{R})$ states that for all $\omega$-sequences $(A_n \mid n < \omega)$ of nonempty sets of reals, there is a function $f \colon \omega \to \mathbb{R}$ such that for all $n < \omega$, $f( n) \in A_n$. 
This small fragment of the axiom of choice is necessary for the definitions of axioms of Blackwell determinacy. Using $\AC_{\omega}(\mathbb{R})$, one can develop the basics of measure theory. If we need more than $\ZF$+$\AC_{\omega}(\mathbb{R})$ for some definitions and statements, we explicitly mention the additional axioms.
We use standard notations from set theory and assume that the reader is familiar with forcing, descriptive set theory, and determinacy of Gale-Stewart games, whose basics can be found in e.g., \cite{Jech} and \cite{new_Moschovakis}. By {\em reals}, we mean elements of the Cantor space  and we use $2^{\omega}$ or $\mathbb{R}$ to denote the Cantor space. 

\subsection{Measure and category}\label{sec:measure and category}

In this subsection, we list some theorems on measure and category we will use throughout the paper. 

For each finite binary sequence $s$, $[s]$ denotes the set $\{ x \in 2^{\omega} \mid x \supseteq s \}$. The family $\{ [s] \mid s \in 2^{<\omega}\}$ forms a basis for the Cantor space. 

\begin{Thm}[Folklore?]\label{fact:category}
Suppose that every set of reals has the Baire property. Then the meager ideal on the Cantor space $2^{\omega}$ and the meager ideal on the Baire space $\omega^{\omega}$ are closed under well-ordered unions. Also, every function $f \colon 2^{\omega} \to 2^{\omega}$ is continuous on a comeager set. 
\end{Thm}

\begin{proof}
For the closure of the meager ideal under well-ordered unions, see e.g., \cite[(8.41)~Theorem and (8.49)~Exercise]{Kechris_Classical}. 

Let $f \colon 2^{\omega} \to 2^{\omega}$ be any function. We show that $f$ is continuous on a comeager set. 
By assumption, for each finite binary sequence $s$, the set $f^{-1}([s])$ has the Baire property. Hence for each $s$, there is an open set $U_s$ such that $f^{-1}([s]) \triangle U_s$ is meager, where $f^{-1}([s]) \triangle U_s$ is the symmetric difference between $f^{-1}([s])$ and $U_s$.  
Let $D = 2^{\omega} \setminus \bigcup_{s \in 2^{<\omega}} f^{-1}([s]) \triangle U_s$. Then $D$ is comeager in the Cantor space and $f$ is continous on $D$, as desired. 
\end{proof}

Let $X$ be a set with at least two elements. The topology of $X^{\omega}$ is given by the product topology where each coordinate (i.e., $X$) is seen as the discrete space. We say $\mu$ is a {\it Borel probability on $X^{\omega}$} if $\mu$ is a probability measure on the collection of Borel subsets of the space $X^{\omega}$. 
\begin{Thm}[Folklore?]\label{fact:LM}
Suppose that every set of reals is Lebesgue measurable. Then for every Borel probability $\mu$ on the Cantor space $2^{\omega}$, every set of reals is $\mu$-measurable and the ideal of $\mu$-null sets is closed under well-ordered unions. Also, there is no injection from $\omega_1$ to the reals. 
\end{Thm}

\begin{proof}
For $\mu$-measurability of every set of reals for every Borel probability $\mu$ on the Cantor space, see \cite[(17.41)~Theorem]{Kechris_Classical}. 
The arguments for the closure under well-ordered unions for the ideal of $\mu$-null sets are in parallel with those for the closure under well-ordered unions for the meager ideal assuming the Baire property for every set of reals, which can be found in \cite[(8.41)~Theorem and (8.49)~Exercise]{Kechris_Classical}. 
For the non-existence of an injection from $\omega_1$ to the reals, see \cite[Theorem~5]{Can_you_Solovay_Raisonnier}. 
\end{proof}

\subsection{Pointclasses}\label{sec:pointclasses}

In this subsection, we duscuss some basics of pointclasses. 
As with Borel sets, one often looks at the properties of a class of sets of reals rather those of a set of reals. Such classes are called pointclasses. 
A {\it pointclass} is the union of sets of subsets of $\omega^m \times \mathbb{R}^n$ for natural numbers $m \ge 0, n\ge1$. If $\Gamma$ is a pointclass, $\Gamma$ is called a {\it boldface pointclass} if it is closed under continuous preimages, i.e., for all natural numbers $m_1, m_2 \ge 0$ and $n_1, n_2 \ge 1$, and continuous functions $f\colon \omega^{m_1} \times \mathbb{R}^{n_1} \to \omega^{m_2} \times \mathbb{R}^{n_2}$, and a subset $A \in \Gamma$ of $\omega^{m_2} \times \mathbb{R}^{n_2}$, $f^{-1}(A)$ is also in $\Gamma$. Closure under recursive preimages is similarly defined with recursive functions. 

A pointclass $\Gamma$ is {\it $\omega$-parametrized} if for all natural numbers $m \ge 0 $ and $n\ge 1$ there is a subset $G^{m,n}$ of $\omega^{m+1} \times \mathbb{R}^n$ in $\Gamma$ such that for any subset $A$ of $\omega^m \times \mathbb{R}^n$ in $\Gamma$, there is a natural number $e$ such that $A = G^{m,n}_e = \{ (x,y ) \mid (e, x,y) \in G^{m,n}\}$. The following lemma is useful: Let $\Gamma$ be a pointclass and $x$ be a real. Then the pointclass $\Gamma (x)$ is the set of all sets $A$ such that there is a set $B \in \Gamma$ such that $A = B_x$. Set $\undertilde{\mathbf{\Gamma}} = \bigcup_{x \in \mathbb{R}} \Gamma (x)$. 

\begin{Lem}\label{good parametrization}
Suppose $\Gamma$ is an $\omega$-parametrized pointclass which is closed under recursive preimages. Then for each natural number $n\ge1$, there is a set $G^n \subseteq \mathbb{R} \times \mathbb{R}^n$ in $\Gamma$ such that the following hold:

\begin{enumerate}
\item For each $n\ge 1$, $G^n$ is universal for subsets of $\mathbb{R}^n$ in $\undertilde{\mathbf{\Gamma}}$, i.e., for any $A \subseteq \R^n$ in $\undertilde{\mathbf{\Gamma}}$, there is a real $x$ such that $A = G^n_x$, where $G^n_x = \{ \vec{y} \mid (x, \vec{y}) \in G^n \}$, 

\item For any $A \subseteq \mathbb{R}^n$ in $\Gamma$, there is a recursive real $x$ such that $A = G^n_x$, and

\item\label{para-item3} For all natural numbers $n, m \ge 1$, there is a recursive function $S^{n,m} \colon \mathbb{R} \times \mathbb{R}^n \to \mathbb{R}$ such that for any real $a$, $x \in \mathbb{R}^n$, and $y\in \mathbb{R}^m$, $G^{m+n} (a,x,y) \iff G^m (S^{n,m}(a,x),y)$. 
\end{enumerate}
\end{Lem}
\begin{proof}
See \cite[3H.1]{new_Moschovakis}. 
\end{proof}

We fix some notations for projections. For natural numbers $m\ge 0$ and $n\ge 1$ and a subset $A$ of $\omega \times \omega^m \times \mathbb{R}^{n}$, let $\exists^{\omega}A = \{ (x, y) \in \omega^{m} \times \mathbb{R}^n \mid (\exists e \in \omega) \ (e, x,y)\in A\}$ and $\forall^{\omega}A = \{ (x, y) \in \omega^{m} \times \mathbb{R}^n \mid (\forall e \in \omega) \ (e, x,y)\in A\}$. The sets $\exists^{\mathbb{R}}A $ and $\forall^{\mathbb{R}}A$ are defined in a similar way. A pointclass $\Gamma$ is {\it closed under $\exists^{\omega}$} if for any $A$ in $\Gamma$, $\exists^{\omega} A $ is in $\Gamma$. Closure under $\forall^{\omega}, \exists^{\mathbb{R}}$, and $\forall^{\mathbb{R}}$ is defined in a similar way. 

\begin{Def}\label{def:Spector}
A pointclass $\Gamma$ is a {\it Spector pointclass} if it satisfies the following:

\begin{enumerate}
\item It contains all the $\Sigma^0_1$ sets and it is closed under recursive substitutions, finite intersections and unions, $\exists^{\omega}$, and $\forall^{\omega}$, 

\item It is $\omega$-parametrized, 

\item It has the substitution property, and 

\item It has the prewellordering property.
\end{enumerate}
\end{Def}

For the definition of the substitution property and the basic theory of $\Gamma$-recursive functions, see \cite[3D \& 3G]{new_Moschovakis}. For the definition of the prewellordering property, see \cite[4B]{new_Moschovakis}. Typical examples of Spector pointclasses are $\Pi^1_1$ and $\Sigma^1_2$. Assuming the determinacy of all the projective sets, one can prove that $\Pi^1_{2n+1}$ and $\Sigma^1_{2n+2}$ are also Spector pointclasses for each natural number $n$. 

A typical example of Spector pointclasses closed under both $\exists^{\mathbb{R}}$ and $\forall^{\mathbb{R}}$ is the collection of {\it inductive sets}. 
Let $A$ be a set of reals. Let $\mathrm{IND}(A)$ be the set of all $\mathrm{pos}\Sigma^1_n(A)$-inductive sets of reals for some natural number $n\ge 1$. For the definition of $\mathrm{pos}\Sigma^1_n(A)$-inductive sets, see \cite[7C]{new_Moschovakis}. All we need on the theory of inductive sets is as follows:
\begin{Thm}\label{inductive}
For any set of reals $A$, $\mathrm{IND}(A)$ is the smallest Spector pointclass containing $A$ and closed under $\exists^{\mathbb{R}}$ and $\forall^{\mathbb{R}}$.
\end{Thm}
\begin{proof}
The argument is the same as \cite[7C.3]{new_Moschovakis}. 
\end{proof}

The following lemma is basic:
\begin{Lem}[Kleene]\label{lem:Kleene}
Let $\Gamma$ be a Spector pointclass. Then there is a $\Gamma$-recursive partial function $(a,x) \mapsto \{ a \} (x)$ from $\R \times \R$ to $\R$ on its domain with the following property: A partial function $f$ from $\R$ to $\R$ is $\undertilde{\mathbf{\Gamma}}$-recursive on its domain if and only if there is some real $a$ such that for all reals $x$, if $f(x)$ is defined, then $f(x) = \{ a \} (x)$. 
\end{Lem}

\begin{proof}
See \cite[7A.1]{new_Moschovakis}. 
\end{proof}

We will use the following general form of Kleene's Recursion Theorem for Spector pointclasses in \S\,\ref{sec:Moschovakis}:
\begin{Thm}[Kleene's Recursion Theorem]\label{recursion theorem}
Let $\Gamma$ be a Spector pointclass and suppose $f$ is a partial function from $\mathbb{R} \times \mathbb{R}$ to $\mathbb{R}$ which is $\undertilde{\mathbf{\Gamma}}$-recursive on its domain. Then there exists a fixed real $a^*$ such that for all reals $x$, if $f(a^*, x)$ is defined, then $f (a^*, x) = \{a^*\}(x)$, where the notation $\{ a^* \} (x)$ is from Lemma~\ref{lem:Kleene}.  
\end{Thm}

\begin{proof}
See \cite[7A.2]{new_Moschovakis}. 
\end{proof}

\subsection{Blackwell games}\label{sec:Blackwell games}

In this subsection, we introduce {\em Blackwell games}, which are infinite games with imperfect information, and compare Blackwell games with Gale-Stewart games. 

We start with the definition of Blackwell games.\footnote{Our definitions of Blackwell games and Blackwell determinacy are different from the original ones given by Blackwell~\cite{Blackwell2} where Blackwell determinacy is formulated as an extension of von Neumann's minimax theorem, but our formulation is equivalent to the original one when it is about the Cantor space (i.e., when $X=2$). For the original formulation of Blackwell games and Blackwell determinacy, see, e.g., \cite[\S\,3 \& \S\,5]{MR2374765}.\label{Blackwell original}}
Let $X$ be a set with at least two elements and assume $\AC_{\omega} (X^{\omega})$. The topology of $X^{\omega}$ is given by the product topology where each coordinate (i.e., $X$) is seen as the discrete space. For each finite sequence $s$ of elements in $X$, $[s]$ denotes the set $\{ x \in X^{\omega} \mid x \supseteq s \}$ and it is a basic open set in the topological space $X^{\omega}$. In Blackwell games, players choose probabilities on $X$ instead of elements of $X$ and with those probabilities, one can derive a Borel probability on $X^{\omega}$, i.e., a measure assigning probability to each Borel subset of $X^{\omega}$. Player I wins if the probability of a given payoff set is $1$ and player II wins if the probability of the payoff set is $0$. Let us formulate this in detail. 
We write $X^{\mathrm{Even}}$ for the set of all finite sequences of elements of $X$ of even length and $X^{\mathrm{Odd}}$ for the set of all finite sequences of elements of $X$ of odd length. 
\begin{Def}
A {\it mixed strategy for player \text{I}} is a function $\sigma \colon X^{\mathrm{Even}} \to \mathrm{Prob}_{\omega} (X)$, where $\mathrm{Prob}_{\omega}(X)$ is the set of functions $\mu \colon X \to [0,1]$ with $\sum_{x \in X} \mu (x) = 1$.\footnote{We use $\mathrm{Prob}_{\omega} (X)$ to denote such functions because they are the same as Borel probabilities $\mu$ on $X$ with countable support, i.e., there is a countable subset $A$ of $X$ with $\mu (A) = 1$.} A {\it mixed strategy for player \text{II}} is a function $\tau \colon X^{\mathrm{Odd}} \to \mathrm{Prob}_{\omega} (X)$. 

Given mixed strategies $\sigma$, $\tau$ for player I and II respectively, let $\nu (\sigma, \tau) \colon X^{<\omega} \to \mathrm{Prob}_{\omega} (X)$ be as follows: For each finite sequence $s$ of elements of $X$,
\begin{align*}
\nu (\sigma, \tau) (s) = \begin{cases}
\sigma (s) \ \ \ \text{ if $ s \in X^{\mathrm{Even}}$,}\\
\tau (s) \ \ \ \text{ if $ s \in X^{\mathrm{Odd}}$.}
\end{cases}
\end{align*}
For each finite sequence $s$ of elements of $X$, define
\begin{align*}
\mu_{\sigma,\tau}([s]) = \prod_{i=0}^{\mathrm{lh}(s)-1} \nu(\sigma, \tau) (s {\upharpoonright} i ) \ \bigl(s(i) \bigr). 
\end{align*}
Recall that $[s]$ denotes the set of $x \in X^{\omega}$ such that $x \supseteq s$ and these sets are basic open sets in the space $X^{\omega}$. 
With the help of $\AC_{\omega} (X^{\omega})$, one can uniquely extend $\mu_{\sigma, \tau}$ to a Borel probability on $X^{\omega}$, i.e., the probability whose domain is the set of all Borel sets in the space $X^{\omega}$. 
Let us also use $\mu_{\sigma, \tau}$ for denoting this Borel probability. 

We sometimes consider a special mixed strategy so-called a {\it blindfolded strategy}: For a $y \in X^{\omega}$, we write $\sigma_y$ for the mixed strategy for player I concentrating on $y$ of measure one, i.e., for all $s \in X^{\mathrm{Even}}$ and $x\in X$, set $\sigma_y (s) (x) = 1$ if $y \bigl(\mathrm{lh}(s) /2 \bigr) = x$ where $\mathrm{lh}(s)$ is the length of $s$, otherwise $\sigma_y (s) (x) =0$. One can define $\tau_y$, the mixed strategy for player II concentrating on $y$ of measure one in a similar manner. We write $\mu_{y , \tau}$ for $\mu_{\sigma_y , \tau}$ and $\mu_{\sigma , y}$ for $\mu_{\sigma , \tau_y}$. 

Let $A$ be a subset of $X^{\omega}$. A mixed strategy $\sigma$ for player I is {\it optimal in $A$} if for any mixed strategy $\tau$ for player II, $A$ is $\mu_{\sigma, \tau}$-measurable and $\mu_{\sigma, \tau} (A) = 1$. A mixed strategy $\tau $ for player II is {\it optimal in $A$} if for any mixed strategy $\sigma$ for player I, $A$ is $\mu_{\sigma, \tau}$-measurable and $\mu_{\sigma, \tau}(A) = 0$. A set $A$ is {\it Blackwell-determined} if one of the players has an optimal strategy in $A$. The axiom $\BlAD_X$ states that every subset of $X^{\omega}$ is Blackwell-determined. We write $\BlAD$ for $\BlAD_{\omega}$. 
\end{Def}
Note that since there is a bijection between $\mathbb{R}$ and $\mathbb{R}^{\omega}$, $\AC_{\omega}(\mathbb{R})$ implies $\AC_{\omega} (\mathbb{R}^{\omega})$ and hence one can formulate Blackwell games in $\mathbb{R}^{\omega}$ and $\BlAD_{\mathbb{R}}$ within $\ZF$+$\AC_{\omega}(\mathbb{R})$.

The following is a useful observation:
\begin{Prop}
${}$
\begin{enumerate}\label{easy_Blackwell}
\item Let $X, Y$ be sets with at least two elements and suppose that there is an injection from $X$ to $Y$ and assume $\AC_{\omega}(Y^{\omega})$. Then $\BlAD_Y$ implies $\BlAD_X$. In particular, $\BlAD_{\mathbb{R}}$ implies $\BlAD$. 

\item The axioms $\BlAD$ and $\BlAD_2$ are equivalent. 
\end{enumerate}
\end{Prop}

\begin{proof}
The first item is easy to see. For the second item, see \cite[Corollary~4.4]{MR1939050}. 
\end{proof}

As for Gale-Stewart games, one could ask what kind of subsets of $X^{\omega}$ are Blackwell-determined for a set $X$. After proving that every $\mathrm{G}_{\delta}$ subset of the Cantor space is Blackwell-determined, Blackwell asked whether every Borel subset of the Cantor space is determined. 
It was Donald Martin who found a general connection between the determinacy of Gale-Stewart games and Blackwell determinacy.\footnote{In \cite{MR1665779}, Martin proved the Blackwell determinacy in the original formulation as mentioned in Footnote~\ref{Blackwell original}, not in our formulation.}
\begin{Thm}[Martin]\label{det implies Blackwell det}
Let $X$ be a set with at least two elements and assume $\AC_{\omega} (X^{\omega})$. If there is a winning strategy for player I (resp., II) in a subset $A$ of $X^{\omega}$, then there is an optimal strategy for player I (resp., II) in $A$. In particular, $\AD$ implies that $\BlAD$ and $\AD_{\mathbb{R}}$ implies that $\BlAD_{\mathbb{R}}$. 
\end{Thm}

\begin{proof}
Given a strategy $\sigma$ for player I (resp., II), one can naturally translate $\sigma$ into a mixed strategy $\hat{\sigma}$ for player I (resp., II) by setting $\hat{\sigma}(s)$ to be the Dirac measure concentrating on $\sigma(s)$. It is easy to see that if $\sigma$ is winning in $A$, then $\hat{\sigma}$ is optimal in $A$. 
\end{proof}

Since every Borel set is determined in $\ZFC$, every Borel subset of the Cantor space is Blackwell-determined in $\ZFC$ and this answers the question of Blackwell. After proving Theorem~\ref{det implies Blackwell det}, Martin conjectured the following:
\begin{Conj}[Martin]\label{Martin-conjecture}
$\BlAD$ implies $\AD$. Hence $\AD$ and $\BlAD$ are equivalent. 
\end{Conj}

This conjecture is still not known to be true. The best known result toward $\AD$ from $\BlAD$ is as follows: A set of reals $A$ is {\it Suslin} if there is a tree $T$ on $2 \times \gamma $ for some ordinal $\gamma$ such that $A = \text{p}[T] = \{ x \in 2^{\omega} \mid (\exists f \in \gamma^{\omega}) \ (x, f) \in [T] \}$, where $[T]$ is the set of all infinite paths through $T$. A set of reals is {\it co-Suslin} if its complement is Suslin. 
\begin{Thm}[Martin, Neeman, and Vervoort]\label{Suslin-co-Suslin-det}
Assume $\BlAD$. Then every Suslin and co-Suslin set of reals is determined. 
\end{Thm}

\begin{proof}
See \cite[Lemma~4.1]{AD_Bl-AD}.
\end{proof}

Together with the following result, one can establish the equiconsistency between $\AD$ and $\BlAD$:
\begin{Thm}[Kechris and Woodin]
Assume that every Suslin and co-Suslin set of reals is determined. Then $\AD^{\mathrm{L}(\mathbb{R})}$ holds. 
\end{Thm}

\begin{proof}
See \cite{MR699440}.
\end{proof}

\begin{Cor}[Martin, Neeman, and Vervoort]\label{equi-cons_AD_Bl-AD}
In $\mathrm{L}(\mathbb{R})$, $\AD$ and $\BlAD$ are equivalent. In particular, $\ZF$+$\AD$ and $\ZF$+$\AC_{\omega}(\mathbb{R})$+$\BlAD$ are equiconsistent. 
\end{Cor}

As a corollary of Theorem~\ref{Suslin-co-Suslin-det}, we obtain strong closure properties of the pointclass of all the Suslin sets:
\begin{Cor}\label{cor:Suslin}
Assume $\BlAD$ and $\DC$. Let $\Gamma$ be the collection of all the Suslin sets. Then $\Gamma$ is a boldface pointclass closed under countable unions and intersections, $\exists^{\mathbb{R}}$, and $\forall^{\mathbb{R}}$. 
\end{Cor} 

\begin{proof}
For the arguments for the assertion that $\Gamma$ is a boldface pointclass as well as the closure of $\Gamma$ under countable unions and intersections, and $\exists^{\mathbb{R}}$, see \cite[2B.2]{new_Moschovakis}. 
For the closure of $\Gamma$ under $\forall^{\mathbb{R}}$, let $\Delta$ be the pointclass of all the Suslin and co-Suslin sets. Then by $\BlAD$ and Theorem~\ref{Suslin-co-Suslin-det}, every set of reals in $\Delta$ is determined. Now by the second periodicity theorem (see \cite[6C.3]{new_Moschovakis}), $\Gamma$ is closed under $\forall^{\mathbb{R}}$. 
\end{proof}

Also, $\BlAD$ has some consequence on regularity properties:
\begin{Thm}[Vervoort]\label{Blackwell Lebesgue}
Assume $\BlAD$. Then every set of reals is Lebesgue measurable. 
\end{Thm}

\begin{proof}
See \cite[Theorem~4.3]{MR1481790}. 
\end{proof}

By Theorem~\ref{fact:LM}, the following holds:
\begin{Cor}\label{cor:BlAD-measure}
Assume $\BlAD$. Then for every Borel probability $\mu$ on the Cantor space, every set of reals is $\mu$-measurable and the ideal of $\mu$-null sets is closed under well-ordered unions. Also, there is no injection from $\omega_1$ to the reals. 
\end{Cor}

In \S\,\ref{sec:infBorel}, we discuss the connection between $\BlADR$ and the Baire property. 

It is not difficult to see that if finite games are Blackwell determined, then they are determined. As a corollary, one can obtain the following:
\begin{Thm}[L\"{o}we]\label{uniformization}
Assume $\BlADR$. Then every relation on the reals can be uniformized by a function. In particular, $\DC_\R$ holds. 
\end{Thm}

\begin{proof}
See \cite[Theorem~9.3]{MR2374765}.
\end{proof}

Since there is a relation on the reals which cannot be uniformized by a function in $\mathrm{L}(\mathbb{R})$, $\BlADR$ does not hold in $\mathrm{L}(\mathbb{R})$. Since $\AD$ implies $\AD^{\LL(\mathbb{R})}$, $\AD$ does not imply $\BlADR$ in $\ZF$. 

Moreover, the consistency of $\BlADR$ is strictly stronger than that of $\AD$:
\begin{Thm}[de Kloet, L\"{o}we, and I.]\label{theorem:Rsharp}
The axiom $\BlADR$ implies that $\mathbb{R}^{\#}$ exists.
\end{Thm}

\begin{proof}
See \cite{Bl_AD_R_sharp}.
\end{proof}

Since $\BlADR$ implies $\BlAD$ by Proposition~\ref{easy_Blackwell} and $\BlAD$ implies $\AD^{\mathrm{L}(\mathbb{R})}$ by Corollary~\ref{equi-cons_AD_Bl-AD}, we obtain the following:
\begin{Cor}[de Kloet, L\"{o}we, and I.]
The consistency of $\ZF$+$\AC_{\omega}(\mathbb{R})$+$\BlADR$ is strictly stronger than that of $\ZF$+$\AD$. 
\end{Cor}

In the proof of Theorem~\ref{theorem:Rsharp}, they proved and used the following:
\begin{Thm}[de Kloet, L\"{o}we, and I.]\label{R sharp}
Assume $\BlADR$. Then there is a fine normal measure on $\wp_{\omega_1}(\mathbb{R})$. 
\end{Thm}

In \S\,\ref{sec:Rsharp}, we give a simpler proof of Theorem~\ref{R sharp}. 

\subsection{Stone space $\mathrm{St}(\mathbb{P})$ and $\mathbb{P}$-Baireness}\label{subsec:P}

In this subsection, we introduce the Stone space $\mathrm{St} (\mathbb{P})$ and $\mathbb{P}$-Baireness, a regularity property for sets of reals, for preorders $\mathbb{P}$. 



For a preorder $\mathbb{P}$, the {\em Stone space} of $\mathbb{P}$ (denoted by $\mathrm{St}(\mathbb{P})$) is the set of all ultrafilters on $\mathbb{P}$ equipped with the topology generated by $\{ O_p \mid p \in \mathbb{P}\}$, where $O_p = \{ u \in \mathrm{St}(\mathbb{P}) \mid u \ni p \}$. 
For example, if $\mathbb{P}$ is Cohen forcing $\mathbb{C} = (2^{<\omega}, \supseteq)$, then $\mathrm{St}(\mathbb{C})$ is homeomorphic to the Cantor space $ 2^{\omega}$.

Let $X, Y$ be topological spaces. A function $f\colon X \to Y$ is {\it Baire measurable} if for any open set $U$ in $Y$, $f^{-1}(U)$ has the Baire property in $X$.
Baire measurable functions are the same as continuous functions modulo meager sets: Let $X, Y$ be topological spaces and assume that $Y$ is second countable. Then a function $f\colon  X\to Y$ is Baire measurable if and only if there is a comeager set $D$ in $X$ such that $f{\upharpoonright} D$ is continuous.

There is a natural correspondence between Baire measurable functions from $\mathrm{St}(\mathbb{P})$ to the reals and $\mathbb{P}$-names for reals:
\begin{Lem}[Feng, Magidor, and Woodin]\label{name-Baire measurable}

Let $\mathbb{P}$ be a preorder.

\begin{enumerate}
\item If $f\colon \mathrm{St}(\mathbb{P}) \to 2^{\omega}$ is a Baire measurable function, then 
\begin{align*}
\tau_f = \bigl\{(\check{n}, p) \mid O_p \setminus \{ u \in \mathrm{St}(\mathbb{P}) \mid f(u)(n) = 1 \} \text{ is meager} \bigr\}
\end{align*}
is a $\mathbb{P}$-name for a real.

\item Let $\tau$ be a $\mathbb{P}$-name for a real. Define $f_{\tau}$ as follows: For $u\in \mathrm{St}(\mathbb{P})$ and $n \in \omega$,
\begin{align*}
f_{\tau}(u) (n) = 1 \iff (\exists p \in u)\ p\Vdash \check{n} \in \tau.
\end{align*}
Then the domain of $f_{\tau}$ is comeager in $\mathrm{St}(\mathbb{P})$ and $f_{\tau}$ is continuous on the domain. Hence it can be uniquely extended to a Baire measurable function from $\mathrm{St}(\mathbb{P})$ to the reals modulo meager sets.

\item If $f\colon \mathrm{St}(\mathbb{P}) \to 2^{\omega}$ is a Baire measurable function, then $f_{(\tau_f)}$ and $f$ agree on a comeager set in $\mathrm{St}(\mathbb{P})$.
Also, if $\tau$ is a $\mathbb{P}$-name for a real, then $\Vdash \tau_{(f_{\tau})} = \tau$.
\end{enumerate}
\end{Lem}

\begin{proof}
See \cite[Theorem~3.2]{MR1233821}. 
\end{proof}

We now define the property $\mathbb{P}$-Baireness. Let $\mathbb{P}$ be a preorder and $A$ be a set of reals. We say $A$ is {\it $\mathbb{P}$-Baire} if for any Baire measurable function $f \colon \mathrm{St}(\mathbb{P}) \to 2^{\omega}$, $f^{-1} (A)$ has the Baire property in $\mathrm{St}(\mathbb{P})$. 
For example, if $\mathbb{P}$ is Cohen forcing $\mathbb{C}$, then $A$ is $\mathbb{P}$-Baire if and only if for all continuous functions $f \colon 2^{\omega} \to 2^{\omega}$, the preimage $f^{-1}(A)$ has the Baire property in the Cantor space. 
It is easy to see that every Borel set of reals is $\mathbb{P}$-Baire for any $\mathbb{P}$. 


If $\mathbb{P}$ is an element of $H_{\omega_1}$, i.e., the transitive closure of $\{ \mathbb{P} \}$ is countable in $V$, then $\text{St}(\mathbb{P})$ is essentially the same as $\text{St}(\mathbb{C})$ where $\mathbb{C}$ is Cohen forcing, hence the Cantor space $2^{\omega}$ by the following lemma:
\begin{Lem}\label{Baireness_Baire}
If $i \colon \mathbb{P} \to \mathbb{Q}$ is dense for preorders $\mathbb{P}$ and $\mathbb{Q}$ in $H_{\omega_1}$, then $\text{St}(\mathbb{P})$ and $\text{St}(\mathbb{Q})$ are isomorphic as Baire spaces, i.e., there is a topological homeomorphism between a comeager set in $\text{St}(\mathbb{P})$ and a comeager set in $\text{St}(\mathbb{Q})$. 
\end{Lem}

\begin{proof}
Since $\mathbb{P}$ and $\mathbb{Q}$ are in $H_{\omega_1}$, so is $i$. Pick a real $x$ coding $(\mathbb{P}, \mathbb{Q}, i)$ and take a countable transitive model $M$ of $\ZF$-$\mathsf{P}$ with $x \in M$, where $\mathsf{P}$ denotes the Power Set Axiom, (e.g., consider $M = \mathrm{L}_{\alpha} [x]$ for a suitable countable $\alpha$). 
Then $i$ induces a natural bijection between $\mathbb{P}$-generic filters over $M$ and $\mathbb{Q}$-generic filters over $M$. Since $M$ is countable, the set of all $\mathbb{P}$-generic filters over $M$ is comeager in $\mathrm{St}(\mathbb{P})$ and the same holds for $\mathbb{Q}$. This natural bijection witnesses the conclusion. 
\end{proof}

By Theorem~\ref{fact:category}, the following also holds:
\begin{Cor}\label{cor:category}
Suppose that every set of reals has the Baire property. Let $\mathbb{P}$ be any preorder in $H_{\omega_1}$.  Then every set of reals is $\mathbb{P}$-Baire, and every function $f \colon \mathrm{St}(\mathbb{P}) \to 2^{\omega}$ is continuous on a comeager set in $\mathrm{St}(\mathbb{P})$, i.e., Baire measurable.
\end{Cor}

\subsection{$\infty$-Borel codes}\label{sec:Borel codes}

In this subsection, we introduce infinitary Borel codes and discuss their basic properties. 
Infinitary Borel codes ($\infty$-Borel codes) are a transfinite generalization of Borel codes: Let $\mathcal{L}_{\infty, 0} (\{\mathbf{a}_{n}\}_{n\in \omega})$ be the language allowing arbitrary many conjunctions and disjunctions and no quantifiers with atomic sentences $\mathbf{a}_{n}$ for each $n \in \omega$. The {\it $\infty$-Borel codes} are the sentences in $\mathcal{L}_{\infty, 0} (\{\mathbf{a}_{n}\}_{n \in \omega})$ belonging to any $\Gamma$ such that
\begin{itemize}
\item the atomic sentence $\mathbf{a}_{n}$ is in $\Gamma$ for each $n\in \omega$,

\item if $\phi $ is in $\Gamma$, then so is $\neg \phi$, and

\item if $\alpha $ is an ordinal and $\langle \phi_{\beta} \mid \beta < \alpha\rangle $ is a sequence of sentences each of which is in $\Gamma$, then $\bigvee_{\beta < \alpha} \phi_{\beta} $ is also in $\Gamma$.
\end{itemize}
To each $\infty$-Borel code $\phi$, we assign a set of reals $B_{\phi}$ in the same way as decoding Borel codes:
\begin{itemize}
\item if $\phi = \mathbf{a}_{n}$, then $B_{\phi} = \{ x \in 2^{\omega} \mid x(n) = 1\}$,

\item if $\phi = \neg \psi$, then $B_{\phi} = 2^{\omega} \setminus B_{\psi}$, and 

\item if $\phi = \bigvee_{\beta < \alpha} \psi_{\beta}$, then $B_{\phi} = \bigcup_{\beta < \alpha} B_{\psi_{\beta}}$. 
\end{itemize}
A set of reals $A$ is called {\it $\infty$-Borel} if there is an $\infty$-Borel code $\phi$ such that $A = B_{\phi}$. 

As Borel codes, one can regard $\infty$-Borel codes as well-founded trees on an ordinal with atomic sentences $\mathbf{a}_{n}$ on terminal nodes and decode them by assigning sets of reals to each node recursively from terminal nodes. (If a node has only one successor, then it means \lq \lq negation" and if a node has more than one successors, then it means \lq \lq disjunction".) The only difference between Borel codes and $\infty$-Borel codes is that trees are on $\omega$ for Borel codes while trees are on an ordinal for $\infty$-Borel codes. From this visualization, it is easy to see that the statement \lq \lq $\phi$ is an $\infty$-Borel code" is absolute between transitive models of $\ZF$. 
Also, it is easy to see that the statement \lq \lq a real $x$ is in $B_{\phi}$" is absolute between transitive models of $\ZF$. 

Considering the visualization of $\infty$-Borel codes described in the last paragraph, we often identify an $\infty$-Borel code $\phi$ with a pair $(S, c)$ where $S$ is a tree on some ordinal $\gamma$ and $c$ is a function from the set of terminal nodes of $S$ to natural numbers: If $c(t) = n$, then the atomic sentence $\mathbf{a}_{n}$ is attached on the terminal node $t$ of $S$. 
If $(S,c)$ is an $\infty$-Borel code, then we often write $S$ to indicate the pair $(S,c)$. 

Recall that $\Theta = \sup\,  \{ \gamma \mid \text{There is a surjection from $\R$ to $\gamma$} \}$. 
By the definition of $\Theta$, it is easy to see that for any $\infty$-Borel code $S$, there is an $\infty$-Borel code $T$ such that $B_S = B_T$ and the underlying tree of the code $T$ is on some ordinal $\gamma$ which is less than $\Theta$.

The following characterization of $\infty$-Borel sets is very useful:
\begin{Lem}\label{equiv_code}
Let $A$ be a set of reals. Then the following are equivalent:
\begin{enumerate}
\item $A$ is $\infty$-Borel, and

\item there are a formula $\phi$ in the language of set theory and a set $S$ of ordinals such that for each real $x$, 
\begin{align*}
x \in A \iff \LL[S,x] \vDash \phi (x).
\end{align*}
\end{enumerate}
\end{Lem}

\begin{proof}
See \cite[Theorem~8.7]{ADplus}. 
\end{proof}

Standard examples of $\infty$-Borel sets are Suslin sets. Recall that a set of reals $A$ is {\em Suslin} if there are an ordinal $\gamma$ and a tree $T$ on $2 \times \gamma$ such that $A = \mathrm{p} [T]$, where $\mathrm{p} [T]$ is the projection of $[T]$ to the first coordinate, i.e., 
\begin{align*}
\mathrm{p} [T] = \{ x \in 2^{\omega} \mid (\exists f \in \gamma^{\omega})\ (x,f) \in [T]\}.
\end{align*}
By Lemma~\ref{equiv_code}, every Suslin set is $\infty$-Borel. Assuming the Axiom of Choice, it is easy to see that every set of reals is Suslin, in particular $\infty$-Borel. Hence the property $\infty$-Borelness is trivial in the $\ZFC$ context while it is nontrivial and powerful in a determinacy world.

\subsection{Wadge games and Blackwell determinacy}

In this subsection, we discuss a weak version of Wadge's Lemma which follows from $\BlAD$: Let $A$ be a set of reals. For a natural number $n \ge 1$, a set of reals $B$ is {\it $\undertilde{\mathbf{\Sigma}}^1_n$ in $A$} if $B$ is definable by a $\Sigma^1_n$ formula in the structure $\mathcal{A}^2_A = (\omega, \wp (\omega) , 0, 1, + , \times , \in , A)$, the second order structure with $A$ as an unary predicate, using  a parameter $x$ for some real $x$. A set of reals $B$ is {\it projective in $A$} if $B$ is $\undertilde{\mathbf{\Sigma}}^1_n$ in $A$ for some $n\ge 1$. 
\begin{Lem}[Weak version of Wadge's Lemma]\label{weak Wadge}
Assume $\BlAD$. Then for any two sets of reals $A$ and $B$, either $A$ is $\undertilde{\mathbf{\Sigma}}^1_2$ in $B$, or $B$ is $\undertilde{\mathbf{\Sigma}}^1_2$ in $A$. 
\end{Lem}

\begin{proof}
We consider the Wadge game $\mathcal{G}_{\mathrm{W}}(A,B)$. By $\BlAD$, one of the players has an optimal strategy in $\mathcal{G}_{\mathrm{W}}(A,B)$. Assume that player II has an optimal strategy $\tau$ in $\mathcal{G}_{\mathrm{W}}(A,B)$. Then for any real $x$,
\begin{align*}
x \in A \iff \mu_{x , \tau} \bigl( \{ (x', y) \mid x' = x \text{ and } y\in B\} \bigr) = 1, 
\end{align*}
where $\mu_{x,\tau} = \mu_{\sigma_x , \tau}$ and $\sigma_x$ is the mixed strategy for player I concentrating on $x$ of measure one as in \S\,\ref{sec:Blackwell games}. 
Notice that the right hand side of the equivalence is $\undertilde{\mathbf{\Sigma}}^1_2$ in $B$: in fact, the right hand side holds if and only if for all $n\in \omega$, there is a tree $T$ on $2$ such that $\mu_{x,\tau} (\{ x \} \times [T]) > 1 - \frac{1}{n+1}$ and for all $y \in [T]$, $y \in B$. (The statement \lq $\mu_{x,\tau} (\{ x \} \times [T]) > 1 - \frac{1}{n+1}$' is Borel in $x$ and $T$.)  
If player I has an optimal strategy in $\mathcal{G}_{\mathrm{W}}(A,B)$, then one can prove that $B$ is $\undertilde{\mathbf{\Sigma}}^1_2$ in $A^{\mathrm{c}} (= 2^{\omega} \setminus A)$ in the same way and hence $B$ is $\undertilde{\mathbf{\Sigma}}^1_2$ in $A$. 
\end{proof}

\subsection{A weak version of Moschovakis' Coding Lemma}\label{sec:Moschovakis}

In this subsection, we discuss a weak version of Moschovakis' Coding Lemma which follows from $\BlAD$.


\begin{Thm}[Weak version of Moschovakis' Coding Lemma]\label{weak coding} 

Assume $\BlAD$. Let $<$ be a strict wellfounded relation on a set of reals $A$ with rank function $\rho \colon A \to \gamma$ onto and let $\Gamma$ be a Spector pointclass containing $<$ and closed under $\exists^{\mathbb{R}}$ and $\forall^{\mathbb{R}}$. Then for any subset $S $ of $\gamma$, there is a set of reals $C\in \undertilde{\mathbf{\Gamma}}$ such that $\rho [C] = S$, where $\rho [C] = \{ \rho (x) \mid x \in C \}$. 
\end{Thm}

By Theorem~\ref{inductive}, $\mathrm{IND}(<)$ satisfies the conditions for $\Gamma$ in Theorem~\ref{weak coding}. 

\begin{proof}
The argument is based on Moschovakis' original argument~\cite[7D.5]{new_Moschovakis}. 

Let $S$ be a subset of $\gamma$. We show that for all $\alpha \le \gamma$, there is a set of reals $C_{\alpha} \in \undertilde{\mathbf{\Gamma}}$ with $\rho [C_{\alpha}] = S \cap \alpha$ by induction on $\alpha$. 

It is trivial when $\alpha = 0$ and it is also easy when $\alpha$ is a successor ordinal because $\undertilde{\mathbf{\Gamma}}$ is a boldface pointclass. So assume that $\alpha $ is a limit ordinal and the above claim holds for each $\xi < \alpha$. We show that there is a $C \in \undertilde{\mathbf{\Gamma}}$ with $\rho [C] = S \cap \alpha$. 

Since $\Gamma$ is $\omega$-parametrized and closed under recursive substitutions, we have $\{G^n \subseteq \mathbb{R} \times \mathbb{R}^n \mid n \ge 1\}$ given in Lemma~\ref{good parametrization}.  Recall that for each real $a$, $G^1_a = \{ x \in \mathbb{R} \mid (a, x) \in G^1\}$. 
We say a set of reals $C$  {\it codes a subset $S'$ of $\gamma$} if $\rho [C] = S'$. 

Let us consider the following game $\mathcal{G}^{\alpha}_{\mathrm{M}}$: Player I and II choose $0$ or $1$ one by one and they produce reals $a$ and $b$ separately and respectively. Player II wins if either ($G^1_a$ does not code $S\cap \xi$ for any $\xi < \alpha$) or ($G^1_a$ codes $S \cap \xi$ for some $\xi < \alpha$ and $G^1_b$ codes $S\cap \eta$ for some $\eta< \alpha $ with $\eta > \xi$). By $\BlAD$, one of the players has an optimal strategy in this game.

\smallskip

\noindent \textbf{Case 1}: Player I has an optimal strategy $\sigma$ in $\mathcal{G}^{\alpha}_{\mathrm{M}}$. 

\smallskip

Recall from \S\,\ref{sec:Blackwell games} that for a real $b$, $\tau_b$ is the mixed strategy for player II concentrating on $b$ of measure one. 
Since $\sigma$ is optimal for player I, for each real $b$, for $\mu_{\sigma, \tau_b}$-measure one many reals $a$, $G^1_a$ codes $S\cap \xi$ for some $\xi< \alpha$. 

Fix a real $b'$. 
By Corollary~\ref{cor:BlAD-measure}, every set of reals is $\mu_{\sigma, \tau_{b'}}$-measurable and the ideal of $\mu_{\sigma, \tau_{b'}}$-null sets is closed under well-ordered unions. Hence there is a unique $\xi_{b'} < \alpha$ such that for $\mu_{\sigma , \tau_{b'}}$-positive measure many reals $a$, $G^1_a$ codes $S \cap \xi_{b'}$ and the set of reals $a$ such that $G^1_a$ codes $S\cap \xi$ for some $\xi < \xi_{b'}$ is $\mu_{\sigma, \tau_{b'}}$-measure zero. 

Let $C = \bigcup \{ G^1_a \mid \text{$G^1_a$ codes $S\cap \xi_b$ for some $b$} \}$. Then for all reals $x$, $x$ is in $C$ if and only if there is a real $b$ such that for $\mu_{\sigma, \tau_b}$-positive measure many reals $a$, they code the same subset $S'$ of $\gamma$, and no proper subsets of $S'$ can be coded by $\mu_{\sigma, \tau_b}$-positive measure many reals, and $x \in G^1_a$ for some real $a$ such that $G^1_a$ codes $S'$. Since $\Gamma$ is closed under $\exists^{\mathbb{R}} $ and $\forall^{\mathbb{R}}$, $C$ is in $\undertilde{\mathbf{\Gamma}}$. By induction hypothesis, for any $\xi< \alpha$, there is a real $b$ such that $G^1_b$ codes $S \cap \xi$. Since $\sigma$ is optimal for player I, $C$ codes $S \cap \alpha$, as desired. 

\smallskip

\noindent \textbf{Case 2}: Player II has an optimal strategy $\tau$ in $\mathcal{G}^{\alpha}_{\mathrm{M}}$. 

\smallskip

Let $ (a, x) \mapsto \{a\} (x)$ be the partial function from $\mathbb{R} \times \mathbb{R}$ to $\mathbb{R}$ which is universal for all the partial functions from $\mathbb{R}$ to itself that are $\undertilde{\mathbf{\Gamma}}$-recursive on their domain as in Lemma~\ref{lem:Kleene}. For reals $a$ and $w$, define a set of reals $A_{a, w}$ as follows: a real $x$ is in $A_{a,w}$ if there exists $v< w$ such that $\{ a\}(v)$ is defined and $\bigl(\{a\}(v), x\bigr) \in G^1$. It is easy to see that the set $\{ (a,w, x) \mid x \in A_{a,w}\}$ is in $\Gamma$. 
Hence by Lemma~\ref{good parametrization}, there is a $\Gamma$-recursive function $\pi \colon \mathbb{R} \times \mathbb{R} \to \mathbb{R}$ such that $A_{a,w} = G^1_{\pi(a,w)}$ for each $a$ and $w$. 

For all reals $a$ and $w$, define a set of reals $C_{a,w}$ as follows: A real $x$ is in $C_{a,w}$ if for $\mu_{\sigma_{\pi(a,w)}, \tau}$-positive measure many $b$, they code the same subset $S'$ of $\gamma$, no proper subsets of $S'$ can be coded by $\mu_{\sigma_{\pi(a,w)}, \tau}$-positive measure many reals, and $x$ is in $G^1_b$ for some real $b$ such that $G^1_b$ codes $S'$. It is easy to see that the set $\{ (a,w,x) \mid x \in C_{a,w}\}$ is in $\undertilde{\mathbf{\Gamma}}$. Hence by Lemma~\ref{good parametrization}, there is a $\undertilde{\mathbf{\Gamma}}$-recursive function $\pi' \colon \mathbb{R} \times \mathbb{R} \to \mathbb{R}$ such that $C_{a,w} = G^1_{\pi'(a,w)}$ for each $a$ and $w$. 

Since the function $(a,w) \mapsto \pi' (a,w) $ is $\undertilde{\mathbf{\Gamma}}$-recursive and total, by Kleene's Recursion Theorem~\ref{recursion theorem}, we can find a fixed $a^*$ such that for all $w$, $\{a^* \} (w) = \pi' (a^*, w)$. Let $g(w) = \{ a^* \} (w)$. 
\begin{Claim}\label{computation bound}
For each $w\in A$ with $\rho (w) < \alpha$, there is some $\eta (w) < \alpha$ with $\rho (w) < \eta (w)$ such that $G^1_{g(w)}$ codes $S \cap \eta (w)$. 
\end{Claim} 

\begin{proof}[Proof of Claim~\ref{computation bound}]
We show the claim by induction on $w$. 
Fix a $w \in A$ with $\rho (w) < \alpha$. Suppose that for all $v < w$, the conclusion of the claim holds for $v$. We will verify the conclusion of the claim for $w$, i.e., there is some $\eta (w) < \alpha$ with $\rho (w) < \eta (w)$ such that $G^1_{g(w)}$ codes $S \cap \eta (w)$. 

We first argue that $A_{a^*,w} = \bigcup_{v < w} G^1_{g(v)}$.  By the definition of $A_{a^*, w}$, for all reals $x$, $x$ is in $A_{a^*, w}$ if and only if for some $v < w$, $\{a^*\} (v)$ is defined and $\bigl( \{a^*\} (v) , x\bigr) \in G^1$, i.e., $x$ is in $G^1_{\{a^*\}(v)}$. Since the function $(a,v') \mapsto \pi' (a,v')$ is total and $\{a^*\} (v') = \pi' (a^* ,v' )$ for all reals $v'$, $\{a^*\} (v)$ is defined for all $v < w$. Also, by the definition of $g$, $ \{ a^* \} (v) = g(v)$ for all $v < w$. Hence, for all reals $x$, $x$ is in $A_{a^*, w}$ if and only if for some $v <w$, $x$ is in $G^1_{\{ a^* \} (v)} = G^1_{g(v)}$. Therefore, $A_{a^*,w} = \bigcup_{v < w} G^1_{g(v)}$, as desired. 

By induction hypothesis, for all $v <w$, there is some $\eta (v) < \alpha$ with $\rho (v) < \eta (v)$ such that $G^1_{g(v)}$ codes $S \cap \eta (v)$. Since $A_{a^*,w} = \bigcup_{v < w} G^1_{g(v)}$, letting $\xi = \mathrm{sup} \{ \eta(v) \mid v < w\}$, $A_{a^*, w}$ codes $S \cap \xi$. Also, since $\eta (v) > \rho (v)$ for all $v < w$, $\xi = \mathrm{sup} \{ \eta(v) \mid v < w\} \ge \mathrm{sup} \{ \rho (v) + 1 \mid v < w \} = \rho (w)$. 

We next argue that $C_{a^* , w}$ codes $S \cap \eta$ for some $\eta > \xi$. By the definition of $C_{a^*, w}$, for some subset $S'$ of $\gamma$, a real $x$ is in $C_{a^*, w}$ if and only if for $\mu_{\sigma_{\pi (a^*, w)}, \tau}$-positive measure many reals $b$, $G^1_b$ codes $S'$, no proper subsets of $S'$ can be coded by $\mu_{\sigma_{\pi (a^*, w)}, \tau}$-positive measure many reals, and $x$ is in $G^1_b$ for some $b$ such that $G^1_b$ codes $S'$. Hence $C_{a^*, w}$ also codes this $S'$. By the definition of $\pi$, $A_{a^*, w} = G^1_{\pi (a^*, w)}$. Since $A_{a^*, w} = G^1_{\pi (a^*,w)}$ codes $S \cap \xi$ and $\tau$ is an optimal strategy for player II, the set $S'$ must be of the form $S \cap \eta$ for some $\eta > \xi$. Hence $C_{a^*,w}$ codes $S \cap \eta$ for some $\eta > \xi$, as desired. 

Set $\eta (w) = \eta$. Then since $\eta > \xi \ge \rho (w)$, $\eta (w) > \rho (w)$. 
By the definition of $g$, $g(w) = \{ a^*\} (w) = \pi' (a^*, w)$. So $G^1_{g(w)} = G^1_{\pi' (a^*, w)} = C_{a^*,w}$, where the second equality follows from the property of $\pi'$. Since $C_{a^* , w}$ codes $S \cap \eta = S \cap \eta (w)$, $G^1_{g(w)}$ codes $S \cap \eta (w)$, as desired. 
\renewcommand{\qedsymbol}{$\square \ (\text{Claim~\ref{computation bound}})$}
\end{proof}

Let $C = \bigcup_{w \in A, \rho(w) < \alpha} G^1_{g(w)}$. Then by Claim~\ref{computation bound}, $C$ codes $S\cap \alpha$ (i.e., $\rho[C] = S \cap \alpha$) and $C$ is in $\undertilde{\mathbf{\Gamma}}$, as desired.

This completes the proof of Theorem~\ref{weak coding}. 
\end{proof}

\section{$\BlADR$ and a fine normal measure on $\wp_{\omega_1} (\mathbb{R})$}\label{sec:Rsharp}

In this section, we give a simple proof of Theorem~\ref{R sharp}, i.e., $\BlADR$ implies that there is a fine normal measure on $\wp_{\omega_1} (\R)$. 

Let us first see what is a fine normal measure. 
Let $X$ be a set and $\kappa$ be an uncountable cardinal. As usual, we denote by $\wp_\kappa(X)$ the set of all subsets of $X$ with cardinality less 
than $\kappa$, i.e., subsets $a$ of $X$ such that there are an $\alpha < \kappa$ and a surjection from $\alpha $ to $a$.
Let $U$ be a family of subsets of $\wp_{\kappa}(X)$. We say that $U$ is {\em $<$$\kappa$-complete} if $U$ is closed under intersections with $<$$\kappa$-many 
elements; we say it is {\em fine} if for any $x\in X$, $\{ a \in \wp_{\kappa}(X)  \mid x \in a\} \in U$; 
we say that $U$ is {\em normal} if for any family $\{ A_x \in U \mid x \in X\}$, the diagonal intersection $\triangle_{x\in X} A_x$ is in $U$, where 
$\triangle_{x\in X} A_x = \{ a \in \wp_{\kappa}(X) \mid (\forall x \in a ) \ a \in A_x \}$. We say that $U$ is a {\em fine measure} if it is a fine 
$<$$\kappa$-complete ultrafilter, and we say that it is a {\em fine normal measure} if it is a fine and normal $<$$\kappa$-complete ultrafilter. 

\begin{proof}[Proof of Theorem~\ref{R sharp}]

The following lemma is the key point: A subset $A $ of $\mathbb{R}^{\omega}$ is {\it range-invariant} if for any $\vec{x}$ and $\vec{y}$ in $\mathbb{R}^{\omega}$ with $\mathrm{ran}(\vec{x}) = \mathrm{ran}(\vec{y})$, $\vec{x} \in A$ if and only if $\vec{y} \in A$. 
\begin{Lem}\label{range-invariant}
Assume $\BlADR$. Then every range-invariant subset of $\mathbb{R}^{\omega}$ is determined. 
\end{Lem}

\begin{proof}[Proof of Lemma~\ref{range-invariant}]
Let $A$ be a range-invariant subset of $\mathbb{R}^{\omega}$. We show that if there is an optimal strategy for player I in $A$, then so is a winning strategy for player I in $A$. The case for player II is similar and we will skip it. 

Let us first introduce some terminology. Given a function $f \colon \mathbb{R}^{<\omega} \to \mathbb{R}$, a countable set of reals $a$ is {\em closed under $f$} if for any finite sequence $s$ of elements in $a$, $f(s)$ is in $a$. For a strategy $\sigma \colon \mathbb{R}^{\text{Even}} \to \mathbb{R}$ for player I, where $\mathbb{R}^{\text{Even} }$ is the set of all finite sequences of reals with even length, a countable set of reals $a$ is {\em closed under $\sigma$} if for any finite sequence $s$ of elements in $a$ of even length, $\sigma (s) $ is in $a$. For a function $F \colon \mathbb{R}^{<\omega} \to \wp_{\omega_1} (\mathbb{R})$, a countable set of reals $a$ is {\em closed under $F$} if for any finite sequence $s$ of elements in $a$, $F (s) $ is a subset of $a$. 

The following two claims are basic:
\begin{Claim}\label{1}
There is a winning strategy for player I in $A$ if and only if there is a function $f \colon \mathbb{R}^{<\omega} \to \mathbb{R}$ such that if $a $ is a countable set of reals and closed under $f$, then any enumeration of $a$ belongs to $A$.
\end{Claim}
\begin{proof}[Proof of Claim~\ref{1}]

We first show the direction from left to right. Given a winning strategy $\sigma$ for player I in $A$, let $f$ be such that if $a$ is closed under $f$, then $a$ is closed under $\sigma$. (Since $\sigma$ is a function from $\mathbb{R}^{\text{Even}}$ to $\mathbb{R}$, any function from $\mathbb{R}^{<\omega}$ to $\mathbb{R}$ extending $\sigma$ will be a desired $f$.) We see this $f$ works for our purpose. Let $a$ be a countable set of reals closed under $f$. Then since $a$ is closed under $\sigma$ and countable, there is a run $\vec{x}$ of the game following $\sigma$ such that its range is equal to $a$. Since $\sigma$ is winning for player I, $\vec{x}$ is in $A$ and by the range-invariance of $A$, any enumeration of $a$ is also in $A$. 

We now show the direction from right to left. Given such an $f$, we can arrange a strategy $\sigma$ for player I such that if $\vec{x}$ is a run of the game following $\sigma$, then the range of $\vec{x}$ is closed under $f$: Given a finite sequence of reals $(a_0, \cdots , a_{2n-1})$, consider the set of all finite sequences $s$ from elements of $\{ a_0, \cdots a_{2n-1}\}$ and all the values $f(s)$ from this set. What we should arrange is to choose $\sigma (a_0, \cdots , a_{2n-1})$ in such a way that the range of any run of the game via $\sigma$ will cover all such values $f(s)$ when $(a_0, \cdots , a_{2n-1})$ is a finite initial segment of the run for any $n$ in $\omega$ moves. But this is possible by a standard book-keeping argument. By the property of $f$, this implies that $\vec{x}$ is in $A$ and hence $\sigma$ is winning for player I. 
\renewcommand{\qedsymbol}{$\square$ (Claim~\ref{1})}
\end{proof}

\begin{Claim}\label{2}
There is a function $f \colon \mathbb{R}^{<\omega} \to \mathbb{R}$ such that if $a $ is a countable set of reals and closed under $f$, then any enumeration of $a$ belongs to $A$ if and only if there is a function $F \colon \mathbb{R}^{<\omega} \to \wp_{\omega_1} (\mathbb{R})$ such that if $a $ is a countable set of reals and closed under $F$, then any enumeration of $a$ belongs to $A$.

\end{Claim}

\begin{proof}[Proof of Claim~\ref{2}]

We first show the direction from left to right: Given such an $f$, let $F (s ) = \{ f(s) \}$. Then it is easy to check that this $F$ works.

We show the direction from right to left: Given such an $F$, it suffices to show that there is an $f $ such that if $a$ is closed under $f$ then $a$ is also closed under $F$. Fix a bijection $\pi \colon \mathbb{R} \to \mathbb{R}^{\omega}$. Let $g \colon \mathbb{R}^{<\omega} \to \mathbb{R}$ be such that $ \text{ran} \big( \pi (g (s) )\bigl) \supseteq F(s)$ for each $s$ (this is possible because every relation on the reals can be uniformized by a function by Theorem~\ref{uniformization}). Let $h \colon  \mathbb{R}^{<\omega} \to \mathbb{R}$ be such that $h (s) = \pi \bigl( s(0) \bigr) (\text{lh}(s)-1 )$, where $\text{lh}(s)$ is the length of $s$ when $s \neq \emptyset$, if $s = \emptyset$ let $h(s)$ be an arbitrary real. 

We verify that if $a$ is closed under $g$ and $h$, then so is under $F$: Fix a finite sequence $s$ of reals in $a$. We have to show that each $x$ in $F(s)$ is in $a$. 
Consider $g(s)$. By the closure under $g$, $g(s)$ is in $a$. By the choice of $g$, we know that $\text{ran}(\pi(g(s))) \supseteq F(s)$, so it is enough to show that $x$ is in $a$ for any $x$ in $\text{ran} (\pi(g(s))$. Suppose $x$ is the $n$th bit of $\pi(g(s))$. Consider the finite sequence $t =\bigl(g(s),...,g(s) \bigr) $ of length $n+1$. Then $h(t) = \pi(t(0))(\text{lh}(t)-1) = \pi(g(s))(n) =x$. 
But $g(s)$ is in $a$ and $a$ was closed under $h$, so $x $ is in $a$. 

Now it is easy to construct an $f$ such that if $a$ is closed under $f$, then so is under $g$ and $h$. 
\renewcommand{\qedsymbol}{$\square$ (Claim~\ref{2})}
\end{proof}

By Claim~\ref{1} and Claim~\ref{2} above, it suffices to show that there is a function $F \colon \mathbb{R}^{<\omega} \to \wp_{\omega_1} (\mathbb{R})$ such that if $a $ is a countable set of reals and closed under $F$, then any enumeration of $a$ belongs to $A$. 

Let $\sigma$ be an optimal strategy for player I in $A$. Let $F $ be as follows:\begin{equation*}
F(s) = \begin{cases}
\emptyset	&	\text{ if $\text{lh}(s)$ is odd,}\\
\{ y \in \mathbb{R} \mid \sigma(s) (y) \neq 0\}	&	\text{ if $\text{lh}(s)$ is even.}
\end{cases}
\end{equation*}

Then $F$ is as desired: If $a$ is closed under $F$, then enumerate $a$ to be $( a_n \mid n \in \omega )$ and let player I follow $\sigma$ and let player II play the Dirac measure concentrating on $a_n$ at her $n$th move. Then the probability of the set $\{ \vec{x} \in \mathbb{R}^{\omega} \mid \text{ran}(\vec{x}) = a \}$ is $1$. Since $\sigma $ is optimal for player I in $A$, there is an $\vec{x}$ such that the range of $\vec{x}$ is $a$ and $x$ is in $A$. But by the range-invariance of $A$, any enumeration of $a$ belongs to $A$. 

This completes the proof of Lemma~\ref{range-invariant}.
\end{proof}

We follow Solovay's original arguments for his theorem~\cite[Lemma~3.1]{Solovay_AD_R} that $\ADR$ implies there is a fine normal measure on $\wp_{\omega_1} (\mathbb{R})$. We define a family $U\subseteq \wp \bigl(\wp_{\omega_1}(\mathbb{R})\bigr)$
as follows: Fix an $A\subseteq\wp_{\omega_1}(\mathbb{R})$ and consider the following game $\mathcal{G}_{\mathrm{S}} (A)$: Players alternately play reals; 
say that they produce an infinite sequence $\vec x = (x_i\mid i\in\omega)$. Then player II wins the game $\mathcal{G}_{\mathrm{S}} (A)$ if $\mathrm{ran}(\vec{x}) \in A$, otherwise player I wins. Since the payoff set of this game is range-invariant, 
by Lemma~\ref{range-invariant}, it is determined. 

We say that $A\in U$ if and only if player II has a winning strategy in the game $\mathcal{G}_{\mathrm{S}} (A)$. We shall show that $U$ is a fine normal measure under the assumption of $\BlADR$, thus finishing the proof of Theorem~\ref{R sharp}.

A few properties of $U$ are obvious: For instance, we see readily that $\emptyset \notin U$ and that $\wp_{\omega_1}(\mathbb{R}) \in U$, as well
as the fact that $U$ is closed under taking supersets. In order to see that $U$ is a fine family, fix a real $x$, and let player II play $x$ in her first move: This is a winning strategy for player II in the game $\mathcal{G}_{\mathrm{S}} \bigl({\{a\mid x\in a\}}\bigr)$.

We next show that for any set $A \subseteq \wp_{\omega_1}(\mathbb{R})$, either $A$ or the complement of $A$ is in $U$. Given any such set $A$, suppose $A$ is not in $U$. We show that the complement of $A$ is in $U$. Since the game $\mathcal{G}_{\mathrm{S}} (A)$ is determined, by the assumption, there is a winning strategy $\sigma$ for I in $\mathcal{G}_{\mathrm{S}} (A)$. Setting $\tau (s) = \sigma \bigl( s {\upharpoonright} (\mathrm{lh}(s)-1) \bigr)$ for $s\in \mathbb{R}^{\mathrm{Odd}}$, it is easy to see that $\tau$ is a winning strategy for player II in the game $\mathcal{G}_{\mathrm{S}} (A^{\mathrm{c}})$, where $A^{\mathrm{c}}$ is the complement of $A$.

We show that $U$ is closed under finite intersections. Let $A_1$ and $A_2$ be in $U$. Since the payoff sets in the games $\mathcal{G}_{\mathrm{S}} (A_1)$ and $\mathcal{G}_{\mathrm{S}} (A_2)$ are range-invariant, by Claim~\ref{1}, there are functions $f_1 \colon \mathbb{R}^{<\omega} \to \mathbb{R}$ and $f_2 \colon \mathbb{R}^{<\omega} \to \mathbb{R}$ such that if $a$ is closed under $f_i$, then $a$ is in $A_i$ for $i=1,2$. Then it is easy to find an $f \colon \mathbb{R}^{<\omega} \to \mathbb{R}$ such that if $a$ is closed under $f$, then $a$ is closed under both $f_1$ and $f_2$. By Claim~\ref{1} again, this $f$ witnesses the existence of a winning strategy for player II in the game $\mathcal{G}_{\mathrm{S}} (A_1\cap A_2)$. 

We have shown that $U$ is an ultrafilter on subsets of $\wp_{\omega_1} (\mathbb{R})$. We show the $<$$\omega_1$-completeness of $U$ as follows: By Theorem~\ref{Blackwell Lebesgue}, every set of reals is Lebesgue measurable assuming $\BlAD$. If there is a non-principal ultrafilter on $\omega$, then there is a set of reals which is not Lebesgue measurable. Hence there is no non-principal ultrafilter on $\omega$, which implies that any ultrafilter on any set is $<$$\omega_1$-complete. In particular, $U$ is $<$$\omega_1$-complete. 

The last to show is that $U$ is normal. Let $\{A_y \mid y \in \mathbb{R}\}$ be a family of sets in $U$. We show that $\triangle_{y\in \mathbb{R} } A_y$ is in $U$. Consider the following game $\tilde{\mathcal{G}}_{\mathrm{S}}$: Player I moves $y$, then player II passes. After that, they play the game $\mathcal{G}_{\mathrm{S}} (A_y)$. This is Blackwell determined and player II has an optimal strategy $\tau$ since each $A_y$ is in $U$. Let $F\colon \mathbb{R}^{<\omega} \to \wp_{\omega_1}(\mathbb{R})$ be as follows: 
\begin{equation*}
F(s) = \begin{cases}
\emptyset	&	\text{ if $\text{lh}(s)$ is even,}\\
\{ z \in \mathbb{R} \mid \tau(s) (z) \neq 0\}	&	\text{ if $\text{lh}(s)$ is odd.}
\end{cases}
\end{equation*}
We claim that if $a$ is closed under $F$, then $a$ is in $\triangle_{y\in \mathbb{R} } A_y$. Then, by Claim~\ref{1} and Claim~\ref{2}, $F$ will witness the existence of a winning strategy for player II in the game $\mathcal{G}_{\mathrm{S}} \bigl(\triangle_{y\in \mathbb{R} } A_y\bigr)$, and we will have proved that $\triangle_{y \in \mathbb{R}} A_y  \in U$. 

Suppose $a$ is closed under $F$. We show that $a \in A_y$ for each $y \in a$. Fix a $y$ in $a$ and enumerate $a$ to be $( x_n \mid n \in \omega )$. In the game $\tilde{\mathcal{G}}_{\mathrm{S}}$, let player I first move $y$ and then they play the game $\mathcal{G}_{\mathrm{S}} (A_y)$. Let player II follow $\tau$ and player I play the Dirac measure concentrating on $x_n$ at the $n$th move. Then the probability of the set $\{ \vec{x} \in \mathbb{R}^{\omega} \mid \text{the first element of $\vec{x}$ is $y$ and } \text{ran}(\vec{x}) = a \}$ is $1$ and since $\tau$ is optimal for player II in the game $\tilde{\mathcal{G}}_{\mathrm{S}}$, there is an $\vec{x}$ whose first element is $y$ such that the range of $\vec{x}$ is $a$ and $\vec{x}$ is a winning run for player II in $\tilde{\mathcal{G}}_{\mathrm{S}}$, hence $a$ is in $A_y$. 

This completes the proof of Theorem~\ref{R sharp}.
\end{proof}

\section{$\BlADR$ and the Baire property}\label{sec:infBorel}

In this section, we show that $\BlADR$ implies that every set of reals has the Baire property. 
By Theorem~\ref{uniformization}, $\BlADR$ implies $\DC_\R$. For the rest of the sections in this paper, we freely use $\DC_\R$ when we assume $\BlADR$ and we fix a fine normal measure $U$ on $\wp_{\omega_1} (\mathbb{R})$, which exists by Theorem~\ref{R sharp}.

We first introduce the Blackwell meager ideal as an analogue of the meager ideal. A set of reals $A$ is {\it Blackwell meager} if player II has an optimal strategy in the Banach-Mazur game $\mathcal{G}_{\mathrm{BM}}(A)$. Let $I_{\mathrm{Bm}}$ denote the set of all Blackwell meager sets of reals. 
\begin{Lem}\label{Blackwell meager}
Assume $\BlAD$. Then every meager set is in $I_{\mathrm{Bm}}$, for each finite binary sequence $s$, $[s] \notin I_{\mathrm{Bm}}$, and $I_{\mathrm{Bm}}$ is a $\sigma$-ideal. Moreover, every set of reals is measurable via $I_{\mathrm{Bm}}$, i.e., for any set of reals $A$ and finite binary sequence $s$, there is a finite binary sequence $t$ extending $s$ such that either $[t] \cap A$ is in $I_{\mathrm{Bm}}$, or $[t] \setminus A$ is in $I_{\mathrm{Bm}}$. 
\end{Lem}

\begin{proof}
If a set of reals $A$ is meager, then player II has a winning strategy in the Banach-Mazur game $\mathcal{G}_{\mathrm{BM}}(A)$ and in particular player II has an optimal strategy in $\mathcal{G}_{\mathrm{BM}}(A)$ by Theorem~\ref{det implies Blackwell det}. Hence $A$ is Blackwell meager. 

It is easy to see that $[s] \notin I_{\mathrm{Bm}}$ for each finite binary sequence $s$ by letting player I first play the Dirac measure concentrating on $s$ in the game $\mathcal{G}_{\mathrm{BM}} ([s])$. 

We show that $I_{\mathrm{Bm}}$ is a $\sigma$-ideal. The closure of $I_{\mathrm{Bm}}$ under subsets is immediate. We prove that it is closed under countable unions. 

In order to prove this, we need to develop the appropriate \emph{transfer technique} (as discussed and applied in \cite{MR1939050}) for the present context. Let $\pi\subseteq\omega$ be an infinite and co-infinite set. We think of $\pi$ as the set of 
rounds in which player I moves. We identify $\pi$ with the increasing enumeration of its members, i.e., $\pi = \{\pi_i\mid i\in\omega\}$.
Similarly, we write $\overline\pi$ for the increasing enumeration of $\omega{\setminus}\pi$, i.e., $\omega{\setminus}\pi = \{\overline\pi_i\mid i\in\omega\}$.
For notational ease, we call $\pi$ a \textbf{I-coding} if no two consecutive numbers are in $\pi$ and $0 \in \pi$ (i.e., the first move is played by I). We call $\pi$ a \textbf{II-coding} if no two consecutive numbers are in $\omega{\setminus}\pi$ and $0\in \pi$.

Fix $A \subseteq 2^{\omega}$ and define two variants of $\mathcal{G}_{\mathrm{BM}} (A)$ with alternative orders of play as determined by $\pi$. If $\pi$ is a I-coding, 
the game $\mathcal{G}_{\mathrm{BM}}^{\pi,\mathrm{I}} (A)$ is played as follows:
\[
\begin{array}{cccccc}
	\text{I}	&	s_{\pi_0} = s_0				&			&		s_{\pi_1}			&		&	\ldots	\\
	\text{II}	&		&	s_{\pi_0+1},	\ldots	,s_{\pi_1 -1}		&		&	s_{\pi_1+1}, \ldots , s_{\pi_2 -1}	&	 \ldots	\\	
\end{array}
\]
If $\pi$ is a II-coding, then they play the game $\mathcal{G}_{\mathrm{BM}}^{\pi,\mathrm{II}} (A)$ as follows:
\[
\begin{array}{cccccc}
	\text{I}	&	s_0,	\ldots	,s_{\overline\pi_0 -1}	&			&	s_{\overline\pi_0+1}, \ldots , s_{\overline\pi_1 -1}	&			 &	\ldots	\\
	\text{II}	&					&	s_{\overline\pi_0}		&					&	s_{\overline\pi_1}		&	\ldots	\\	
\end{array}
\]
In both cases, player II wins the game if $s_0^{\frown} s_1^{\frown} \ldots ^{\frown} s_n ^{\frown} \ldots \notin A$.
Obviously, we have $$\mathcal{G}_{\mathrm{BM}} (A)=\mathcal{G}_{\mathrm{BM}}^{\mathrm{Even},\mathrm{II}} (A)$$ where $\mathrm{Even}$ is the set of all even numbers.

\begin{Sublem}\label{key lemma}
Let $A \subseteq 2^{\omega}$ and $\pi$ be a I-coding.
Then there is a translation $\sigma\mapsto\sigma_\pi$ of mixed strategies for player I such that if $\sigma$ is an optimal strategy for player I in
$\mathcal{G}_{\mathrm{BM}} (A)$, then $\sigma_{\pi}$ is an optimal strategy for player I in $\mathcal{G}_{\mathrm{BM}}^{\pi,\mathrm{I}} (A)$.

Similarly, if $\pi$ is a II-coding, there is a translation $\tau\mapsto\tau_\pi$ of mixed strategies for player II such that
if $\tau$ is an optimal strategy for player II in $\mathcal{G}_{\mathrm{BM}} (A)$, then $\tau_\pi$ is an optimal strategy  for player II in $\mathcal{G}_{\mathrm{BM}}^{\pi,\mathrm{II}} (A)$.
\end{Sublem}

\begin{proof}[Proof of Sublemma~\ref{key lemma}]

We prove only the lemma for the games $\mathcal{G}_{\mathrm{BM}}^{\pi,\mathrm{I}} (A)$, the other proof being similar.
If $\vec s = ( s_i\mid i\in\omega )$ is an
infinite sequence of finite binary sequences, we define
$$b^{\vec s}_i =
s_{\pi_i +1 } ^{\frown} \ldots ^{\frown} s_{\pi_{i+1} -1}.$$
Note that in order to compute $b^{\vec s}_i$, we only need the first $\pi_{i+1}$ bits of $\vec s$.
The idea is that now the $\mathcal{G}_{\mathrm{BM}} (A)$-run
\begin{equation}\tag{$*$}\label{eq:raveled}
\begin{array}{cccccccc}
	\text{I}	&	s_{\pi_0} &			      &	s_{\pi_1}  &		         & s_{\pi_2} &              & \ldots\\
	\text{II}	&		      &	b^{\vec s}_0  &		       &	b^{\vec s}_1 &           & b^{\vec s}_2 & \ldots
\end{array}
\end{equation}
yields the same output in terms of the concatenation of all played finite sets as the run $\vec s$ in the game $\mathcal{G}_{\mathrm{BM}}^{\pi,\mathrm{I}} (A)$. We can define
a map $\pi^*$ on infinite sequences of finite binary sequences by
$$(\pi^*(\vec s))_i = \left\{
\begin{array}{ll}
s_{\pi_k} & \mbox{if $i = 2k$,}\\
b^{\vec s}_k & \mbox{if $i=2k+1$,}
\end{array}\right.$$
and see that $s_0^{\frown} s_1^{\frown} \ldots  = (\pi^*(\vec s))_0^{\frown} (\pi^*(\vec s))_1^{\frown} \ldots$.

Now, given a mixed strategy $\sigma$ for player I in $\mathcal{G}_{\mathrm{BM}} (A)$ and a run $\vec s$ of the game $\mathcal{G}_{\mathrm{BM}} (A)$, we define $\sigma_{\pi}$ via $\pi^*$ as follows:
\begin{align*}
\sigma_{\pi} (s_0, \ldots , s_{\pi_m -1} ) =   \sigma (s_{\pi_0}, b^{\vec s}_0, \ldots , s_{\pi_i}, b^{\vec s}_i,  \ldots , s_{\pi_{m-1}} ,
b^{\vec s}_{m-1}).
\end{align*}
Assume that $\sigma$ is an optimal strategy for player I in $\mathcal{G}_{\mathrm{BM}} (A)$ and fix an arbitrary mixed strategy $\tau$ in the game $\mathcal{G}_{\mathrm{BM}}^{\pi,\mathrm{I}} (A)$.
We show that the payoff set for $A$ in $\mathcal{G}_{\mathrm{BM}}^{\pi,\mathrm{I}} (A)$ is $\mu_{\sigma_{\pi}, \tau}$-measurable and $\mu_{\sigma_{\pi}, \tau} (A) = 1$. In order to do so, we
construct a mixed strategy $\tau_{\pi^{-1}}$ for player II in $\mathcal{G}_{\mathrm{BM}} (A)$ so that the game played by $\sigma_{\pi}$ and $\tau$ is essentially the same as the
game played by $\sigma$ and $\tau_{\pi^{-1}}$.

Given a sequence $\vec b$ of moves in $\mathcal{G}_{\mathrm{BM}} (A)$, we need to unravel it into a sequence of moves in $\mathcal{G}_{\mathrm{BM}}^{\pi,\mathrm{I}} (A)$ in an
inverse of the maps $\vec s\mapsto b^{\vec s}_i$ according to (\ref{eq:raveled}), i.e., $b_{2i+1} = b^{\vec s}_i$.
Thus, we define
\begin{eqnarray*}
A^{\vec b}_{2i+1}  & =  & \{\vec s\mid b^{\vec s}_i = b_{2i+1}\}\mbox{,}\\
A^{\vec b}_{\leq 2i+1}  & =  & \bigcap_{j\leq i} A^{\vec b}_{2j+1}.
\end{eqnarray*}
Note that only a finite fragment of $\vec s$ is needed to check whether $b^{\vec s}_i = b_{2i+1}$, and thus we think of $A^{\vec b}_{\leq 2i+1}$ as a 
set of $(\pi_{i+1}-(i+1))$-tuples of finite binary sequences. In the following, when we quantify over all ``$\vec s\in A^{\vec b}_{\leq 2i+1}$'', we think 
of this as a collection of finite strings of finite binary sequences. In order to pad the moves made in $\mathcal{G}_{\mathrm{BM}}^{\pi,\mathrm{I}} (A)$, we define the following
notation: For infinite sequences $\vec s$ and $\vec b$, we write
\begin{align*}
x^{\vec s,\vec b}_i  =  (b_{2i},s_{\pi_i+1},...,s_{\pi_{i+1}-1}).
\end{align*}
Compare (\ref{eq:raveled}) to see that if $\vec s$ corresponds to moves in $\mathcal{G}_{\mathrm{BM}}^{\pi,\mathrm{I}} (A)$ and $\vec b$ to the moves in $\mathcal{G}_{\mathrm{BM}} (A)$,
then these are exactly the finite sequences that player II will have to respond to in $\mathcal{G}_{\mathrm{BM}}^{\pi,\mathrm{I}} (A)$.
Moreover, for a given sequence $\vec z$ of finite binary sequences, we let
$$P_\tau(z_0,...,z_n) = \prod_{k\leq n,k\notin\pi} \tau(z_0,...,z_{k-1})(z_k).$$

Fix a sequence $\vec b$ of finite binary sequences with even length and define
$\tau_{\pi^{-1}}$ as follows:
\begin{eqnarray*}
\tau_{\pi^{-1}} (b_0) (b_1) & = & \sum_{\vec s\in A^{\vec b}_1} P_\tau(x^{\vec s,\vec b}_0)\mbox{, and}\\
\tau_{\pi^{-1}} (b_0, \ldots , b_{2m}) (b_{2m+1}) & = & \frac{\sum_{\vec s\in A^{\vec b}_{\leq 2m+1}} P_\tau (x^{\vec s,\vec
b}_0{}^\smallfrown\ldots^\smallfrown x^{\vec s,\vec b}_m)}{\prod_{i=1}^m \tau_{\pi^{-1}} (b_0, \ldots , b_{2i-2}) (b_{2i-1})}.
\end{eqnarray*}

Using the two operations $\sigma\mapsto\sigma_\pi$ and $\tau\mapsto\tau_{\pi^{-1}}$, since the
payoff set for $\mathcal{G}_{\mathrm{BM}} (A)$ is invariant under
$\pi^*$, it now suffices to prove for all basic open
sets $[t]$ induced by a finite sequence $t = (b_0,...,b_{\mathrm{lh}(t)-1})$ that
$\mu_{\sigma, \tau_{\pi^{-1}}} ([t]) = \mu_{\sigma_{\pi}, \tau} ( (\pi^{\ast})^{-1} ([t]))$.
We prove this by induction on the length of $t$, and have to consider three different cases:

\smallskip

\noindent \textbf{Case 1.} $\text{lh}(t) =0$.
This is immediate.

\smallskip

\noindent\textbf{Case 2.} $\text{lh}(t) = 2m+1$ with $m \ge 0$. By induction hypothesis, we have that
$X = \mu_{\sigma , \tau_{\pi^{-1}}} ([b_0, \ldots , b_{2m-1}]) = \mu_{\sigma_{\pi}, \tau} ( (\pi^{\ast})^{-1} ( [b_0, \ldots , b_{2m-1}]))$. Thus,
\begin{eqnarray*}
\mu_{\sigma , \tau_{\pi^{-1}}} ( [b_0, \ldots , b_{2m}])
&= &  X \cdot \sigma (b_0, \ldots , b_{2m-1}) (b_{2m})\\
& = & \mu_{\sigma_{\pi}, \tau} ( (\pi^{\ast})^{-1}( [b_0, \ldots , b_{2m}]) ).
\end{eqnarray*}

\smallskip

\noindent\textbf{Case 3.} $\text{lh}(t) = 2m+2$ with $m \ge 0$.
\begin{eqnarray*}
\mu_{\sigma , \tau_{\pi^{-1}}} (t)
& = & \prod_{i=0}^{m} \sigma (b_0, \ldots  , b_{2i-1}) (b_{2i}) \cdot \sum_{\vec s\in A^{\vec b}_{\leq 2m+1}} P_{\tau} (x^{\vec s,\vec
b}_0{}^\smallfrown\ldots^\smallfrown x^{\vec s,\vec b}_m) \\
& = & \ \mu_{\sigma_{\pi}, \tau} \bigl( (\pi^{\ast})^{-1} ( [b_0, \ldots , b_{2m+1}] ) \bigr).
\end{eqnarray*}
This calculation finishes the proof of Sublemma~\ref{key lemma}.
\end{proof}

We now show that $I_{\mathrm{Bm}}$ is closed under countable unions. Let $\{A_n \mid n \in \omega\}$ be a family of sets in $I_{\mathrm{Bm}}$. Using $\AC_{\omega}(\R)$, take an optimal strategy $\tau_n$ in the game $\mathcal{G}_{\mathrm{BM}}(A_n)$ for each $n$. We prove that $\bigcup_{n\in \omega} A_n$ is also in $I_{\mathrm{Bm}}$. 

Fix a bookkeeping bijection $\rho$ from $\omega\times\omega$ to $\omega$ such that $\rho(n, m) < \rho (n, m +1)$ and $\rho (n,0) \ge n$. We are 
playing infinitely many games in a diagram where the first coordinate is for the index of the game we are playing,
and the second coordinate is for the number of moves. Hence the pair $(n,m)$ stands for \lq \lq $m$th move in the $n$th game".
Define a II-coding $\pi_n = \omega{\setminus}\{2\rho(n,i)+1\mid i\in\omega\}$ corresponding to the following game diagram:
\[
\begin{array}{cccccc}
\text{I}	&	s_0, \ldots , s_{2 \rho (n,0)}	&			&	s_{2 \rho(n,0) +2}, \ldots , s_{2 \rho (n,1)}		&			 &	\ldots	\\
\text{II}	&			&	s_{2 \rho (n,0) +1}	&			&	s_{2 \rho (n, 1) + 1}	&	\ldots
\end{array}
\]
By Sublemma~\ref{key lemma}, we know that for each $n \in \omega$, the mixed strategy
$(\tau_n)_{\pi_n}$ is optimal for player II in the game $\mathcal{G}_{\mathrm{BM}}^{\pi_n,\mathrm{II}} (A_n)$.
Let $\tau$ be the following mixed strategy
$$\tau(s_0, \ldots , s_{2\rho(n,m)}) = (\tau_n)_{\pi_n} (s_0, \ldots , s_{2\rho(n,m)}).$$
The properties of $\rho$ make sure that this strategy is well-defined. We shall now prove that $\tau$ is an optimal strategy for player II in 
$\mathcal{G}_{\mathrm{BM}} \bigl(\bigcup_{n \in \omega} A_n \bigr)$.

Pick any mixed strategy $\sigma$ for player I in
$\mathcal{G}_{\mathrm{BM}} \bigl(\bigcup_{n \in \omega} A_n\bigr)$ and define strategies $\sigma_n$ for $\mathcal{G}_{\mathrm{BM}}^{\pi_n,\mathrm{II}} (A_n)$. Let $m = \rho(k,\ell)$, then
\begin{eqnarray*}
\sigma_n (s_0, \ldots , s_{2m -1}) &=&   \sigma (s_0, \ldots , s_{2m-1})\mbox{, and}\\
\sigma_n (s_0, \ldots , s_{2m}) &=&  (\tau_k)_{\pi_k} (s_0, \ldots , s_{2m}) \text{ (if $k \neq n$).}
\end{eqnarray*}
Note that for each $n \in \omega $, $\mu_{\sigma, \tau}=\mu_{\sigma_n, (\tau_n)_{\pi_n}}$. 

The payoff set (for player II) in $\mathcal{G}_{\mathrm{BM}} \bigl(\bigcup_{n \in \omega} A_n\bigr)$ is $A=\{\vec s \mid s_0^{\frown} s_1^{\frown} \ldots \notin
\bigcup_{n \in \omega} A_n\}$. We show that $\mu_{\sigma, \tau}(A)=1$. Since
$A = \bigcap_{n \in \omega} \left\{ \vec s\mid s_0^{\frown} s_1^{\frown} \ldots \notin A_{n}\right\}$,
it suffices to check that the sets
$B_n = \{ \vec s\mid s_0^{\frown} s_1^{\frown} \ldots\notin A_{n}\}$ has $\mu_{\sigma, \tau}$-measure~$1$. But
$$\mu_{\sigma, \tau} (B_n) = \mu_{\sigma_n , (\tau_n)_{\pi_n}} (B_n) = 1.$$
Thus we have shown that $I_{\mathrm{Bm}}$ is a $\sigma$-ideal. 

We finally show that every set of reals $A$ is measurable with respect to $I_{\mathrm{Bm}}$, i.e., for any finite binary sequence $s$, there is a finite binary sequence $t$ extending $s$ such that either $[t] \cap A$ is in $I_{\mathrm{Bm}}$, or $[t] \setminus A$ is in $I_{\mathrm{Bm}}$. Fix such $A$ and $s$. If $[s] \cap A$ is in $I_{\mathrm{Bm}}$, we are done. So suppose not. Then player II does not have an optimal strategy in the game $\mathcal{G}_{\mathrm{BM}}([s] \cap A)$. By $\BlAD$, there is an optimal strategy $\sigma$ for player I in the game $\mathcal{G}_{\mathrm{BM}}([s] \cap A)$. Let $t$ be any $s'$ with $\sigma (\emptyset) (s') \neq 0$. Then since $\sigma$ is optimal, $t$ extends $s$ and the strategy $\sigma$ easily gives us an optimal strategy for player II in the game $\mathcal{G}_{\mathrm{BM}}([t] \setminus A)$. Hence $[t] \setminus A$ is in~$I_{\mathrm{Bm}}$. 

This completes the proof of Lemma~\ref{Blackwell meager}. 
\end{proof}

Recall the definitions of Stone space $\text{St}(\mathbb{P})$ and $\mathbb{P}$-Baireness for a preorder $\mathbb{P}$ from \S\,\ref{subsec:P}. The base set of $\text{St}(\mathbb{P})$ was the set of all ultrafilters on $\mathbb{P}$ and without the Axiom of Choice, it might be empty if $\mathbb{P}$ is big. But in this section, we only consider preorders $\mathbb{P}$ which are elements of $H_{\omega_1}$ in $V$, i.e., the transitive closure of $\{ \mathbb{P} \}$ is countable in $V$. If $\mathbb{P}$ is an element of $H_{\omega_1}$, then by Lemma~\ref{Baireness_Baire}, $\text{St}(\mathbb{P})$ is isomorphic to $\text{St}(\mathbb{C})$ as Baire spaces where $\mathbb{C}$ is Cohen forcing, hence $\text{St}(\mathbb{C})$ is homeomorphic to the Cantor space $2^{\omega}$. 
%
%

Since every meager set is Blackwell meager as we have seen in Lemma~\ref{Blackwell meager}, if $\mathbb{P}$ is in $H_{\omega_1}$, then one can define Blackwell meagerness for subsets of $\mathrm{St}(\mathbb{P})$ via an isomorphism between the Baire space and $\mathrm{St}(\mathbb{P})$ as Baire spaces and identify them as structures of topological spaces together with Blackwell meager ideals. From now on, we will use this identification without any notice. 

We are now ready to prove the Baire property for every set of reals from~$\BlADR$. 
\begin{Thm}\label{Blackwell Baire property}
Assume $\BlADR$. Then every set of reals has the Baire property. 
\end{Thm}

\begin{proof}
Take any set of reals $A$. We show that $A$ has the Baire property. Let $\mathcal{A}^2_A = (\omega , 2^{\omega} , 0, 1, + , \times, \in , A)$ be the second-order arithmetic structure with $A$ as a unary predicate. Since any relation on the reals can be uniformized by a function by Theorem~\ref{uniformization}, we can construct a Skolem function $F$ for $\mathcal{A}^2_A$ and by a simple coding of finite sequences of reals and formulas via reals, we regard it as a function from the reals to themselves. Let $\Gamma_F$ be the graph of $F$, i.e., $\Gamma_F = \{ (x, s) \in \mathbb{R} \times 2^{<\omega} \mid F(x) \supseteq s \}$. The following are the key objects for the proof (they are called {\it term relations}): Recall from Lemma~\ref{name-Baire measurable} that for a $\mathbb{P}$-name $\tau$ for a real, $f_{\tau}$ is the Baire measurable function (which is continuous on a comeager set) corresponding to $\tau$. We write $A^{\mathrm{c}}$ for the complement of $A$ and the same for $(\Gamma_F)^{\mathrm{c}}$. 
\begin{align*}
\tau_A	= 	\{ (\mathbb{P}, p, \sigma) \in H_{\omega_1} \mid &\ \text{$\sigma$ is a $\mathbb{P}$-name for a real and }	\\
&\bigl(\forall^{\infty}_{\mathrm{B}} G \in \mathrm{St}(\mathbb{P}) \bigr) \ p\in G \implies f_{\sigma}(G) \in A \}, \\
\tau_{A^{\mathrm{c}}}	=  \{ (\mathbb{P}, p, \sigma) \in H_{\omega_1} \mid &\ \text{$\sigma$ is a $\mathbb{P}$-name for a real and }	\\
\displaybreak[0]
&\bigl(\forall^{\infty}_{\mathrm{B}} G \in \mathrm{St}(\mathbb{P}) \bigr) \ p\in G \implies f_{\sigma}(G) \in A^{\mathrm{c}} \}, \\
\tau_{\Gamma_F} =   \{ (\mathbb{P}, p, \sigma , s) \in H_{\omega_1} \mid &\  \text{$\sigma$ is a $\mathbb{P}$-name for a real and }	\\
&\bigl(\forall^{\infty}_{\mathrm{B}} G \in \mathrm{St}(\mathbb{P}) \bigr) \ p\in G \implies \bigl( f_{\sigma}(G), s\bigr) \in \Gamma_F\}, \\
\tau_{(\Gamma_F)^{\mathrm{c}}} =  \{ (\mathbb{P}, p, \sigma , s) \in H_{\omega_1} \mid &\  \text{$\sigma$ is a $\mathbb{P}$-name for a real and }	\\
&\bigl(\forall^{\infty}_{\mathrm{B}} G \in \mathrm{St}(\mathbb{P}) \bigr) \ p\in G \implies \bigl( f_{\sigma}(G), s\bigr) \in (\Gamma_F)^{\mathrm{c}}\}, 
\end{align*}
where $\bigl(\forall^{\infty}_{\mathrm{B}} G \in \mathrm{St}(\mathbb{P}) \bigr)$ means \lq\lq for all $G$ modulo a Blackwell meager set in $\mathrm{St}(\mathbb{P})$\ldots". 
Let $\vec{\tau} = \bigl(\tau_A, \tau_{A^{\mathrm{c}}}, \tau_{\Gamma_F}, \tau_{(\Gamma_F)^{\mathrm{c}}}\bigr)$ and set $M = \mathrm{HOD}^{\mathrm{L}[\vec{\tau}]}_{\vec{\tau}}$. For $G \in \mathrm{St}(\mathbb{P})$, let $A_G = \{ f_{\sigma}(G) \mid (\exists p \in G) \ (\mathbb{P}, p, \sigma) \in \tau_A \cap M\}$.  
Note that for any set $x$ in $M \cap H_{\omega_1}$, $\wp(x) \cap M$ is countable: Since $M$ is a transitive model of $\ZFC$, if $\wp(x) \cap M = \wp^M (x)$ was uncountable, then there would be an injection from $\omega_1$ to the reals. However, by Corollary~\ref{cor:BlAD-measure} and Proposition~\ref{easy_Blackwell}, $\BlADR$ implies there is no injection from $\omega_1$ to the reals. Hence for any $\mathbb{P} \in H_{\omega_1} \cap M$, the set of $\mathbb{P}$-generic filters over $M$ is comeager, in particular Blackwell comeager (i.e., its complement is Blackwell meager). Therefore, when we discuss statements starting from $\bigl( \forall^{\infty}_{\mathrm{B}} G \in \mathrm{St}(\mathbb{P})\bigr)$, we may assume that $G$ is $\mathbb{P}$-generic over $M$. 
\begin{Claim}\label{local Baire}
${}$ 

\begin{enumerate}
\item Let $\mathbb{P}$ be a preorder in $M$. Then $\bigl( \forall^{\infty}_{\mathrm{B}} G \in \mathrm{St}(\mathbb{P})\bigr) \ A_G = A \cap M[G] \in M[G]$ and $M[G]$ is closed under $F$. 

\item Let $\mathbb{P} = \mathrm{Coll}(\omega, 2^{\omega})^M$, where $\mathrm{Coll}(\omega, 2^{\omega})$ is the forcing collapsing $2^{\omega}$ into countable with finite conditions. Then $\bigl( \forall^{\infty}_{\mathrm{B}} G \in \mathrm{St}(\mathbb{P})\bigr)$ $A_G$ has the Baire property in $M[G]$. 
\end{enumerate}
\end{Claim}
\begin{proof}
We first show that $A_G = A \cap M[G]$ for Blackwell comeager many $G$. 
Since for any set $x$ in $M \cap H_{\omega_1}$, $\wp (x) \cap M$ is countable and $I_{\mathrm{Bm}}$ is a $\sigma$-ideal, it follows that for Blackwell comeager many $G$, $G$ is $\mathbb{P}$-generic over $M$ and if $(\mathbb{P}, p , \sigma) \in \tau_A \cap M$ (resp., $\tau_{A^{\mathrm{c}}} \cap M$) and $p\in G$, then $f_{\sigma}(G) = \sigma^G \in A$ (resp., $A^{\mathrm{c}}$). We show that $A_G = A \cap M[G]$ for any such $G$. 

Fix such a $G$. We first see that $A_G \subseteq A \cap M[G]$. Take any real $x$ in $A_G$. Then there are a $p\in G$ and a $\sigma$ such that $(\mathbb{P}, p , \sigma) \in \tau_A \cap M$ and $\sigma^G = x$. Then by the property of $G$, $x = \sigma^G = f_{\sigma}(G) \in A$, as desired. 
We show that $A \cap M[G] \subseteq A_G$. Let $x$ be a real in $M[G]$ which is not in $A_G$. We prove that $x$ is also not in $A$. Since $x$ is in $M[G]$, there is a $\mathbb{P}$-name $\sigma$ for a real in $M$ such that $\sigma^G = x$. Since every subset of $\mathrm{St} (\mathbb{P})$ is measurable with respect to $I_{\mathrm{Bm}}$ by Lemma~\ref{Blackwell meager}, the set $\{ p \in \mathbb{P} \mid \text{ either } (\mathbb{P}, p , \sigma) \in \tau_A\cap M \text{ or } (\mathbb{P}, p , \sigma) \in \tau_{A^{\mathrm{c}}}\cap M\}$ is dense and it is in $M$. Since $G$ is $\mathbb{P}$-generic over $M$, there is a $p\in G$ such that either $(\mathbb{P}, p, \sigma) \in \tau_A$ or $(\mathbb{P}, p , \sigma) \in \tau_{A^{\mathrm{c}}}$. But $(\mathbb{P}, p, \sigma) \in \tau_A$ cannot hold because it would imply $x = \sigma^G \in A_G$. Hence $(\mathbb{P}, p, \sigma) \in \tau_{A^{\mathrm{c}}}$ and $x = \sigma^G = f_{\sigma}(G) \in A^{\mathrm{c}}$ by the property of $G$, as desired. 

Let $\rho_A = \{ (\sigma, p ) \mid (\mathbb{P}, p, \sigma) \in \tau_A \cap M\}$. Since the comprehension axioms with $\tau_A$ as a unary predicate hold in $M$, $\rho_A$ is a $\mathbb{P}$-name for a set of reals in $M$ and $\rho_A^G = A_G\in M[G]$. Hence $A_G = A \cap M[G] \in M[G]$ for Blackwell comeager many $G$, as desired. 

Next, we show that $M[G]$ is closed under $F$ for Blackwell comeager many $G$. We prove this for any $G$ which is $\mathbb{P}$-generic over $M$ such that if $(\mathbb{P}, p, \sigma, s) \in \tau_{\Gamma_F}$ (resp., $\tau_{\Gamma_{F^{\mathrm{c}}}}$) and $p$ is in $G$, then $F(\sigma^G) \supseteq s$ (resp., $F(\sigma^G) \nsupseteq s$). 
Fix such a $G$ and let $x$ be a real in $M[G]$. We show that $F(x)$ is also in $M[G]$. Since $x$ is in $M[G]$, there is a $\mathbb{P}$-name $\sigma$ for a real in $M$ such that $\sigma^G = x$. Since every subset of $\mathrm{St}(\mathbb{P})$ is measurable with respect to $I_{\mathrm{Bm}}$, the function $G' \mapsto F\bigl(f_{\sigma}(G')\bigr)$ is continuous modulo a Blackwell meager set in $\mathrm{St}(\mathbb{P})$. Hence for any finite binary sequence $s$, the set of all $p\in \mathbb{P}$ such that either $\bigl( \forall^{\infty}_{\mathrm{B}} G' \in \mathrm{St}(\mathbb{P}) \bigr)\  p \in G' \implies  F \bigl(f_{\sigma}(G') \bigr) \supseteq s$, \text{ or } $\bigl( \forall^{\infty}_{\mathrm{B}} G' \in \mathrm{St}(\mathbb{P}) \bigr)\  p \in G' \implies  F \bigl(f_{\sigma}(G') \bigr) \nsupseteq s $, is dense and is in $M$, given that $M$ satisfies Comprehension axioms with predicates $\tau_{\Gamma_F}$ and $\tau_{(\Gamma_F)^{\mathrm{c}}}$. By the genericity and the property of $G$, for any $s$, there is a $p\in G$ such that $F(\sigma^G) \supseteq s$ if and only if $\bigl( \forall^{\infty}_{\mathrm{B}} G' \in \mathrm{St}(\mathbb{P}) \bigr)\  p \in G' \implies  F \bigl(f_{\sigma}(G') \bigr) \supseteq s$ if and only if $(\mathbb{P}, p, \sigma, s) \in \tau_{\Gamma_f} \cap M$. Hence $ F(x) = F(\sigma^G) = \bigcup \{ s \mid (\exists p \in G ) \ (\mathbb{P}, p, \sigma, s) \in \tau_{\Gamma_f} \cap M\}$, which is in $M[G]$, as desired. 

Finally, we show that $A_G$ has the Baire property in $M[G]$ for Blackwell comeager many $G$ when $\mathbb{P} = \mathrm{Coll}(\omega, 2^{\omega})^M$. Actually, we show that $A_G$ has the Baire property in $M[G]$ for any $\mathbb{P}$-generic $G$ over $M$. Let $G$ be any $\mathbb{P}$-generic filter over $M$ and $s$ be a finite binary sequence. We show that there is a $t$ extending $s$ such that either $[t] \cap A_G$ or $[t] \setminus A_G$ is meager in $M[G]$. Let $\dot{c}$ be a canonical name for a Cohen real. Since one can embed Cohen forcing into $\mathrm{Coll}(\omega, 2^{\omega})^M$ in a natural way in $M$, we may regard $\dot{c}$ as a $\mathbb{P}$-name for a Cohen real. Since $2^{\omega} $ in $M$ is countable in $M[G]$, the set of Cohen reals over $M$ is comeager in $M[G]$. Take any Cohen real $c$ over $M$ with $s \subseteq c$ in $M[G]$. We may assume $c $ is in $A_G$ (the case $c \notin A_G$ can be dealt with in the same way). Recall that $\rho^G = A_G$ and hence by the forcing theorem, there is a $p\in G$ and a $\sigma$ such that $M\vDash p \Vdash \lq\lq \dot{c} = \sigma \supseteq \check{s}"$ and $(\mathbb{P}, p, \sigma) \in \tau_A\cap M$, which implies $(\mathbb{P}, p, \dot{c}) \in \tau_A \cap M$, namely $(\dot{c}, p) \in \rho_A$. But the value of $\dot{c}$ will be decided within Cohen forcing and by the definition of $\tau_A$, we may assume that $p$ is a condition of Cohen forcing extending $s$. Hence for any Cohen real $c'$ over $M$ with $p \subseteq c'$ in $M[G]$, $c'$ is in $A_G$. Since the set of all Cohen reals over $M$ is comeager in $M[G]$, this is what we desired. 
\renewcommand{\qedsymbol}{$\square$ (Claim~\ref{local Baire})}
\end{proof}

We now finish the proof of Theorem~\ref{Blackwell Baire property} by showing that $A$ has the Baire property. Let $G$ be an element of $\mathrm{St}(\mathbb{P})$ in $V$ such that the conclusions of Claim~\ref{local Baire} hold. By the first item of Claim~\ref{local Baire}, the structure $( \omega, 2^{\omega} \cap M[G], 0, 1, + , \times , \in ,  A_G )$ is an elementary substructure of $\mathcal{A}^2_A$. Since the Baire property for $A$ can be described in the structure $\mathcal{A}^2_A$ in this language and $A_G$ has the Baire property in $M[G]$, $A$ also has the Baire property, as desired. 
\renewcommand{\qedsymbol}{$\square$ (Theorem~\ref{Blackwell Baire property})}
\end{proof}

\section{$\BlADR$ and $\infty$-Borel codes}\label{sec:infBorelcodes}

In this section, we show that $\BlADR$ implies every set of reals is $\infty$-Borel, and that $\BlADR$ and $\DC$ imply that every set of reals is strongly $\infty$-Borel, an important step towards the equivalence between $\ADR$ and $\BlADR$ in $\ZF$+$\DC$. 

Throughout this section, we fix a fine normal measure $U$ on $\wp_{\omega_1}(\R)$ which exists by Theorem~\ref{R sharp}.

We first introduce the Vop\v{e}nka algebra and its variant, which is a main tool for our arguments. 
We define the Vop\v{e}nka algebra and its variant for $\mathrm{HOD}_X$, where $X$ is an arbitrary set. Here $\mathrm{OD}_X$ is the class of all sets that are ordinal definable with a parameter $X$, and $\mathrm{HOD}_X$ is the class of all sets $a$ where any element of the transitive closure of $\{ a\}$ is in $\mathrm{OD}_X$. 

Take an arbitrary set $X$ and fix an injection $i_X\colon \text{OD}_X \to \text{HOD}_X$ which is ordinal definable with the parameter $X$. Then consider the {\em Vop\v{e}nka algebra} $\mathbb{P}_{\mathrm{V}, X}$ in $\mathrm{HOD}_X$ as follows: $\mathbb{P}_{\mathrm{V}, X} = \{ i_X( A) \mid A \in \mathrm{OD} \text{ and } A \subseteq \wp(\omega)\}$. For $p, q \in \mathbb{P}_{\mathrm{V}, X} $, $p\le q $ if $i_X^{-1}(p) \subseteq i_X^{-1}(q)$. 
It is easy to see that the properties of $\mathbb{P}_{\mathrm{V}, X}$ do not depend on the choice of $i_X$, i.e., if there are two such injections, then the corresponding two preorders are isomorphic in $\mathrm{HOD}_X$. 
Vop\v{e}nka~\cite{MR0444473} proved that $\mathbb{P}_{\mathrm{V}, \emptyset}$ is a complete Boolean algebra in HOD (when $X = \emptyset$) and each real in $V$ induces a $\mathbb{P}_{\mathrm{V}, \emptyset}$-generic filter over HOD in the following way: For each real $x$ in $V$, the set $G_x = \{ p \in \mathbb{P}_{\mathrm{V}, \emptyset} \mid x \in i_{\emptyset}^{-1}(p) \}$ is a $\mathbb{P}_{\mathrm{V}, \emptyset}$-generic filter over HOD and $x \in \text{HOD}[G_x]$.\footnote{Here is the argument for $x \in \text{HOD}[G_x]$: For each $n\in \omega$ and $i=0,1$, let $N_{n,i} = \{ y \in 2^{\omega} \mid y(n)=i\}$. Then the function $(n,i) \mapsto i_{\emptyset} (N_{n,i})$ is OD and each $i_{\emptyset} (N_{n,i})$ is in HOD. So the function $(n,i) \mapsto i_{\emptyset} (N_{n,i})$ is in HOD. Now for all $n$ and $i$, $x(n) = i \iff x \in N_{n,i} \iff i_{\emptyset} (N_{n,i}) \in G_x$. Hence $x$ can be computed from $G_x$ and the function $(n,i) \mapsto i_{\emptyset} (N_{n,i})$. Therefore, $x \in \text{HOD}[G_x]$.} Conversely, if $G$ is a $\mathbb{P}_{\mathrm{V}, \emptyset}$-generic filter over HOD, then the set $\bigcap \{ i_{\emptyset}^{-1} (p) \mid p \in G\}$ is a singleton. We call the element of the singleton a {\it Vop\v{e}nka real over HOD} and denote it by $y_G$. Then $y_{G_x} = x$ for each real $x$ in $V$. The analogue of the above results holds for $\mathrm{HOD}_X$ for an arbitrary set $X$.

Recall the definition and basic properties of $\infty$-Borel codes from \S\,\ref{sec:Borel codes}. 
As in \S\,\ref{sec:Borel codes}, we often identify an $\infty$-Borel code $\phi$ with a pair $(S, c)$ where $S$ is a tree on some ordinal $\gamma$ and $c$ is a function from the set of terminal nodes of $S$ to natural numbers: If $c(t) = n$, then the atomic sentence $\mathbf{a}_{n}$ is attached on the terminal node $t$ of $S$. 
If $(S,c)$ is an $\infty$-Borel code, then we often write just $S$ to indicate the pair $(S,c)$. 

We now introduce a variant of the Vop\v{e}nka algebra, namely the {\em Vop\v{e}nka algebra with $\infty$-Borel codes}. 
Given a set $X$, consider the following preorder $\mathbb{P}_{\mathrm{V},X}^{\ast} $ in $\mathrm{HOD}_X$: Conditions of $\mathbb{P}_{\mathrm{V},X}^{\ast} $ are $\infty $-Borel codes in $\mathrm{HOD}_X$ where the underlying tree of the code is on some ordinal $\gamma$ which is less than $\Theta$, 
and for $\phi , \psi$ in $\mathbb{P}_{\mathrm{V},X}^{\ast} $, $\phi \le \psi$ if $B_{\phi} \subseteq B_{\psi}$ in $V$.\footnote{As in \S\,\ref{sec:Borel codes}, for any $\infty$-Borel code $S$ in $\mathrm{HOD}_X$, there is an $\infty$-Borel code $T$ in $\mathrm{HOD}_X$ such that $B_S = B_T$ and the underlying tree of the code $T$ is on some ordinal $\gamma$ which is less than $\Theta$.  Hence the restriction of ordinals for $\infty$-Borel codes will not affect the structure of this preorder.}
Then we can prove the analogue of Vop\v{e}nka's theorem in exactly the same way: for $\phi , \psi \in \mathbb{P}_{\mathrm{V}, X}^{\ast}$, we write $\phi \equiv \psi$ if $\phi \le \psi$ and $\psi \le \phi$. 
\begin{Thm}(Folklore?)\label{folk}
Let $X$ be an arbitrary set. 

\begin{enumerate}
\item The quotient $\mathbb{P}_{\mathrm{V},X}^{\ast} \, / \equiv$ is a complete Boolean algebra in $\mathrm{HOD}_X$. 

\item For each real $x$ in $V$, the set $G_x = \{ \phi \in \mathbb{P}_{\mathrm{V},X}^{\ast} \mid x \in B_{\phi}\}$ is $\mathbb{P}_{\mathrm{V},X}^{\ast} $-generic over $\mathrm{HOD}_X$ and $\text{HOD}_X[x] = \text{HOD}_X[G_x]$. Conversely, if $G$ is a $\mathbb{P}^{\ast}_{\mathrm{V},X}$-generic filter over $\mathrm{HOD}_X$, then the set $\bigcap \{ B_{\phi} \mid \phi \in G \} $ in $V$ is a singleton and we call the real in the singleton a {\it Vop\v{e}nka real over $\mathrm{HOD}_X$} and denote it by $y_G$. Then $\text{HOD}_X[y_G] = \text{HOD}_X[G]$ and $y_{G_x} = x$ for each $G$ and $x$.
\end{enumerate}
\end{Thm}

\begin{proof}
The proof is exactly the same as for the Vop\v{e}nka algebra which can be found, e.g., in Jech's textbook~\cite[Theorem~15.46]{Jech}.
\end{proof}

The difference between $\mathbb{P}_{\mathrm{V},X}$ and $\mathbb{P}_{\mathrm{V},X}^{\ast} $ is that $\text{HOD}_X[y_G] \neq \text{HOD}_X[G]$ in general 
for $\mathbb{P}_{\mathrm{V},X}$ while it is always the case that $\text{HOD}_X[y_G] = \text{HOD}_X[G]$ for $\mathbb{P}_{\mathrm{V},X}^{\ast} $. This is because the injection $i_X$ is not in $\mathrm{HOD}_X$ in general while the definition of $\mathbb{P}_{\mathrm{V},X}^{\ast} $ does not refer to $i_X$. 
For our purpose, we will use $\mathbb{P}_{\mathrm{V},X}^{\ast} $.

\begin{Thm}\label{Bl-reg}

Assume $\BlADR$. Then every set of reals is $\infty$-Borel. 
\end{Thm}

\begin{proof}
We use the arguments for the following theorem by Woodin:
\begin{Thm}[Woodin]\label{Woodin}
Assume $\AD$ and that every relation on the reals can be uniformized. Then every set of reals is $\infty$-Borel.
\end{Thm}

In the proof of Theorem~\ref{Woodin} (which can be found in \cite[Theorem~11.18]{ADplus}), Woodin used the Martin measure on Turing degrees, that exists assuming $\AD$. Since we do not know if the Martin measure on Turing degrees exists assuming only $\BlADR$, we instead use the fine normal measure $U$ on $\wp_{\omega_1} (\R)$, that we fixed at the beginning of this section, to implement Woodin's arguments. 

Let $A$ be an arbitrary set of reals. We show that $A$ is $\infty$-Borel.

By Theorem~\ref{Blackwell Baire property}, every set of reals has the Baire property. Hence by Corollary~\ref{cor:category}, every set of reals is $\mathbb{P}$-Baire for any $\mathbb{P} \in H_{\omega_1}$. We freely use this fact later. 
We fix a simple coding of elements of $H_{\omega_1}$ by reals and if we say \lq \lq a real $x$ codes $\ldots$", then we refer to this coding. 

Let $\tau_A$ and $R_A$ be as follows:
\begin{align*}
\tau_A = \{ (\mathbb{P}, p, \sigma) \in H_{\omega_1} \mid& \ \text{$\sigma$ is a $\mathbb{P}$-name for a real and}	\\
&\bigl(\forall^{\infty} G \in \mathrm{St}(\mathbb{P}) \bigr) \ p\in G \implies f_{\sigma}(G) \in A \}, \\
R_A = \{ (x,y) \mid &\text{ if $x$ codes a $(\mathbb{P}, p , \sigma) \in \tau_A$, then $y$ codes a $(D_i \mid i < \omega )$} \\
& \text{ such that } (\forall i) \ \text{$D_i$ is dense in $\mathbb{P}$ and }\\
& \bigl(\forall G \in \mathrm{St}(\mathbb{P}) \bigr) \ \bigl(p \in G, (\forall i ) \ G \cap D_i \neq \emptyset \implies f_{\sigma}(G) \in A\bigr)\},
\end{align*}
where $\lq \lq \bigl(\forall^{\infty} G \in \mathrm{St}(\mathbb{P}) \bigr)\ldots "$ means \lq \lq For comeager many $G$ in $\mathrm{St}(\mathbb{P}), \ldots$". Note that the term relation $\tau_A$ defined here is different from the one in the proof of Theorem~\ref{Blackwell Baire property} in the sense that now we use comeagerness for the quantifier $\forall^{\infty}$ instead of Blackwell comeagerness for the quantifier $\forall^{\infty}_{\mathrm{B}}$. 

Let $F_A$ uniformize $R_A$ and $\Gamma_A$ be the graph of $F_A$, i.e., $\Gamma_A = \{ (x, s) \mid s \in 2^{<\omega}, F_A(x) \supseteq s \}$. Define $\tau_{\Gamma_A}$ as follows:
\begin{align*}
\tau_{\Gamma_A} = \{ (\mathbb{P}, p, \sigma , s) \in H_{\omega_1} \mid &\  \text{$\sigma$ is a $\mathbb{P}$-name for a real and }	\\
&\bigl(\forall^{\infty} G \in \mathrm{St}(\mathbb{P}) \bigr) \ p\in G \implies \bigl( f_{\sigma}(G), s\bigr) \in \Gamma_A\}, 
\end{align*}
here we also use comeagerness for the quantifier $\forall^{\infty}$. 

Recall that $A^{\mathrm{c}}$ is the complement of $A$ (i.e., $A^{\mathrm{c}} = 2^{\omega} \setminus A$) and we define and construct $\tau_{A^{\mathrm{c}}}, R_{A^{\mathrm{c}}}, F_{A^{\mathrm{c}}}, \Gamma_{A^{\mathrm{c}}}$, and $\tau_{\Gamma_{A^{\mathrm{c}}}}$ similarly as above. 

The following is the key point:
\begin{Claim}[Woodin]\label{key claim}

Let $M$ be a transitive subset of $H_{\omega_1}$ and $(M, \in , \tau_A, \tau_{\Gamma_A})$ is a model of $\ZFC$.\footnote{Here it satisfies Comprehension scheme and Replacement scheme for formulas in the language of set theory with predicates for $\tau_A$ and $\tau_{\Gamma_A}$.} 
Let $(\mathbb{P}, p, \sigma) \in M \cap \tau_A$. 
Then for every $\mathbb{P}$-generic filter $G$ over $M$, if $p$ is in $G$, then $\sigma^G \in A$. The same holds for $A^{\mathrm{c}}$.
\end{Claim}

\begin{proof}[Proof of Claim~\ref{key claim}]

Let $\mathbb{Q} = \left[\text{Coll}\bigl(\omega, \text{TC} (\mathbb{P}) \bigr) \right]^M$, where $ \text{TC} (\mathbb{P})$ is the transitive collapse of $\{\mathbb{P} \}$ and $\text{Coll}\bigl(\omega, \text{TC} (\mathbb{P}) \bigr)$ is the standard forcing collapsing $\text{TC} (\mathbb{P})$ into a countable set with finite sets as conditions. 
Since $\mathbb{P}, p, \sigma$ are countable in $M^{\mathbb{Q}}$, there is a $\mathbb{Q}$-name $\sigma'$ for a real in $M$ coding the triple $(\mathbb{P}, p, \sigma)$. 

\begin{Subclaim}\label{subclaim}
There is a $\mathbb{Q}$-name $\rho$ for a real in $M$ such that in $V$, for comeager many $H$ in $\mathrm{St}(\mathbb{Q})$, $f_{\rho}(H) = F_A (f_{\sigma'}(H))$.
\end{Subclaim}
\begin{proof}[Proof of Subclaim~\ref{subclaim}]
\renewcommand{\qedsymbol}{$\square \ (\text{Subclaim~\ref{subclaim})}$}

Since every set of reals has the Baire property, by Corollary~\ref{cor:category}, the map $f \colon H \mapsto F_A \bigl(f_{\sigma'}(H)\bigr)$ is continuous on a comeager set in $\mathrm{St}(\mathbb{Q})$, i.e., Baire measurable. 
Let $\rho = \tau_f$ where the notation $\tau_f$ is from Lemma~\ref{name-Baire measurable}.
Then $\rho$ is a $\mathbb{Q}$-name for a real because the map $f$ is Baire measurable as we observed. 
Moreover, $\rho$ is in $M$ because 
\begin{align*}
(\check{n} , q ) \in \rho \iff (\exists s \in 2^{<\omega}) \  \bigl( s(n) = 1 \text{ and } \bigl(\mathbb{Q}, q , (\sigma' , s) \bigr)\in \tau_{\Gamma_A} \bigr)
\end{align*}
and the right hand side of the equivalence is definable in $(M, \tau_A, \tau_{\Gamma_A})$, which is a model of $\ZFC$ by assumption. 
Finally, by Lemma~\ref{name-Baire measurable}, it is easy to see that for comeager many $H$ in $\mathrm{St}(\mathbb{Q})$, $f_{\rho}(H) = F_A (f_{\sigma'}(H))$.
\end{proof}

Now let $G$ be a $\mathbb{P}$-generic filter over $M$ with $p \in G$. We show that $\sigma^G = f_{\sigma} (G) \in A$. 
Take a $\mathbb{Q}$-generic filter $H$ over $M[G]$ with $\rho^H = F_A (\sigma'^H)$. This is possible by Subclaim~\ref{subclaim} and that $M[G] \subseteq H_{\omega_1}$. 
Then $G$ is also a $\mathbb{P}$-generic filter over $M[H]$ and $F_A (\sigma'^H ) = \rho^H \in M[H]$. 
But by the definition of $F_A$, $F_A (\sigma'^H)$ codes a sequence $(D_i \mid i \in \omega )$ such that $D_i$ is a dense subset of $\mathbb{P}$ in $M[H]$ for each $i\in \omega$ and for any $G'$ in $ \mathrm{St}(\mathbb{P})$, if $G' \cap D_i \neq \empty$ for each $i$, then $f_{\sigma}(G') \in A$. 
But $G$ is a $\mathbb{P}$-generic filter over $M[H]$ and each $D_i$ is in $M[H]$. Hence $G\cap D_i \neq \emptyset$ for each $i\in \omega$ and $f_{\sigma} (G) \in A$, as desired.
\renewcommand{\qedsymbol}{$\square \ (\text{Claim~\ref{key claim}})$}
\end{proof}

Let $X = (A, \tau_A, \tau_{\Gamma_A}, \tau_{A^{\mathrm{c}}}, \tau_{\Gamma_{A^{\mathrm{c}}}})$. Recall that $U$ is the fine normal measure on $\wp_{\omega_1} (\R)$ we fixed at the beginning of this section. Let $N= \mathrm{L}(X, \mathbb{R})[U]$. Since the statement \lq \lq a real is in the decode of an $\infty$-Borel code" is absolute between transitive models of $\ZF$ as in \S\,\ref{sec:Borel codes} and $N$ contains all the reals, if $A$ is $\infty$-Borel in $N$, so is in $V$. 

From now on, we work in $N$ and prove that $A$ is $\infty$-Borel in $N$, which completes the proof of the theorem. The benefit of working in $N$ is that we have $\DC$ in $N$ because $\DC_{\mathbb{R}}$ implies $\DC$ in $N$ while $\DC$ might fail in $V$ in general. Note that $U \cap N$ is a fine normal measure on $\wp_{\omega_1}(\mathbb{R})$ in $N$. For simplicity, we write $U$ to denote $U \cap N$.

We will find a set of ordinals $S$ and a formula $\phi$ such that for any real $x$, 
\begin{align}\label{equiv_wanted}
x \in A \iff \LL[S,x] \vDash \phi (x).
\end{align}
By Lemma~\ref{equiv_code}, this would imply that $A$ is $\infty$-Borel.

For each $a$ in $\wp_{\omega_1}(\mathbb{R})$, let $M_a, \mathbb{Q}^{\ast}_a$, and $b_a$ be as follows:
\begin{align*}
M_a = & \,  \text{HOD}^{\LL_{\omega_1}[X] (a)}_X,\\
\mathbb{Q}^{\ast}_a =  & \,  \mathbb{P}_{\mathrm{V},X}^{\ast} \text{ in }M_a,\\
b_a = & \sup \, \{ q \in \mathbb{Q}^{\ast}_a \mid (\mathbb{Q}^{\ast}_a, q , \dot{y_G} ) \in \tau_A \}  \text{ in }M_a,
\end{align*}
where $\dot{y_G}$ is a canonical $\mathbb{Q}^{\ast}_a$-name for a Vop\v{e}nka real given in Theorem~\ref{folk}.

Note that $M_a$ is a transitive subset of $H_{\omega_1}$ and $(M_a, \tau_A, \tau_{\Gamma_{A}})$ and $(M_a, \tau_{A^{\mathrm{c}}}, \tau_{\Gamma_{A^{\mathrm{c}}}})$ are models of $\ZFC$ because $\LL_{\omega_1}[X] (a)$ is a transitive model of $\ZF$ (to check the power set axiom, we use the condition that there is no injection from $\omega_1$ to the reals by $\BlAD$ and Corollary~\ref{cor:BlAD-measure}). 
Note also that $b_a$ exists for each $a$ because $\mathbb{Q}^{\ast}_a \, / \! \equiv$ is a complete Boolean algebra in $M_a$ by Theorem~\ref{folk}. 

Then we claim that for each $a \in \wp_{\omega_1}(\mathbb{R})$ and real $x$ which induces the filter $G_x$ that is $\mathbb{Q}^{\ast}_a$-generic filter over $M_a$, $x \in A \iff b_a \in G_x$. 
Fix such $a$ and $x$. Assume $b_a \in G_x$. We show that $x \in A$. If we apply Claim~\ref{key claim} to $M= M_a, (\mathbb{P}, p, \tau)= (\mathbb{Q}^{\ast}_a, b_a, \dot{y_G})$, and $G = G_x$, then we get $x \in A$ because $y_{G_x} = x$ as in Theorem~\ref{folk}. 
For the converse, we assume $b_a$ is not in $G_x$ and prove that $x$ is not in $A$.  Let $\overline{b_a}$ be the one corresponding to $b_a$ for $A^{\mathrm{c}}$ instead of for $A$, i.e., 
\begin{align*}
\overline{b_a} = \text{sup} \, \{ q\in \mathbb{Q}^{\ast}_a \mid (\mathbb{Q}^{\ast}_a, q, \dot{y_G} ) \in \tau_{A^{\mathrm{c}}}\}.
\end{align*}
Then $b_a \vee \overline{b_a} = \mathbf{1}$. This is because $f_{\dot{y_G}}^{-1} (A)$ has the Baire property in $\mathrm{St}(\mathbb{Q}^{\ast}_a)$, since $\mathbb{Q}^{\ast}_a$ is countable in $V$ and $A$ is $\mathbb{Q}^{\ast}_a$-Baire. Since $b_a \notin G_x$ and $G_x$ is $\mathbb{P}_{\mathrm{V},X}^{\ast}$-generic over $M_a$, $\overline{b_a}$ is in $G_x$. 
Hence we can apply Claim~\ref{key claim} to $M_a, A^{\mathrm{c}}, (\mathbb{Q}^{\ast}_a, \overline{b_a}, \dot{y_G})$, and $G_x$ to get $x \in A^{\mathrm{c}}$, i.e., $x$ is not in $A$, as desired. 

Fix an $a\in \wp_{\omega_1}(\mathbb{R})$. Note that since $\mathbb{P}_{\mathrm{V},X}^{\ast}$ is the Vop\v{e}nka algebra with $\infty$-Borel codes defined in $M_a$, any real $x$ in $\mathrm{L}_{\omega_1}[X](a)$ induces a $\mathbb{P}_{\mathrm{V},X}^{\ast}$-generic filter $G_x$ over $M_a$. Hence for any real $x$ in $\mathrm{L}_{\omega_1}[X](a)$, $x \in A  \iff b_a \in G_x$. 

Now we use this local equivalence in $\LL_{\omega_1}[X] (a)$ to get the global equivalence~(\ref{equiv_wanted}) by taking the ultraproduct of $M_a$ via $U$. 
Let $M_{\infty}, \mathbb{Q}_{\infty}, b_{\infty}$ be as follows:
\begin{align*}
M_{\infty} = \prod_{a \in \wp_{\omega_1} (\R)} M_a  \, /\, U, \ \mathbb{Q}_{\infty} =  [a \mapsto \mathbb{Q}^{\ast}_a]_U , b_{\infty} = [a \mapsto b_a]_U. 
\end{align*}
Note that \L o\'s's theorem holds for $M_{\infty}$ because there is a canonical function mapping $a$ to a well-order on $M_a = \text{HOD}^{\LL_{\omega_1}[X] (a)}_X$. 
By $\DC$ (in $N$), $M_{\infty}$ is wellfounded. So we may assume $M_{\infty}$ is transitive. Hence, $M_{\infty}$ is a transitive model of $\ZFC$, $\mathbb{Q}_{\infty}$ is a preorder consisting of $\infty$-Borel codes, and $b_{\infty} \in \mathbb{Q}_{\infty}$. 

We claim that for each real $x$, $x \in A \iff x\in B_{b_{\infty}}$. 
This will establish the equivalence~(\ref{equiv_wanted}) because the pair $(\mathbb{Q}_{\infty},  b_{\infty})$ can be seen as a set of ordinals since they are objects in the transitive model $M_{\infty}$ of $\ZFC$.

Let us fix a real $x$. By the fineness of $U$, $x \in a$ for almost all $a$ w.r.t. $U$. Then
\begin{align*}
x \in A \iff & b_a \in G_x \text{ for almost all $a$} \\
\iff & x \in B_{b_a} \text{ for almost all $a$} \\
\iff & x \in B_{b_{\infty}},
\end{align*}
where the first equivalence is by the local equivalence we have seen, the second equivalence is by the definition of $G_x$ in $\mathbb{Q}^{\ast}_a$, and the third equivalence follows from \L o\'s's theorem for $\prod_{a\in \wp_{\omega_1} (\R)} M_a[x] \, / \, U$ (note that $M_a[x] = M_a[G_x]$ by Theorem~\ref{folk} and we can prove \L o\'s's theorem for $\prod_{a \in \wp_{\omega_1}(\R)} M_a[x] \, / \, U$ in the same way as for $\prod_{a \in \wp_{\omega_1}(\R) } M_a \, / \, U$). 

This completes the proof of Theorem~\ref{Bl-reg}.
\end{proof}

Next, we show that every set of reals is strongly $\infty$-Borel assuming $\BlADR$ and $\DC$. Before giving a definition of strongly $\infty$-Borel codes, we start with a lemma:
\begin{Lem}\label{bound j ADR}
Assume $\BlADR$ and $\DC$. Let $j \colon V \to \mathrm{Ult}(V,U)$ be the ultrapower map via $U$. Then $j (\omega_1) = \Theta$. 
\end{Lem}

\begin{proof}
We first show that $j (\omega_1) \ge \Theta$. Let $\alpha $ be an ordinal less than $\Theta$ and $R$ be a prewellorder on the reals of length $\alpha$. Define $f \colon \wp_{\omega_1}(\mathbb{R}) \to \omega_1$ as follows: For $a \in \wp_{\omega_1}(\mathbb{R})$, $f(a)$ is the length of the prewellorder $R \cap (a\times a)$ on $a$. Since $a$ is countable, $f(a)$ is also countable. Hence $f \in_U c_{\omega_1}$, where $\in_U$ is the membership relation for $\mathrm{Ult}(V,U)$ and $c_{\omega_1}$ is the constant function on $\wp_{\omega_1}(\mathbb{R})$ with value $\omega_1$. 

We show that the structure $([f]_U , \in)$ is isomorphic to $(\alpha , \in )$ and hence $[f]_U = \alpha$, which implies $\alpha < j(\omega_1)$ because $f \in_U c_{\omega_1}$. For any $a\in \wp_{\omega_1}(\mathbb{R})$, let $\pi(a)$ be the transitive collapse of $\bigl(a , R\cap (a \times a) \bigr)$ into $\bigl(f(a) , \in \bigr)$. Then by \L o\'s's Theorem for simple formulas, $[\pi]_U$ is an isomorphism between $\bigl([\mathrm{id}]_U, j (R) \cap ([\mathrm{id}]_U \times [\mathrm{id}]_U)\bigr)$ and $([f]_U, \in)$, where id is the identity function on $\wp_{\omega_1}(\mathbb{R})$. 

\begin{Claim}\label{basic j}
The identity function id represents $\mathbb{R}$, i.e., $[\mathrm{id}]_U = \mathbb{R}$. 
\end{Claim}

\begin{proof}[Proof of Claim~\ref{basic j}]
By the fineness of $U$, for any real $x$, $\{ a \mid x \in a \} \in U$. Hence $[c_x]_U \in [\mathrm{id}]_U$. By the countable completeness of $U$, $[c_x]_U = x$ and hence $x \in [\mathrm{id}]_U$ for any real $x$. 
Suppose $g$ is a function on $\wp_{\omega_1}(\mathbb{R})$ with $g \in_U \mathrm{id}$. Then by the normality of $U$, there is a real $x$ such that $\{ a \mid x = g(a) \} \in U$, i.e., $c_x =_U g$. Hence $[g]_U =x $ and $[g]_U$ is a real, which finishes the proof. 
\renewcommand{\qedsymbol}{$\square \ (\text{Claim~\ref{basic j}})$}
\end{proof}

By Claim~\ref{basic j}, we have $[\mathrm{id}]_U = \mathbb{R}$ and $j(R) \cap ([\mathrm{id}]_U \times [\mathrm{id}]_U)\bigr) = R$. Since $\bigl([\mathrm{id}]_U, j (R) \cap ([\mathrm{id}]_U \times [\mathrm{id}]_U)\bigr)$ and $([f]_U, \in)$ are isomorphic, $([f]_U, \in)$ is isomorphic to $(\mathbb{R}, R)$, which is isomorphic to $(\alpha, \in )$, as desired. Hence $\alpha < j(\omega_1)$ and $j(\omega_1) \ge \Theta$. 

Next, we show that $j(\omega_1) \le \Theta$. Let $f$ be a function from $\wp_{\omega_1}(\mathbb{R})$ to $\omega_1$. We show that $[f]_U < \Theta$. By uniformization for every set of reals, there is a function $e$ from the reals to themselves such that if a real $x$ codes an $a \in \wp_{\omega_1}(\mathbb{R})$, then $e(x)$ codes $f(a)$. Let $S$ be an $\infty$-Borel code for the graph $\Gamma_e$ of $e$ which exists by Theorem~\ref{Bl-reg}. 
\begin{Claim}\label{Woodin calculation}
For all $a \in \wp_{\omega_1}(\mathbb{R})$, $f(a) < \Theta^{\mathrm{L}[S](a)}$. 
\end{Claim}

\begin{proof}[Proof of Claim~\ref{Woodin calculation}]
Note that $\wp(x) \cap \mathrm{L}[S](a)$ is countable in $V$ for any $x \in H_{\omega_1} \cap \mathrm{L}[S](a)$: In fact, since $a$ is countable in $V$, the model $\mathrm{L}[S](a)$ is well-orderable in $V$, and so is $\wp (x) \cap \mathrm{L}[S](a)$ in $V$, which cannot be uncountable because that would give us an injection from $\omega_1$ to $\R$ in $V$, which cannot exist assuming $\BlADR$ by Corollary~\ref{cor:BlAD-measure} and Proposition~\ref{easy_Blackwell}. Hence there is a $\mathrm{Coll}(\omega, a)$-generic $g$ over $\mathrm{L}[S](a)$ in $V$. Fix such a $g$. Let $x_g$ be a real coding $a$ from $g$. Then since $S$ is an $\infty$-Borel code for $\Gamma_e$, one can compute whether $e(x_g) \supseteq s$ or not for each finite binary sequence $s$ in $\mathrm{L}[S](a,g)$, hence $e(x_g) \in \mathrm{L}[S](a,g)$. Therefore $f(a)$ is countable in $\mathrm{L}[S](a,g)$. But $\Theta^{\mathrm{L}[S](a)}$ stays an uncountable cardinal in $\mathrm{L}[S](a,g)$. Hence $f(a) < \Theta^{\mathrm{L}[S](a)}$, as desired. 
\renewcommand{\qedsymbol}{$\square \ (\text{Claim~\ref{Woodin calculation}})$}
\end{proof}

By the normality of $U$, the following choice principle holds: For any function $F$ from $\wp_{\omega_1} (\mathbb{R})$ to $V$ such that $\emptyset \neq F(a) \in \mathrm{L}[S](a)$ for almost $a$ with respect to $U$, there is a function $f \colon \wp_{\omega_1}(\mathbb{R}) \to V$ such that $f(a) \in F(a)$ for almost all $a$ with respect to $U$. 
This implies \L o\'s's Theorem for the ultraproduct $\prod_{a \in \wp_{\omega_1} (\R)} \mathrm{L}[S] (a) / U$. 

Let $S^{*} = j(S)$. Then since $[\mathrm{id}]_U = \mathbb{R}$ by Claim~\ref{basic j}, by the \L o\'s's Theorem for the ultraproduct $\prod_{a \in \wp_{\omega_1}(\R)} \mathrm{L}[S] (a) / U$, it follows that $\bigl(\prod_{a \in \wp_{\omega_1}(\R)} \mathrm{L}[S](a) / U, \in_U \bigr)$ is isomorphic to $\bigl(\mathrm{L}[S^{*}](\mathbb{R}) , \in \bigr)$. 
(Note that $\mathrm{Ult}(V,U)$ is wellfounded by $\DC$.) Hence 
\begin{align*}
[f]_U < [a \mapsto \Theta^{\mathrm{L}[S](a)}]_U = \Theta^{\mathrm{L}[S^*](\mathbb{R})} \le \Theta^V,
\end{align*}
as desired. 

This completes the proof of Lemma~\ref{bound j ADR}. 
\end{proof}

As in \S\,\ref{sec:Borel codes}, we often identify an $\infty$-Borel code $\phi$ with a pair $(S, c)$ where $S$ is a tree on some ordinal $\gamma$ and $c$ is a function from the set of terminal nodes of $S$ to natural numbers: If $c(t) = n$, then the atomic sentence $\mathbf{a}_{n}$ is attached on the terminal node $t$ of $S$. 
If $(S,c)$ is an $\infty$-Borel code, then we often write $S$ to indicate the pair $(S,c)$. 

We now introduce strongly $\infty$-Borel codes. 
An $\infty$-Borel code $S = (S, c)$ is {\it strong} if $S$ is a tree on $\gamma$ for some $\gamma < \Theta$ and for any $f \colon \mathbb{R}^{<\omega}  \to \mathbb{R}$ and surjection $\pi \colon \mathbb{R} \to \gamma$, there is an $a\in \wp_{\omega_1} (\R)$ such that $a$ is closed under $f$, $S {\upharpoonright} \pi[a]$ is an $\infty$-Borel code, and $B_{S{\upharpoonright}\pi[a]} \subseteq B_S$, where $S{\upharpoonright}\pi[a]$ is the $\infty$-Borel code of the form $\bigl(S \cap (\pi[a])^{<\omega}, c \upharpoonright S \cap (\pi[a])^{<\omega}\bigr)$. A set of reals $A$ is {\it strongly $\infty$-Borel} if $A = B_S$ for some strong $\infty$-Borel code $S$. The following is a refinement of Lemma~\ref{equiv_code}:
\begin{Thm}\label{equiv strong codes}
${}$

\begin{enumerate}
\item Let $S$ be a strong $\infty$-Borel code and $\gamma < \Theta$ be such that $S$ is a tree on $\beta$ for some $\beta < \gamma$ and $\mathrm{L}_{\gamma}[S,x] \vDash \lq \lq\mathsf{KP}\mathrm{+\Sigma_1\text{-Separation}}"$ for any real $x$. Let $\phi (S, x)$ be a $\Sigma_1$-formula expressing \lq\lq$x \in B_S$" in $\mathrm{L}_{\gamma}[S,x]$. Then for any function $f \colon \mathbb{R}^{<\omega} \to \mathbb{R}$ and surjection $\pi \colon \mathbb{R} \to \gamma$, there is an $a\in \wp_{\omega_1} (\mathbb{R})$ such that $a$ is closed under $f$ and for any real $x$, if $\mathrm{L}_{\overline{\gamma}}[\overline{S},x] \vDash \phi (\overline{S}, x)$, then $\mathrm{L}_{\gamma}[S,x] \vDash \phi (S, x)$, where $\mathrm{L}_{\overline{\gamma}}[\overline{S}]$ is the transitive collapse of the Skolem hull of $\pi[a] \cup \{S\}$ in $\mathrm{L}_{\gamma}[S]$. 

\item Let $\gamma$ be an ordinal with $\gamma < \Theta$, $\phi$ be a $\Sigma_1$-formula, and $S$ be a bounded subset of $\gamma$ such that $\mathrm{L}_{\gamma}[S,x] \vDash \lq\lq \mathsf{KP}\mathrm{+\Sigma_1\text{-Separation}}"$ for any real $x$. Set $A = \{ x \in \mathbb{R} \mid \mathrm{L}_{\gamma}[S,x] \vDash \phi (S,x) \}$. Assume that for any function $f \colon \mathbb{R}^{<\omega} \to \mathbb{R}$ and surjection $\pi \colon \mathbb{R} \to \gamma$, there is an $a\in \wp_{\omega_1}(\mathbb{R})$ such that $a$ is closed under $f$ and for any real $x$, if $\mathrm{L}_{\overline{\gamma}}[\overline{S},x] \vDash \phi (\overline{S}, x)$, then $\mathrm{L}_{\gamma}[S,x] \vDash \phi (S, x)$, where $\mathrm{L}_{\overline{\gamma}}[\overline{S}]$ is the transitive collapse of the Skolem hull of $\pi[a] \cup \{S\}$ in $\mathrm{L}_{\gamma}[S]$. Then $A$ is strong $\infty$-Borel. 
\end{enumerate}
\end{Thm}

Therefore, a given $\infty$-Borel code $S$ is strong if and only if there are \lq stationarily many' $a$ in $\wp_{\omega_1} (\R)$ such that for any real $x$, if $\mathrm{L}_{\overline{\gamma}}[\overline{S},x] \vDash \phi (\overline{S}, x)$, then $\mathrm{L}_{\gamma}[S,x] \vDash \phi (S, x)$, where $\gamma$, $\overline{\gamma}$, and $\overline{S}$ are as in Theorem~\ref{equiv strong codes}.

\begin{proof}[Proof of Theorem~\ref{equiv strong codes}]
This can be done by closely looking at the argument for Lemma~\ref{equiv_code}. The details of the argument for Lemma~\ref{equiv_code} are given in \cite[Theorem~8.7]{ADplus}. 
\end{proof}

\begin{Thm}\label{strong infty}
Assume $\BlADR$ and $\DC$. Then every set of reals is strongly $\infty$-Borel.
\end{Thm}

\begin{proof}
Let $A$ be any set of reals. We show that $A$ is strongly $\infty$-Borel. 
Let $\bigl( ( M_a, \mathbb{Q}^{\ast}_a, b_a) \mid a \in \wp_{\omega_1}(\mathbb{R})\bigr)$ and $(M_{\infty}, \mathbb{Q}^{\ast}_{\infty}, b_{\infty})$ be  as in the proof of Theorem~\ref{Bl-reg}, but we construct them in $V$, not in $N$. Since we have $\DC$ in $V$ now, we can prove the following equivalences in exactly the same way as in the proof of Theorem~\ref{Bl-reg}: For all $a \in \wp_{\omega_1}(\mathbb{R})$ and all real $x$ inducing the filter $G_x$ which is $\mathbb{Q}^{\ast}_a$-generic over $M_a$,
\begin{align*}
 x \in A \iff& b_a \in G_x \ (\text{in } \mathbb{Q}^{\ast}_a).
\end{align*}
Also, 
\begin{align*}
(\forall x \in \mathbb{R}) \ x \in A \iff& b_{\infty} \in G_x \ (\text{in }\mathbb{Q}^{\ast}_{\infty}). 
\end{align*}
For each $a \in \wp_{\omega_1}(\R)$, let $\mathcal{D}_a $ be the set of all dense subsets of $\mathbb{Q}^{\ast}_a$ in $M_a$ and let $\mathcal{D}_{\infty} = [a \mapsto \mathcal{D}_a]_U$. Let $\phi$ be a $\Sigma_1$-formula such that for all $a$, 
\begin{align*}
\phi (\mathbb{Q}^{\ast}_a, b_a, \mathcal{D}_a,x) \iff& \text{ $x$ determines the filter $G_x \subseteq \mathbb{Q}^{\ast}_a$ such that } \\
&\text{ $(\forall D \in \mathcal{D}_a) \ G_x \cap D \neq \emptyset$ and $b_a \in G_x$,}\\
\phi (\mathbb{Q}^{\ast}_{\infty}, b_{\infty}, \mathcal{D}_{\infty},x) \iff& \text{ $x$ determines the filter $G_x \subseteq \mathbb{Q}^{\ast}_{\infty}$such that } \\
&\text{ $(\forall D \in \mathcal{D}_{\infty}) \ G_x \cap D \neq \emptyset$ and $b_{\infty} \in G_x$}. 
\end{align*}

For each $a \in \wp_{\omega_1}(\R)$, let $S_a$ be a set of ordinals coding the triple $(\mathbb{Q}^{\ast}_a, b_a, D_a)$. Let $S_{\infty} = [a \mapsto S_a]_U$. Then $S_{\infty}$ codes the triple $(\mathbb{Q}^{\ast}_{\infty} , b_{\infty}, D_{\infty})$. 
For each $a\in \wp_{\omega_1}(\mathbb{R})$, let $\alpha_a$ be the least ordinal $\alpha$ such that $S_a$ is a bounded subset of $\alpha$ and for all $x \in a$, $\mathrm{L}_{\alpha} [S_a, x]$ is a model of $\mathsf{KP}$+$\Sigma_1$-Separation and let $\alpha_{\infty}$ be the least ordinal $\alpha$ such that $S_{\infty}$ is a bounded subset of $\alpha$ and for all $x \in \mathbb{R}$, $\mathrm{L}_{\alpha} [S_{\infty}, x]$ is a model of $\mathsf{KP}$+$\Sigma_1$-Separation. Note that for every real $x$, by \L o\'s's Theorem for the ultraproduct $\prod_{a\in \wp_{\omega_1}(\R)} M_a [x] / U$, the structure $(\prod_{a\in \wp_{\omega_1}(\R)} \mathrm{L}_{\alpha_a}[S_a,x], \in_U)$ is isomorphic to $(\mathrm{L}_{\alpha_{\infty}} [S_{\infty},x], \in)$. Since each $\alpha_a$ is countable, by Lemma~\ref{bound j ADR}, $\alpha_{\infty} < \Theta$. Also, by the above equivalences, for all $a\in \wp_{\omega_1}(\mathbb{R})$ and all reals $x$, 
\begin{align*}
 x \in A \iff& \mathrm{L}_{\alpha_a}[S_a, x] \vDash  \phi \, (S_a, x) \text{ and }\\
x \in A \iff& \mathrm{L}_{\alpha_{\infty}} [S_{\infty}, x] \vDash \phi \, (S_{\infty}, x).
\end{align*}

By the second item of Theorem~\ref{equiv strong codes}, it suffices to show the following: For any function $f \colon \mathbb{R}^{<\omega} \to \mathbb{R}$ and surjection $\pi \colon \mathbb{R} \to \alpha_{\infty}$, there is an $a\in \wp_{\omega_1}(\mathbb{R})$ such that $a$ is closed under $f$ and for any real $x$, if $\mathrm{L}_{\overline{\alpha_{\infty}}}[\overline{S_{\infty}},x] \vDash \phi (\overline{S_{\infty}}, x)$, then $\mathrm{L}_{\alpha_{\infty}}[S_{\infty},x] \vDash \phi (S_{\infty}, x)$, where $\mathrm{L}_{\overline{\alpha_{\infty}}}[\overline{S_{\infty}}]$ is the transitive collapse of the Skolem hull of $\pi[a] \cup \{S_{\infty}\}$ in $\mathrm{L}_{\alpha_{\infty}}[S_{\infty}]$. 

Let us fix $f\colon \mathbb{R}^{<\omega} \to \mathbb{R}$ and $\pi \colon \mathbb{R} \to \alpha_{\infty}$. Since $ x \in A \iff \mathrm{L}_{\alpha_b}[S_b, x] \vDash  \phi \, (S_b, x)$ for each real $x$ and $b\in \wp_{\omega_1}(\mathbb{R})$, the following claim completes the proof:
\begin{Claim}\label{local global}
There are $a$ and $b $ in $\wp_{\omega_1}(\mathbb{R})$ such that $a$ is closed under $f$ and $(X_a, \in)$ is isomorphic to $(\mathrm{L}_{\alpha_b}[S_b], \in)$, where $X_a$ is the Skolem hull of $\pi[a] \cup \{S_{\infty}\}$ in~$\mathrm{L}_{\alpha_{\infty}}[S_{\infty}]$. 
\end{Claim}
\begin{proof}[Proof of Claim~\ref{local global}]
Let $\Gamma_f$ be the graph of $f$, i.e., $\Gamma_f = \{ (x , s) \in \mathbb{R} \times 2^{<\omega} \mid f(x) \supseteq s\}$. For each $b$, consider the following game $\mathcal{G}^b_{\text{isom}}$ in $\mathrm{L}[S_b, S_{\infty}, \Gamma_f, \pi]$ (which is a model of $\ZFC$ closed under $f$): In $\omega$ rounds,
\begin{enumerate}
\item player I and II jointly produce a countable elementary substructure $X$ of $\mathrm{L}_{\alpha_b}[S_b]$, 

\item player II produces an $a \in \wp_{\omega_1}(\mathbb{R})$ which is closed under $f$, and 

\item player II tries to construct an isomorphism between $(X, \in )$ and $(X_a, \in)$, where $X_a$ is the Skolem hull of $\pi[a] \cup \{S_{\infty}\}$ in~$\mathrm{L}_{\alpha_{\infty}}[S_{\infty}]$.
\end{enumerate}
Player II wins if she succeeds to construct an isomorphism between $(X, \in )$ and $(X_a, \in~)$. This is an open game on some 
set in $\mathrm{L}[S_b, S_{\infty}, \Gamma_f, \pi]$, which is a model of $\ZFC$. Hence it is determined in $\mathrm{L}[S_b, S_{\infty}, \Gamma_f, \pi]$. 
\begin{Subclaim}\label{local global sub}
There is a $b\in \wp_{\omega_1}(\mathbb{R})$ such that player II has a winning strategy in the game $\mathcal{G}^b_{\text{isom}}$ in $\mathrm{L}[S_b, S_{\infty}, \Gamma_f, \pi]$. 
\end{Subclaim}

\begin{proof}[Proof of Subclaim~\ref{local global sub}]
To derive a contradiction, suppose there is no $b$ such that player II has a winning strategy in the game $\mathcal{G}^b_{\text{isom}}$ in $\mathrm{L}[S_b, S_{\infty}, \Gamma_f, \pi]$. For each $b$, by the determinacy of the game $\mathcal{G}^b_{\text{isom}}$, player I has a winning strategy in the game $\mathcal{G}^b_{\text{isom}}$ in $\mathrm{L}[S_b, S_{\infty}, \Gamma_f, \pi]$. Let $j\colon V \to \mathrm{Ult}(V,U)$ be the ultrapower map. Then by \L o\'s's Theorem, the structure $\bigl(\prod_{b} \mathrm{L}[S_b, S_{\infty}, \Gamma_f, \pi], \in_U, \Gamma_f, \pi \bigr)$ is isomorphic to $\bigl(\mathrm{L}[S_{\infty}, j (S_{\infty}), \Gamma_f , j(\pi)], \in , \Gamma_f, j(\pi)\bigr)$. Then the game $\mathcal{G}^{\infty}_{\text{isom}} = [b \mapsto \mathcal{G}^b_{\text{isom}}]_U$ is an open game on some set 
in a $\ZFC$ model $\mathrm{L}[S_{\infty}, j (S_{\infty}), \Gamma_f , j(\pi)]$ such that in $\omega$ rounds, 
\begin{enumerate}
\item players I and II jointly produce a countable elementary substructure $Y$ of $\mathrm{L}_{\alpha_{\infty}}[S_{\infty}]$, 

\item player II produces an $a \in \wp_{\omega_1}(\mathbb{R})$ which is closed under $f$, and 

\item player II tries to construct an isomorphism between $(Y, \in )$ and $(Y_a, \in)$, where $Y_a$ is the Skolem hull of $j(\pi)[a] \cup \{ j(S_{\infty})\}$ in $\mathrm{L}_{j(\alpha_{\infty})}[j(S_{\infty})]$. 
\end{enumerate}
Player II wins if she succeeds to construct an isomorphism between $Y$ and $Y_a$. By \L o\'s's Theorem, player I has a winning strategy $\sigma$ in $\mathrm{L}[S_{\infty}, j (S_{\infty}), \Gamma_f , j(\pi)]$. 
Since the game is open, $\sigma$ is also winning in $V$. In $V$, let player II move in such a way that she can arrange that $a$ is closed under $f$, $j[Y] = Y_a$, and $j{\upharpoonright}Y$ is the candidate for the isomorphism. This is possible by a bookkeeping argument. But then player II wins because $j{\upharpoonright}Y$ is an isomorphism between $Y$ and $j[Y]$ and defeats the strategy $\sigma$, contradiction!
\renewcommand{\qedsymbol}{$\square \ (\text{Subclaim~\ref{local global sub}})$}
\end{proof}
Hence there is a $b\in\wp_{\omega_1}(\mathbb{R})$ such that player II has a winning strategy $\tau$ in the game $\mathcal{G}^b_{\text{isom}}$ in $\mathrm{L}[S_b, S_{\infty}, \Gamma_f, \pi]$. Since the game is open, $\tau$ is also winning in $V$. Since $\mathrm{L}_{\alpha_b}[S_b]$ is countable in $V$, we can let player I move in such a way that $X = \mathrm{L}_{\alpha_b}[S_b]$ and let player II follow $\tau$. Since $\tau$ is winning in $V$, there is an $a\in \wp_{\omega_1}(\mathbb{R})$ such that $a$ is closed under $f$ and $\mathrm{L}_{\alpha_b}[S_b] = X $ is isomorphic to $X_a$, as desired. 
\renewcommand{\qedsymbol}{$\square \ (\text{Claim~\ref{local global}})$}
\end{proof}
This completes the proof of Theorem~\ref{strong infty}.
\end{proof}

The definition of strong $\infty$-Borel codes says that given $S, \gamma$ and $\pi$ as in the definition, for \lq stationarily many' $a$ in $\wp_{\omega_1}(\mathbb{R})$, $a$ is \lq good' i.e., $B_{S \upharpoonright \pi [a]} \subseteq B_S$. Using the range-invariant determinacy (that follows from $\BlADR$ by Lemma~\ref{range-invariant}), one can make it for \lq club many' $a$ in $\wp_{\omega_1}(\R)$, $a$ is \lq good' in the same sense:
\begin{Lem}\label{stronger code}
Assume $\BlADR$. Let $S$ be a strong $\infty$-Borel code whose tree is on $\gamma$ for some $\gamma < \Theta$. Let $\pi \colon \R \to \gamma$ be a surjection. Then there is a function $f\colon \mathbb{R}^{<\omega} \to \mathbb{R}$ such that for all $a \in \wp_{\omega_1}( \R)$ closed under $f$,  $S {\upharpoonright} \pi [a]$ is an $\infty$-Borel code and $B_{S{\upharpoonright}\pi [a]} \subseteq B_S$. 
\end{Lem}
\begin{proof}[Proof of Lemma~\ref{stronger code}]
Let us consider the following game: Player I and II choose reals one by one and produce an $\omega$-sequence $\vec{x}$ of reals. Setting $a = \mathrm{ran}(\vec{x})$, player I wins if $S {\upharpoonright} \pi [a]$ is an $\infty$-Borel code and $B_{S{\upharpoonright} \pi [a]} \subseteq B_S$. Since $S$ is a strong $\infty$-Borel code, player I can defeat any strategy for player II, because strategies can be seen as functions from $\mathbb{R}^{<\omega}$ to $\mathbb{R}$ by Claim~\ref{1}. Since the payoff set of this game is range-invariant, by Lemma~\ref{range-invariant}, this game is determined. Hence player I has a winning strategy and by Claim~\ref{1}, there is a function $f$ as desired.
\renewcommand{\qedsymbol}{$\square \ (\text{Lemma~\ref{stronger code}})$}
\end{proof}

The following corollary will be used in the next section:
\begin{Cor}\label{strong code countable intersection}
Assume $\BlADR$+$\DC$. Let $(A_n \mid n < \omega)$ an $\omega$-sequence of sets of reals such that each $A_n$ has a strong $\infty$-Borel code $S_n$ whose tree is on $\gamma$ for some $\gamma < \Theta$. Let $\bigwedge_{n< \omega} S_n$ be the $\infty$-Borel code which codes the intersection $\bigcap_{n < \omega} A_n$ using $(S_n \mid n < \omega)$. Then $\bigwedge_{n<\omega} S_n$ is a strong $\infty$-Borel code whose tree is on $\gamma$. 
The same holds for $\bigvee_{n< \omega} S_n$, the $\infty$-Borel code which codes the union $\bigcup_{n < \omega} A_n$. 
\end{Cor}

\begin{proof}
We argue only for $\bigwedge_{n< \omega} S_n$ (the $\infty$-code $\bigvee_{n< \omega} S_n$ is similarly taken care of). 
Let $\pi \colon \R \to \gamma$ be a surjection. By Lemma~\ref{stronger code}, for each $n < \omega$, there is a function $f_n \colon \R^{<\omega} \to \R$ such that if $a$ is closed under $f_n$, then $S_n \upharpoonright \pi [a]$ is an $\infty$-Borel code and $B_{S_n \upharpoonright \pi [a]} \subseteq B_{S_n}$. Using $\DC$, we choose an $\omega$-sequence $(f_n \mid n < \omega)$ of such functions. Then one can find an $f \colon \R^{<\omega} \to \R$ such that if $a$ is closed under $f$, then $a$ is closed under each $f_n$. Then this $f$ witnesses that $\bigwedge_{n< \omega} S_n$ is a strong $\infty$-Borel code. 
\end{proof}

\section{Equivalence between $\ADR$ and $\BlADR$ assuming $\DC$}\label{sec:equiv}

In this section, we prove the main theorem of this paper: 
\begin{Thm}[Theorem~\ref{main theorem}]\label{almost theorem}
The axioms $\ADR$ and $\BlADR$ are equivalent in $\ZF$+$\DC$. 
\end{Thm}

Since $\ADR$ implies $\BlADR$ by Theorem~\ref{det implies Blackwell det}, the question is whether $\BlADR$ implies $\ADR$ in $\ZF$+$\DC$. 
Woodin proved the following:
\begin{Thm}[Woodin]\label{equiv AD and ADR}
Assume $\AD$ and $\DC$. Then the following are equivalent:
\begin{enumerate}
\item Every set of reals is Suslin, 

\item The axiom $\ADR$ holds, and

\item Every relation on the reals can be uniformized. 
\end{enumerate}
\end{Thm}

\begin{proof}[Proof of Theorem~\ref{equiv AD and ADR}]
See \cite[Theorem~12.12 and Theorem~13.1]{ADplus}. 
\end{proof}

Hence, to prove Theorem~\ref{almost theorem}, it suffices to show that every set of reals is Suslin from $\BlADR$: If every set of reals is Suslin, then by Theorem~\ref{Suslin-co-Suslin-det}, $\AD$ holds. Now by Theorem~\ref{equiv AD and ADR} and Theorem~\ref{uniformization}, $\ADR$ holds assuming $\BlADR$ and $\DC$. 

We try to simulate the arguments for the implication from uniformization to Suslinness in Theorem~\ref{equiv AD and ADR} using $\BlADR$. 
Throughout this section, we fix $U$ as a fine normal measure on $\wp_{\omega_1}(\mathbb{R})$, which exists by Theorem~\ref{R sharp}.

We first prove the key statement towards Theorem~\ref{almost theorem}: Recall that for a natural number $n$ with $n \ge 1$ and a subset $A$ of $\mathbb{R}^{n+1}$, $\exists^{\mathbb{R}}A = \{ x \in \mathbb{R}^n \mid (\exists y \in \mathbb{R}) \ (x,y)\in A\}$.  
\begin{Thm}\label{Becker's bound}
Assume $\BlADR$ and $\DC$. Let $A$ be a subset of $\mathbb{R}^3$ which has a strong $\infty$-Borel code $S$ where the tree of $S$ is on $\gamma < \Theta$. If $\exists^{\mathbb{R}} A$ is a strict well-founded relation on a set of reals, then the length of $\exists^{\mathbb{R}} A$ is less than $\gamma^+$. 
\end{Thm}

\begin{proof}
Let $A, S$, and $\gamma$ be as in the assumptions. We show that the length of $\exists^{\mathbb{R}} A$ is less than $\gamma^+$. 
Fix a surjection $\pi \colon \mathbb{R} \to \gamma$. 

We fix an $f_0$ witnessing the conclusion of Lemma~\ref{stronger code} for the rest of this proof. 
Recall that $U$ is the fine normal measure on $\wp_{\omega_1}(\mathbb{R})$ we fixed at the beginning of this section. Using $\pi$, we can transfer this measure to a fine normal measure on $\wp_{\omega_1} (\gamma)$ as follows: Let $\pi_{\ast} \colon \wp_{\omega_1}(\mathbb{R}) \to \wp_{\omega_1}(\gamma)$  be such that $\pi_{\ast}(a) = \pi [a]$ for each $a\in \wp_{\omega_1}(\mathbb{R})$. For $A \subseteq \wp_{\omega_1}(\gamma)$, $A \in U_{\pi}$ if $\pi_{\ast}^{-1}(A) \in U$. It is easy to check that $U_{\pi}$ is a fine normal measure on $\wp_{\omega_1}(\gamma)$. 

Here is the key lemma for the theorem:
\begin{Lem}\label{key generic embedding}
Let $G$ be $\mathrm{Coll}(\omega, \gamma)$-generic over $V$. Then in $V[G]$, there is an elementary embedding $j \colon \mathrm{L}(\mathbb{R}, S, f_0, \pi) \to \mathrm{L}\bigl(j(\mathbb{R}), j(S), j(f_0), j (\pi) \bigr)$ such that all the reals in $V[G]$ are contained in $\mathrm{L}\bigl(j(\mathbb{R}), j(S), j(f_0), j (\pi) \bigr)$. 
\end{Lem}

\begin{proof}[Proof of Lemma~\ref{key generic embedding}]
The argument is the same as that for the result of Kechris and Woodin~\cite[Theorem~6.2]{MR2463619}. 

In the result of Kechris and Woodin, they assumed $V = \mathrm{L}(\R)$ and $\gamma$ is a reliable, and considered an elementary embedding $j \colon \mathrm{L}(\R) \to \mathrm{L}\bigl( j(\R) \bigr)$ in $V[G]$. In our setting, we assume $\BlADR$ and $\DC$ while $\gamma$ is an arbitrary ordinal less than $\Theta$, and consider an elementary embedding $j \colon \mathrm{L}(\mathbb{R}, S, f_0, \pi) \to \mathrm{L}\bigl(j(\mathbb{R}), j(S), j(f_0), j (\pi) \bigr)$ in $V[G]$. We do not need to assume that $\gamma$ is reliable because we have uniformization for every relation on the reals by Theorem~\ref{uniformization}, which ensures that \L o\'s's Theorem holds for the ultrapower witnessing the desired embedding $j$. 
\renewcommand{\qedsymbol}{$\square \ (\text{Lemma~\ref{key generic embedding}})$}
\end{proof}

We go back to the proof of Theorem~\ref{Becker's bound}. 
Let $G$ and $j$ be as in Lemma~\ref{key generic embedding}. Let $M$ be the model $\mathrm{L}\bigl(j(\mathbb{R}), j(S), j(f_0), j (\pi) \bigr)$. So $j \colon  \mathrm{L}(\mathbb{R}, S, f_0, \pi)  \to M$ is in $V[G]$. 

We first claim that $S$ and $j[S]$ are in $M$. Since $\gamma$ is countable in $V[G]$, there is a real $x$ coding $S$ in $V[G]$. But by Lemma~\ref{key generic embedding}, such an $x$ is in $M$. Hence $S$ is also in $M$. Since $\gamma$ is countable in $V[G]$, there is a $b \in \wp_{\omega_1}(\mathbb{R}^V)$ in $V[G]$ such that $\pi[b] = S$ and hence $j(\pi)[b] = j(\pi) \bigl[j[b]\bigr] = j\bigl[ \pi [b] \bigr] = j[S]$ in $V[G]$. But since $j(\pi) \in M$ and $b \in M$ by Lemma~\ref{key generic embedding}, $j[S] = j(\pi)[b]$ is also in $M$, as desired. 

By Lemma~\ref{stronger code} and elementarity of $j$, the following is true in $M$: For any $a \in \wp_{\omega_1}(\R)$ closed under $j(f_0)$, $j(S) {\upharpoonright} j(\pi)[a]$ is an $\infty$-Borel code and $B_{j(S){\upharpoonright}j(\pi) [a]} \subseteq B_{j(S)}$. Also, by elementarity of $j$, $\exists^{\mathbb{R}}B_{j(S)}$ is a well-founded relation on a set of reals in $M$.

We now argue that $B_{j[S]} \subseteq B_{j(S)}$ in $M$. 
Let $b \in \wp_{\omega_1}(\R^V)$ be such that $\pi [b] = S$ in $M$. Let $a \in \wp_{\omega_1}(\R^V)$ be the closure of $b$ under $f_0$ in $V[G]$. Then since $a$ is countable in $V[G]$, by Lemma~\ref{key generic embedding}, $a$ is also in $M$.  Then $j(f_0) [a] = j(f_0) \bigl[ j[a] \big] =j \bigl[ f_0 [a] \bigr] \subseteq j[a] = a$. So $a$ is closed under $j(f_0)$ in $M$. 
Hence $j(S) {\upharpoonright} j(\pi) [a]$ is an $\infty$-Borel code and $B_{j(S){\upharpoonright} j(\pi) [a]} \subseteq B_{j(S)}$ in $M$. Since $ j[S] =  j\bigl[ \pi [b] \bigr] \subseteq j \bigl[ \pi [a] \bigr] = j(\pi) \bigl[ j[a] \bigr] =  j(\pi) [a] $ while $j (\pi) [a] = j \bigl[ \pi [a] \bigr] \subseteq j[ \pi [\R^V] ] = j [\gamma]$, we have $ j[S] \subseteq j (\pi)[a] \subseteq j[\gamma]$. Hence $j(S) \upharpoonright j[S] \subseteq j(S) \upharpoonright j(\pi) [a] \subseteq j(S) \upharpoonright j[\gamma]$. However, $j(S) \upharpoonright j[S] = j(S) \upharpoonright j[\gamma]$. Therefore, $j(S) \upharpoonright j[S]= j(S) \upharpoonright j(\pi) [a] $. Since $j (S) \upharpoonright j[S]$ is the same as $j[S]$ as $\infty$-Borel codes, $B_{j[S]} = B_{j(S){\upharpoonright} j(\pi) [a]} \subseteq B_{j(S)}$ in $M$, as desired. 

Hence $\exists^{\mathbb{R}} B_{j[S]}$ is also a wellfounded relation on a set of reals in $M$. Since $j[S]$ is countable in $M$, the relation $\exists^{\mathbb{R}} B_{j[S]}$ is $\undertilde{\mathbf{\Sigma}}^1_1$ and hence by Kunen-Martin Theorem (see \cite[2G.2]{new_Moschovakis}), its rank is less than $\omega_1$ in $M$ which is the same as $\gamma^+$ in $V$. Finally, since $S$ and $j[S]$ are equivalent as $\infty$-Borel codes, $\exists^{\mathbb{R}} B_S$ has length less than $\omega_1$ in $M$. Since $M$ has more reals than $V$, $\bigl(\exists^{\mathbb{R}} B_S \bigr)^V \subseteq \bigl(\exists^{\mathbb{R}} B_S \bigr)^M$. Therefore, the length of $\bigl(\exists^{\mathbb{R}} B_S \bigr)^V$ is less than $\omega_1^M = (\gamma^+)^V$, as desired. 

This completes the proof of Theorem~\ref{Becker's bound}. 
\end{proof}

Becker proved the following:
\begin{Thm}[Becker]
Assume $\AD$, $\DC$, and the uniformization for every relation on the reals. Suppose that the conclusion of Theorem~\ref{Becker's bound} holds, i.e., let $A$ be a subset of $\mathbb{R}^3$ which has a strong $\infty$-Borel code $S$ where the tree of $S$ is on $\gamma < \Theta$. If $\exists^{\mathbb{R}} A$ is a strict well-founded relation on a set of reals, then the length of $\exists^{\mathbb{R}} A$ is less than $\gamma^+$. 
Then every set of reals is Suslin. 
\end{Thm}
\begin{proof}
See \cite[2.3~Main~Theorem]{Becker}.
\end{proof}

We simulate Becker's argument using $\BlADR$ to show that every set of reals is Suslin. 
The exposition of our arguments is closer to the one in Larson~\cite[Section~12.3]{ADplus} than to the original one by Becker~\cite{Becker}.

For the rest of this section, we assume $\BlADR$ and $\DC$. We fix a set of reals $A$ and prove that $A$ is Suslin. 

We start with constructing an $\omega$-sequence of \lq good' pointclasses. 
Recall the definition of Spector pointclasses from Definition~\ref{def:Spector}. 
\begin{Claim}\label{sequence of ptclass}
There is a sequence $\bigl( (\Lambda_n, <_n, \lambda_n, ) \mid n < \omega \bigr)$ such that for all $n$, 
\begin{enumerate}
\item $\Lambda_n$ is a Spector pointclass closed under $\exists^{\mathbb{R}}$ and $\forall^{\mathbb{R}}$, $\Lambda_n \subseteq \Lambda_{n+1}$, and $A \in \Lambda_0$, 

\item every relation on the reals which is projective in a set in $\undertilde{\mathbf{\Lambda}}_n$ can be uniformized by a function in $\undertilde{\mathbf{\Lambda}}_{n+1}$, and 

\item $<_n$ is in $\undertilde{\mathbf{\Lambda}}_n$ and a strict wellfounded relation on the reals of length $\lambda_n$ and every set of reals which is projective in a set in $\undertilde{\mathbf{\Lambda}}_n$ has a strong $\infty$-Borel code whose tree is on $\lambda_{n+1}$. 
\end{enumerate}
\end{Claim}

\begin{proof}[Proof of Claim~\ref{sequence of ptclass}]
We construct them by induction on $n$. For $n=0$, let $\Lambda_0$ be any Spector pointclass closed under $\exists^{\mathbb{R}}$ and $\forall^{\mathbb{R}}$ containing $A$ which exists by Theorem~\ref{inductive}, and $<_0$ be any strict wellfounded relation on the reals in $\Lambda_0$. Set $\lambda_0$ to be the length of $<_0$. Then they satisfy all the items above.

Suppose we have constructed $(\Lambda_n, <_n, \lambda_n)$ with the above properties. We construct $\Lambda_{n+1}, <_{n+1}$, and $\lambda_{n+1}$. First note that there is a set of reals $B_n$ which is not projective in any set in $\undertilde{\mathbf{\Lambda}}_n$ by uniformization for every relation on the reals. Then by the weak version of Wadge Lemma (Lemma~\ref{weak Wadge}), every set projective in a set in $\undertilde{\mathbf{\Lambda}}_n$ is $\undertilde{\mathbf{\Sigma}}^1_2$ in $B_n$. Let $H_n^1$ and $H_n^2$ be universal sets for $\undertilde{\mathbf{\Sigma}}^1_2 (B_n)$ sets of reals and $\undertilde{\mathbf{\Sigma}}^1_2 (B_n)$ subsets of $\mathbb{R}^2$, respectively. By uniformization for every relation on the reals, there is a function $f_n$ uniformizing $H_n^2$. By Theorem~\ref{strong infty}, there is a $\gamma < \Theta$ such that $H_n^1$ has a strong $\infty$-code whose tree is on $\gamma$. Let $\lambda_{n+1} = \gamma$, $<_{n+1}$ be a strict wellfounded relation on the reals of length $\lambda_{n+1}$, and let $\Lambda_{n+1}$ be a Spector pointclass closed under $\exists^{\mathbb{R}}$ and $\forall^{\mathbb{R}}$ containing $\Lambda_n \cup \{ H_n^1, H_n^2, f_n , <_{n+1}\}$. We show that they satisfy all the items above for $n+1$. The first item is trivial. The second item is easy by noting that if $f_n$ uniformizes $H_n^2$ then $(f_n)_a$ uniformizes $(H_n^2)_a$ for any real $a$. The third item follows from that if $H_n^1$ has a strong $\infty$-code whose tree is on $\lambda_{n+1}$, then $(H_n^1)_a$ has a strong $\infty$-code whose tree is on $\lambda_{n+1}$ for every real $a$. 
\renewcommand{\qedsymbol}{$\square \  (\text{Claim~\ref{sequence of ptclass}})$}
\end{proof}
Note that in the proof of Claim~\ref{sequence of ptclass}, we have used $\DC$ to obtain the sequence $\bigl( (\Lambda_n, <_n, \lambda_n, ) \mid n < \omega \bigr)$.

We will introduce an integer-real game $\mathcal{G}_{\mathrm{Be}} (A)$ for the verification of Suslinness of $A$. 
We first prepare some notations. 
We fix a simple surjection $\rho$ from the reals to $\{ 0, 1\}$, e.g., $x \mapsto x(0)$. 
We fix $\bigl( (\Lambda_n, <_n, \lambda_n) \mid n < \omega \bigr)$ as above and let $\Gamma^{\mathrm{I}}_n = \Lambda_{2n}, \Gamma^{\mathrm{II}}_n = \Lambda_{2n+1}, <^{\mathrm{I}}_{n} $ be the prewellordering on the reals induced by $\rho$, $<^{\mathrm{II}}_n = <_{2n+1}$, $\gamma^{\mathrm{I}}_n = 2$, and  $\gamma^{\mathrm{II}}_n = \lambda_{2n+1}$. Let $\rho^{\mathrm{I}}_n = \rho$ and $\rho^{\mathrm{II}}_n$ be the surjection from the reals onto $(\gamma^{\mathrm{II}}_n)^{n+1}$ induced by $<_{2n+1}$. Let $\pi^{\mathrm{II}}_n$ be the function $a \mapsto \rho^{\mathrm{II}}_n [ G^n_a]$ where $G^n$ is a universal set for $\undertilde{\mathbf{\Gamma}}^{\mathrm{II}}_n$ sets of reals (we do not use $\pi^{\mathrm{I}}_n$). Then by the weak version of Moschovakis' Coding Lemma (Theorem~\ref{weak coding}), $\pi^{\mathrm{II}}_n$ is a surjection from the reals onto $\wp \bigl( (\gamma^{\mathrm{II}}_n)^{n+1} \bigr)$. Consider the following game $\mathcal{G}_{\mathrm{Be}} (A)$: Player I plays $0$ or $1$ and player II plays a real one by one in turn and they produce a real $z$ and a sequence $t \in \R^{\omega}$, respectively. Setting $T_n = \pi^{\mathrm{II}}_n \bigl( t(n) \bigr)$, which is a subset of $(\gamma^{\mathrm{II}}_n)^{n+1}$, player II wins in the game $\mathcal{G}_{\mathrm{Be}} (A)$ if for all $n \in \omega$, $T_{n+1} {\upharpoonright} (n +1) \subseteq T_n $   
and $z \in A \iff \text{the tree } T = \{ \emptyset \} \cup \bigcup_{n \in \omega} \bigcap_{m > n} T_m \upharpoonright (n+1)$ is ill-founded, 
where $T_{m} {\upharpoonright} (n+1) = \{ s  {\upharpoonright} (n+1) \mid s \in T_m\}$. This is an integer-real game in the sense that player I chooses integers and player II chooses reals. 

We will introduce an integer-integer game $\tilde{\mathcal{G}}_{\mathrm{Be}} (A)$ simulating the game $\mathcal{G}_{\mathrm{Be}} (A)$. For each perfect set of reals $P$, we fix a homeomorphism $c_{P} \colon P \to (2^{\omega})^3$ which is recursive in $P$. We say a real $y$ is {\it identified with a triple of reals $(u, x' ,y')$ via $P$} if $c_P (y) = (u, x',y')$.  In the game $\tilde{\mathcal{G}}_{\mathrm{Be}} (A)$, players choose pairs of $0$ or $1$ one by one and produce pairs of reals $(x_0, y_0)$ and $(a_0, b_0)$ in $\omega$ rounds respectively. From $(x_0, y_0)$ and $(a_0, b_0)$, we \lq \lq decode" a real $ z$ and an $\omega$-sequence of reals $t$ respectively as follows: For each pointclass $\Gamma^{\mathrm{I}}_n$ (resp. $\Gamma^{\mathrm{II}}_n$) above, we fix a set $U^{\mathrm{I}, n}$ in $\undertilde{\mathbf{\Gamma}}^{\mathrm{I}}_n$ (resp. $U^{\mathrm{II},n}$ in $\undertilde{\mathbf{\Gamma}}^{\mathrm{II}}_n$)  which is universal for relations in $\undertilde{\mathbf{\Gamma}}^{\mathrm{I}}_n$ (resp. $\undertilde{\mathbf{\Gamma}}^{\mathrm{II}}_n$). Setting $F_0 = U^{\mathrm{I}, 0}_{x_0}$, $F_0$ is a function from the reals to perfect sets of reals (otherwise player I loses). Let $P_{x_0} = F_0(x_0)$. Then $y_0$ is an element of $P_{x_0}$ (otherwise player I loses) and is identified with a triple $(u_0, x_1, y_1)$ of reals via $P_{x_0}$. Then setting $F_1 = U^{\mathrm{I}, 1}_{x_1}$, $F_1$ is a function from the reals to perfect sets of reals (otherwise player I loses). Let $P_{x_1} = F_1(x_1)$. Then $y_1$ is an element of $P_{x_1}$ (otherwise player I loses) and is identified with a triple of reals $(u_1, x_2, y_2)$ via $P_{x_1}$. Continuing this process, from the pair $(x_0, y_0)$, one can unwrap $(x_n , y_n)$ and obtain $(u_n , x_{n+1}, y_{n+1})$ for each $n$ and get an $\omega$-sequence $(u_n \mid n < \omega)$. Let $z (n) = \rho (u_n)$. In the same way, one can obtain an $\omega$-sequence $(t_n \mid n< \omega)$ of reals from $(a_0, b_0)$. Setting $T_n = \pi^{\mathrm{II}}_n ( t_n )$, player II wins in the game $\tilde{\mathcal{G}}_{\mathrm{Be}} (A)$ if for all $n \in \omega$, $T_{n+1} {\upharpoonright} (n+1) \subseteq T_n $  
and $z \in A \iff 
\text{the tree } T = \{ \emptyset \} \cup \bigcup_{n \in \omega} \bigcap_{m > n} T_m \upharpoonright (n+1)$ is ill-founded. 

Becker proved the following:
\begin{Lem}[Becker]\label{Becker lemma}
${}$

\begin{enumerate}
\item If player I has a winning strategy in the game $\tilde{\mathcal{G}}_{\mathrm{Be}} (A)$, then player I has a winning strategy $\sigma$ in the game $\mathcal{G}_{\mathrm{Be}} (A)$ such that $\sigma$ is a countable union of sets in $\bigcup_{n < \omega} \undertilde{\mathbf{\Gamma}}^{\mathrm{II}}_n$. 

\item If player II has a winning strategy in the game $\tilde{\mathcal{G}}_{\mathrm{Be}} (A)$, then player II has a winning strategy in the game $\mathcal{G}_{\mathrm{Be}} (A)$. 
\end{enumerate}
\end{Lem}

\begin{proof}
See \cite[Lemma~A \& B]{Becker} or \cite[Section~12.3]{ADplus}.  
\end{proof}

We will show the following two lemmas corresponding to the above Becker's lemma (Lemma~\ref{Becker lemma}): First recall the notation $\tau_y$ for a blindfolded strategy in \S\,\ref{sec:Blackwell games}. Let $ B \subseteq \R^{\omega}$. A mixed strategy $\sigma$ for player I is {\it weakly optimal in $B$} if for any $s \in \mathbb{R}^{\mathrm{Even}}$, the set $\{ x \mid \sigma (s) (x) \neq 0\}$ is finite and for any $\omega$-sequence of reals $y $, $\mu_{\sigma, \tau_y} (B) > 1/2$. One can introduce the weak optimality for mixed strategies for player II in the same way. Note that if player I has an optimal strategy in some payoff set, then player I has a weakly optimal strategy in the same payoff set. The same holds for player II. 
\begin{Lem}\label{lem A}
If player I has an optimal strategy $\tilde{\sigma}$ in the game $\tilde{\mathcal{G}}_{\mathrm{Be}} (A)$, then player I has a weakly optimal strategy $\sigma$ in the game $\mathcal{G}_{\mathrm{Be}} (A)$ such that $\sigma$ is a countable union of sets in $\bigcup_{n < \omega} \undertilde{\mathbf{\Gamma}}^{\mathrm{II}}_n$. 
\end{Lem}

\begin{Lem}\label{small conj}
If player II has an optimal strategy in the game $\tilde{\mathcal{G}}_{\mathrm{Be}} (A)$, then player II has a weakly optimal strategy in the game $\mathcal{G}_{\mathrm{Be}} (A)$. 
\end{Lem}

Before seeing the proofs of Lemma~\ref{lem A} and Lemma~\ref{small conj}, we show that Lemma~\ref{lem A} and Lemma~\ref{small conj} imply that $A$ is Suslin (and hence Theorem~\ref{almost theorem}).

Using Lemma~\ref{lem A}, Corollary~\ref{strong code countable intersection}, and Theorem~\ref{Becker's bound}, we first argue that there is no optimal strategy for player I in the game $\tilde{\mathcal{G}}_{\mathrm{Be}} (A)$:
\begin{Lem}\label{no I}
There is no optimal strategy for player I in the game $\tilde{\mathcal{G}}_{\mathrm{Be}} (A)$. 
\end{Lem}

\begin{proof}
To derive a contradiction, suppose that player I has an optimal strategy in the game $\tilde{\mathcal{G}}_{\mathrm{Be}} (A)$. Then by Lemma~\ref{lem A}, player I has a weakly optimal strategy $\sigma$ in the game $\mathcal{G}_{\mathrm{Be}} (A)$ such that $\sigma$ is a countable union of sets in $\bigcup_{n < \omega}\undertilde{\mathbf{\Gamma}}^{\mathrm{II}}_n$. 

Consider the following set:
\begin{align*}
X = \{ (t, s) \in \R^{\omega} \times \mathbb{R}^{<\omega} \mid & \ \mu_{\sigma, \tau_t} \bigl(\{ (z, t') \mid t' = t \text{ and } z \in A\}\bigr) > 1/2 \text{ and}\\
& \bigl(\forall i < \mathrm{lh} (s) \bigr) \ \bigl( |s(0)|_{<^{\mathrm{II}}_0}, \ldots , |s(i)|_{<^{\mathrm{II}}_i} \bigr) \in T_i \}, 
\end{align*}
where $|s(i)|_{<^{\mathrm{II}}_i} $ is the rank of $s(i)$ with respect to the wellfounded relation $<^{\mathrm{II}}_i$  on $\R$ and $T_i = \rho^{\mathrm{II}}_i \bigl(t(i) \bigr)$. For $(t, s)$ and $(t', s')$ in $X$, $(t,s) < (t' , s')$ if $t =t'$ and $s$ codes a node in the tree $T = \{ \emptyset \} \cup \bigcup_{n \in \omega} \bigcap_{m > n} T_m \upharpoonright (n+1)$ extending a node coded by $s'$. Note that for any $(t,s)$ in $X$, if $T$ is the tree coded by $t$ as above, $T$ is wellfounded because $\sigma$ is weakly optimal for player I in the game $\mathcal{G}_{\mathrm{Be}} (A)$. Hence $<$ is a strict wellfounded relation on $X$. 
Let $\gamma_{\omega} = \mathrm{sup} \{ \gamma^{\mathrm{II}}_n \mid n \in \omega \}$. By $\DC$, the cofinality of $\Theta$ is greater than $\omega$. Hence $\gamma_{\omega} < \Theta$. 
Note that for any ordinal $\alpha < \gamma_{\omega}^+$, there is a wellfounded tree $T$ coded by some real $t$ as above such that the rank of $T$ is $\alpha$. Hence the length of $(X,<)$ is at least $\gamma_{\omega}^+$. 

Since $\sigma$ is a countable union of sets in $\bigcup_{n < \omega} \undertilde{\mathbf{\Gamma}}^{\mathrm{II}}_n$, the relation $<$ on $X$ is in the pointclass $\exists^{\mathbb{R}} \bigwedge^{\omega}  \bigvee^{\omega} \bigcup_{n \in \omega} \undertilde{\mathbf{\Gamma}}^{\mathrm{II}}_n$, i.e., it is a projection of a countable intersection of a countable union of sets in $\bigcup_{n < \omega} \undertilde{\mathbf{\Gamma}}^{\mathrm{II}}_n$. Since every set in $\undertilde{\mathbf{\Gamma}}^{\mathrm{II}}_n$ has a strong $\infty$-Borel code whose tree is on $\gamma^{\mathrm{II}}_n$ for every $n$, by Corollary~\ref{strong code countable intersection}, every set in $\bigwedge^{\omega} \bigvee^{\omega} \bigcup_{n \in \omega} \undertilde{\mathbf{\Gamma}}^{\mathrm{II}}_n$ has a strong $\infty$-Borel code whose tree is on $\gamma_{\omega}$. By Theorem~\ref{Becker's bound}, the length of $<$ on $X$ must be less than $\gamma_{\omega}^+$, which is not possible because it was at least $\gamma_{\omega}^+$. Contradiction!

This finishes the proof of Lemma~\ref{no I}. 
\end{proof}

We now show that $A$ is Suslin assuming Lemma~\ref{lem A} and Lemma~\ref{small conj}: 
\begin{proof}[Proof of Theorem~\ref{almost theorem} from Lemma~\ref{lem A} and Lemma~\ref{small conj}]

By Lemma~\ref{no I}, player I does not have an optimal strategy in the game $\tilde{\mathcal{G}}_{\mathrm{Be}} (A)$. Hence by $\BlAD$, player II has an optimal strategy in the game $\tilde{\mathcal{G}}_{\mathrm{Be}} (A)$. By Lemma~\ref{small conj}, player II has a weakly optimal strategy $\tau$ in the game $\mathcal{G}_{\mathrm{Be}} (A)$. Note that $\tau$ can be seen as a real because each measure on the reals given by $\tau$ is of finite support by the weak optimality of $\tau$. For each finite binary sequence $s$ of length $n$, let $$T_s = \left\{ u \in \R^n \middle| (\forall i < n ) \ \tau \Bigl( (s {\upharpoonright} i) * \bigl(u {\upharpoonright} (i -1) \bigr) \Bigr) \bigl( u(i) \bigr) \neq 0 \right\},$$ where $(s {\upharpoonright} i) * \bigl(u {\upharpoonright} (i -1) \bigr)$ is the concatenation of $s {\upharpoonright} i$ and $u {\upharpoonright} (i -1)$ bit by bit, i.e., $(s {\upharpoonright} i) * \bigl(u {\upharpoonright} (i -1) \bigr) = \bigl(s(0), u(0), s(1), u(1), \ldots , s(i-2), u(i-2), s(i-1) \bigr)$. 
Since $\tau$ is of finite support, for each binary sequence $s$, the set $a_s = \bigcup \{ \mathrm{ran} (u) \mid u \in T_s\}$ is a finite set of reals. Hence $a = \bigcup \{ a_s \mid s \in 2^{<\omega} \}$ is a countable set of reals. Let $\iota \colon a \to \omega$ be an injection. 
For a binary sequence $s$ of length $n$, let $\pi_s \colon T_s \to \omega^n$ be such that $\pi_s (u) = \iota [u]$.  
For any real $x$, $T_x = \bigcup_{n \in \omega} T_{x{\upharpoonright}n}$ is a tree on $a$. Let $\pi_x (T_x) = \bigcup_{n < \omega} \pi_{x \upharpoonright n}[T_{x \upharpoonright n}]$. Then $\pi_x (T_x)$ is a tree on $\omega$, and the map $\pi_x$ induces a homeomorphism $(\pi_x)_{\ast} \colon [T_x] \to [\pi_x (T_x)]$. 
Consider the following tree:
\begin{align*}
W = \left\{ (s,t,v) \in \bigcup_{n\in \omega} \bigl(2^n \times \omega^n \times (\gamma_{\omega})^n \bigr)\middle| \ t \in \pi_s (T_s) \text{ and } \bigl(\forall i < \mathrm{lh}(s) \bigr)\ v(i) = |y_i|_{<^{\mathrm{II}}_i} \right\},
\end{align*}
where $y_i = \iota^{-1} \bigl( t(i) \bigr)$. 
Let $B = \{ (x,y) \in \mathbb{R} \times \omega^{\omega} \mid \bigl(\exists f \in (\gamma_{\omega})^{\omega} \bigr)\ (x,y,f) \in [W]\}$.  
Then by the weak optimality of $\tau$, the following holds: for all reals $x$,
\begin{align*}
x \in A \iff &\mu_{\sigma_x, \tau} \bigl( (\pi_x)_{\ast}^{-1} (B_x) \bigr) > 1/2\\
\iff & (\exists T': \text{a tree on $\omega$}) \ [T'] \subseteq B_x \text{ and } \mu_{\sigma_x, \tau} \bigl( (\pi_x)_{\ast}^{-1} ([T']) \bigr) > 1/2.
\end{align*}

Hence $A$ is definable from $B$ using the quantifiers $\exists^{\R}$ and $\forall^{\R}$, and coutable unions and intersections with projective sets. Since $B$ is Suslin, by Corollary~\ref{cor:Suslin}, $A$ is also Suslin, as desired. 

Since $A$ was an arbitrary set of reals, we have shown that every set of reals is Suslin. Then by Theorem~\ref{Suslin-co-Suslin-det}, $\AD$ holds. Now by Theorem~\ref{equiv AD and ADR} and Theorem~\ref{uniformization}, $\ADR$ holds. 

This completes the proof of Theorem~\ref{almost theorem} from Lemma~\ref{lem A} and Lemma~\ref{small conj}.
\end{proof}

We now prove Lemma~\ref{lem A} and Lemma~\ref{small conj}: 
\begin{proof}[Proof of Lemma~\ref{lem A}]
We first topologize $\mathrm{Prob}(\mathbb{R})$, the set  of all Borel probabilities on the reals. Consider the following map $\iota \colon \mathrm{Prob}(\mathbb{R}) \to [0,1]^{(2^{<\omega})}$: Given a Borel probability $\mu$ on the reals, for each finite binary sequence $s$, $\iota (\mu) (s) = \mu ([s])$. We topologize $[0,1]^{(2^{<\omega})}$ by the product topology where each coordinate $[0,1]$ is equipped with the relative topology of the real line. We then identify $\mathrm{Prob}(\mathbb{R})$ with its image via $\iota$ and topologize it with the relative topology of $[0,1]^{(2^{<\omega})}$. Then the space $\mathrm{Prob}(\mathbb{R})$ is compact. 

Notice that since every set of reals is Lebesgue measurable, by Theorem~\ref{fact:LM}, for all $\mu \in \text{Prob}(\R)$, every set of reals is $\mu$-measurable. 
\begin{Claim}\label{easy measure value}
For any set of reals $B$, the map $\mu \mapsto \mu (B)$ is a continuous map from $\mathrm{Prob}(\mathbb{R}) $ to $[0,1]$.
\end{Claim}
\begin{proof}[Proof of Claim~\ref{easy measure value}]
This is easy when $B$ is closed or open. In general, it follows from the following equations: For any $\mu \in \mathrm{Prob}(\mathbb{R})$, 
\begin{align*}
\mu (B) = & \ \mathrm{sup} \{ \mu (C) \mid  C\subseteq B \text{ and $C $ is closed}\}\\
= & \ \mathrm{inf} \{ \mu (O) \mid  O \supseteq B \text{ and $O$ is open}\}.
\end{align*}
\renewcommand{\qedsymbol}{$\square \ (\text{Claim~\ref{easy measure value}})$}
\end{proof}

Next, we introduce a complete metric $d$ on $\mathrm{Prob}(\mathbb{R})$ compatible with the topology we consider. Let $(s_n \mid n \in \omega)$ be an injective enumeration of finite binary sequences. For $\mu $ and $\mu'$ in $\mathrm{Prob}(\mathbb{R})$, $d(\mu , \mu') = \sum_{n\in \omega} | \mu ([s_n]) - \mu' ([s_n])| / 2^{n+1}$. Then $d$ is a complete metric compatible with our topology. Since $\mathrm{Prob}(\mathbb{R})$ is compact, the map $\mu \mapsto \mu (A)$ is uniformly continuous with the metric $d$. Hence there is an $\varepsilon > 0$ such that if $d (\mu , \mu') < \varepsilon $, then $|\mu (A) - \mu' (A) | < 1/2$. Let us fix a sequence $( \varepsilon_n \mid n \in \omega )$ of positive real numbers such that $\sum_{n \in \omega} \varepsilon_n / 2^{n+1}  < \varepsilon$. For any finite binary sequence $v$, let $k_{v} $ be the natural number such that $s_{k_v} = v$. 

Let $\tilde{\sigma}$ be an optimal strategy for player I in the game $\tilde{\mathcal{G}}_{\mathrm{Be}} (A)$. 
We will construct a weakly optimal strategy $\sigma$ for player I in the game $\mathcal{G}_{\mathrm{Be}} (A)$ such that $\sigma$ is a countable union of sets in $\bigcup_{n < \omega} \undertilde{\mathbf{\Gamma}}^{\mathrm{II}}_n$. 

We first prove the following claim which is corresponding to \cite[Lemma~12.6]{ADplus}: 
\begin{Claim}\label{claim:BeckerA}
There is a function $E$ mapping each $(t_0, \ldots , t_{n-1}) \in \R^{n}$ for each $n \in \omega$ to a sequence of reals $E(t_0, \ldots t_{n-1}) = (a_0, a_1, \ldots , a_{n})$ such that setting $H_m = U^{\mathrm{II}, m}_{a_m}$ and $Q_{m} = H_m (a_m)$ for each $m \le n$, the following hold: 
\begin{enumerate}[(1)]
\item\label{BeckerA1} For all $m \le n$, $Q_{m}$ is a perfect set of reals, and for all $b \in Q_{m}$, there is a unique sequence $(b_0, \ldots , b_{m-1}) \in \R^m$ such that for all $\ell < m-1$, $b_{\ell}$ is identified with $(t_{\ell} , a_{\ell + 1}, b_{\ell+1})$ via $Q_{\ell}$, and $b_{m-1}$ is identified with $(t_{m-1}, a_m, b)$ via $Q_{m-1}$. 
We call the unique sequence $(b_0, \ldots , b_{m-1})$ for $b$ the {\it induced sequence from $b$ via $(Q_{\ell} \mid \ell < m)$}. 

\item\label{BeckerA2} For all $m \le n$ and $s \in 2^m$, let $h^s \colon Q_m \to [0,1]^2$ be the following function: For each $b\in \R$ and $i \in 2$, $h^s (b) (i) = \mu_{\tilde{\sigma}, \tau_{(a_0, b_0)}} \bigl(\{ \bigl( (x_0, y_0) , (a'_0, b'_0) \big) \mid (a'_0, b'_0) = (a_0, b_0) \text{ and for all $\ell < m$}, \rho (u_{\ell}) = s(\ell) \text{ and } \rho (u_m) = i \} \bigr)$, where $(b_0, \ldots , b_{m-1})$ is the induced sequence from $b$ via $(Q_{\ell} \mid \ell < m)$, and for each $\ell \le m$, setting $F_{\ell} = U^{\mathrm{I},\ell}_{x_{\ell}}$ and $P_{\ell} = F_{\ell} (x_{\ell})$, the sequence $(y_0, \ldots , y_{m})$ is the induced sequence from $y_{m+1}$ via $(P_{\ell} \mid \ell \le m)$, and each $y_{\ell}$ is identified with $(u_{\ell}, x_{\ell +1}, y_{\ell +1})$ via $P_{\ell}$.   

Then for all $m \le n$, $s\in 2^m$, $i < 2$, and $b, b' \in Q_m$, $|h^s (b) (i) - h^s (b') (i)| < \min \{ \varepsilon_{k_0} , \varepsilon_{k_1} \}$ for $k_i = k_{s^{\frown} \langle i \rangle}$. 
\end{enumerate}
\end{Claim}

\begin{proof}[Proof of Claim~\ref{claim:BeckerA}]
We construct $E(t_0, \ldots t_{n-1}) = (a_0, a_1, \ldots , a_{n})$ by induction on $n$. 

Let $n=0$. We will construct $E( \emptyset ) = (a_0)$. 
Fix a real $a$. Consider the function $h^{\emptyset}_a \colon \mathbb{R} \to [0,1]^2$ as follows: Given a real $b$ and $i < 2$, let $h^{\emptyset}_a (b) (i)  = \mu_{\tilde{\sigma}, \tau_{(a,b)}} \bigl(\{ \bigl( (x_0, y_0) ,\\  (a'_0, b'_0)\bigr) \mid (a'_0, b'_0) = (a,b) \text{ and }\rho (u_0) = i\} \bigr) $, where setting $F_{0} = U^{\mathrm{I},0}_{x_0}$ and $P_{0} = F_{0} (x_0)$, $y_0 $ is identified with $(u_0, x_1, y_1)$ via $P_0$. Since every set of reals has the Baire property, by Theorem~\ref{fact:category}, the function $h^{\emptyset}_a$ is continuous on a comeager set. Then there is a perfect set of reals $Q$ such that for all $b , b'$ in $Q$ and $i < 2$, $|h^{\emptyset}_a (b) (i)  - h^{\emptyset}_a (b') (i) | < \min \{ \varepsilon_{k_0}, \varepsilon_{k_1} \}$ for $k_i = k_{v_i}$, where $v_i$ is the sequence $( i )$ of length $1$. Since the set $X_0 = \{ (a, Q )\mid (\forall b , b' \in Q)\ (\forall i < 2)\,  |h^{\emptyset}_a (b)  (i)  - h^{\emptyset}_a (b')  (i) | < \min \{ \varepsilon_{k_{v_0}} , \varepsilon_{k_{v_1}} \}\}$ is projective in $\Gamma^{\mathrm{I}}_0$, there is a real $a_0$ such that the function $H_0 = U^{\mathrm{II},0}_{a_0}$ uniformizes $X_0$. 

Set $E(\emptyset ) = (a_0)$. Then letting $Q_0 = H_0 (a_0)$, $E(\emptyset)$ satisfies \ref{BeckerA1} and \ref{BeckerA2} in the claim for $m = n=0$ and $s = \emptyset$. 

Now let $n > 0$. Take any $(t_0, \ldots , t_{n-1}) \in 2^n$. We will construct $E(t_0, \ldots , t_{n-1})$. By induction hypothesis, we already obtained $E(t_0, \ldots , t_{n-2}) = (a_0 , \ldots , a_{n-1})$ and the sequence  $(Q_m \mid m \le n-1)$  that satisfy \ref{BeckerA1} and \ref{BeckerA2} in the claim for all $m \le n-1$ and $s \in 2^{m}$. 
We will find a real $a_n$ such that the sequences $(a_0 , \ldots , a_{n})$ and $(Q_m \mid m \le n)$ satisfy \ref{BeckerA1} and \ref{BeckerA2} in the claim for $m =n$ and all $s \in 2^{m}$. 

Fix a real $a$. Then fix an $s \in 2^n$. Consider the function $h^s_a \colon \mathbb{R} \to [0,1]^2 $ as follows: For a real $b$ and $i < 2$, set $h^s_a (b) (i) = \mu_{\tilde{\sigma}, \tau_{( a_0, b_0)}} \bigl( \{ \bigl( (x_0,y_0) , (a'_0 , b'_0)\bigr)\mid (a'_0, b'_0) = ( a_0, b_0) \text{ and for all $\ell < n$}, \rho (u_{\ell}) = s(\ell) \text{ and } \rho (u_n) = i\} \bigr)$, where $(b_0, \ldots , b_{n-1})$ is the induced sequence from $b$ via $(Q_{\ell} \mid \ell < n)$, and for each $\ell \le n$, setting $F_{\ell} = U^{\mathrm{I},\ell}_{x_{\ell}}$ and $P_{\ell} = F_{\ell} (x_{\ell})$, the sequence $(y_0, \ldots , y_{n})$ is the induced sequence from $y_{n+1}$ via $(P_{\ell} \mid \ell \le n)$, and each $y_{\ell}$ is identified with $(u_{\ell}, x_{\ell +1}, y_{\ell +1})$ via $P_{\ell}$.  

Since every set of reals has the Baire property,  by Theorem~\ref{fact:category}, for all $s\in 2^n$, the function $h^s_a$ is continuous on a comeager set. Hence there is a perfect set of reals $Q$ such that for all $s \in 2^n$, $b, b' \in Q$, and $i < 2$, $|h^s_a (b) (i) - h^s_a(b') (i)| < \min \{ \varepsilon_{k_0}, \varepsilon_{k_1}\}$ for $k_i = k_{v_i}$, where $v_i = s^{\frown} \langle i \rangle$. 
Since the set $X_n = \{ (a, Q )\mid (\forall s \in 2^n)\  (\forall b , b' \in Q)\ (\forall i < 2)\,  |H^s_a (b)  (i)  - H^s_a (b')  (i) | < \min \{ \varepsilon_{k_{v_0}}, \varepsilon_{k_{v_1}}\} \text{ for }v_i = s^{\frown} \langle i \rangle\}$ is projective in $\Gamma^{\mathrm{I}}_n$, there is a real $a_n$ such that the function $H_n = U^{\mathrm{II},n}_{a_n}$ uniformizes $X_n$.  

Set $E(t_0, \ldots , t_{n-1} ) = (a_0, \ldots , a_n)$. Then letting $Q_n = H_n (a_n)$, $E(t_0, \ldots , t_{n-1})$ satisfies \ref{BeckerA1} and \ref{BeckerA2} in the claim for $m = n$ and all $s \in 2^m$ as desired. 

This completes the proof of Claim~\ref{claim:BeckerA}. 
\end{proof}

We can now define a mixed strategy $\sigma$ for player I in the game $\mathcal{G}_{\mathrm{Be}} (A)$ as follows: Let $n < \omega$, $s\in 2^n$, and $t = (t_0, \ldots , t_{n-1}) \in \R^n$. Let $s * t$ be the concatenation of $s$ and $t$ bit by bit, i.e., $s*t = (s_0, t_0, s_1, t_1, \ldots , s_{n-1}, t_{n-1})$. We define $\sigma (s * t) \in \text{Prob}(2)$ by induction on $n$. 
Suppose that we already obtained $\sigma \bigl(s \upharpoonright m) * (t\upharpoonright m)\bigr)$ for all $m < n$. Then let $\sigma (s*t)$ be as follows: Let $E(t_0, \ldots , t_{n-1}) = (a_0, \ldots , a_n)$, $(Q_m \mid m \le n)$, and $h^s \colon Q_n \to [0,1]^2$ be as in Claim~\ref{claim:BeckerA}. Then set 
\begin{align*}
\sigma (s * t) (0) & = \frac{\sup \{ h^s (b) (0) \mid b \in Q_n \}}{\prod_{m< n} \sigma \bigl( (s \upharpoonright m) * (t\upharpoonright m) \bigr) \bigl( s (m) \bigr)},  \\
\sigma (s * t) (1) & = 1 - \sigma (s * t) (0).  
\end{align*}


We show that $\sigma$ is weakly optimal for player I in the game $\mathcal{G}_{\mathrm{Be}} (A)$. Let $(t_n \mid n < \omega)$ be an $\omega$-sequence of reals such that the tree $T = \{ \emptyset \} \cup \bigcup_{n< \omega} \bigcap_{m > n} T_{m} {\upharpoonright} (n+1) $ is well-founded. We show that the probability of the payoff set via $\mu_{\sigma,\tau_{t}}$ for $t =(t_n \mid n < \omega)$ is greater than $1/2$. (The case when the tree $T$ is ill-founded is dealt with in the same way.) 

First note that together with $(t_n \mid n < \omega)$, $\sigma$ produces a Borel probability $\mu$ on the reals such that for all $n < \omega$ and $s\in 2^n$, $\mu ([s]) = \prod_{m < n } \sigma \bigl( (s\upharpoonright m) * (t \upharpoonright m) \bigr) \bigl( s (m) \bigr)$. Since the tree from $(t_n \mid n< \omega)$ is well-founded, it suffices to show that $\mu(A) > 1/2$. 

Let $(a_n \mid n < \omega )$ be such that for all $n < \omega$, $(a_0 , \ldots a_n) = E(t_0, \ldots , t_{n-1} )$ as in Claim~\ref{claim:BeckerA}.  
Fix an $n < \omega$. Let $H_n = U^{\mathrm{II}, n}_{a_n}$ and $Q_{n} = H_n (a_n)$ as in Claim~\ref{claim:BeckerA}. 
Now set $C^n = \bigl\{ b'_0 \in Q_{0} \mid \text{for some}$ $b'_{n} \in Q_{n}, (b'_m \mid m < n) \text{ is the induced sequence from $b'_n$}  \text{ via } \\ (Q_m \mid m < n)\bigr\}$. 
Then $C^n$ is nonempty and each $b'_0 \in C^n$ uniquelly determines a sequence $(b'_m \mid m \le n)$ such that for all $m < n$, $b'_m$ is identified with $(t_m, a_{m+1}, b'_{m+1})$ via $Q_m$ and $b'_n \in Q_n$. 
We also note that for all $m , n \in \omega$ with $m < n$, $C^m \supseteq C^n$. Since each $C^n$ is perfect in the Cantor space $2^{\omega}$, each $C^n$ is compact and nonempty. Hence the intersection $\bigcap_{n < \omega} C^n$ is nonempty.

We fix an element $b_0$ of $\bigcap_{n < \omega} C^n$ from now on. Then $b_0$ uniquely determines a sequence $(b_n \mid n < \omega)$  such that for all $n< \omega$, $b_n$ is in $Q_n$ and $b_n$ is identified with $(t_n , a_{n+1} , b_{n+1})$ via $Q_n$. 
Since $\tilde{\sigma}$ is optimal for player I in the game $\tilde{\mathcal{G}}_{\mathrm{Be}} (A)$, we have $\mu_{\tilde{\sigma} , \tau_{(a_0, b_0)}} \bigl( \{ \bigl( (x_0, y_0) , (a'_0, b'_0) \bigr) \mid (a'_0, b'_0) = (a_0, b_0) \text{ and } z \in A \iff \text{$T$ is well-founded}\} \bigr) =1$, where $T$ is the tree from $t = (t_n \mid n < \omega)$ as in the description of the game $\tilde{\mathcal{G}}_{\mathrm{Be}} (A)$.  

Also, the measure $\mu_{\tilde{\sigma}, \tau_{(a_0, b_0)}}$ induces a Borel probability measure $\nu$ on the reals as follows: For all $n < \omega$ and $s \in 2^n$, $\nu ([s]) = \mu_{\tilde{\sigma}, \tau_{(a_0,b_0)}} \bigl( \{ \bigl( (x_0, y_0) , (a'_0, b'_0\bigr) \mid (a'_0 , b'_0) = (a_0, b_0) \text{ and } (\forall m < n) \ \rho (u_m) = s(m)\}$, where for each $m < n$, setting $F_{m} = U^{\mathrm{I},m}_{x_{m}}$ and $P_{m} = F_{m} (x_{m})$, the sequence $(y_0, \ldots , y_{n-1})$ is the induced sequence from $y_{n}$ via $(P_{m} \mid m < n)$, and each $y_{m}$ is identified with $(u_{m}, x_{m +1}, y_{m +1})$ via $P_{m}$. 
Since $T$ is well-founded and $\tilde{\sigma}$ is optimal for player I in the game $\tilde{\mathcal{G}}_{\mathrm{Be}} (A)$, $\nu (A) = 1$. 

We now verify that $\mu (A) > 1/2$. By $\nu (A) = 1$ and the choice of $\varepsilon$, it is enough to argue that $d(\mu , \nu) < \epsilon$. 
By the definition of $d$ and the choice of $(\varepsilon_n \mid n < \omega)$, it is enough to see that for all $n < \omega$, $|\mu ([s_n]) - \nu ([s_n])| < \varepsilon_k$ for $k = k_{s_n}$. 

Take any $n < \omega$. We show that $|\mu ([s_n]) - \nu ([s_n])| \le \varepsilon_k$ for $k = k_{s_n}$. If the length of $s_n$ is $0$, then $s_n = \emptyset$, and $\mu ([s_n]) = \mu ([\emptyset]) = 1 = \nu ([\emptyset]) = \nu ([s_n])$. So $| \mu ([s_n]) - \nu ([s_n])| < \varepsilon_k$ for $k = k_{s_n}$, as desired. So assume that the length of $s_n$ is positive,  and let $m +1$ be the length of $s_n$ for some $m \in \omega$. 
Set $s = s_n \upharpoonright m$. 
Then $\mu ([s_n]) =  \prod_{\ell < m+1 } \sigma \bigl( (s_n \upharpoonright \ell) * (t \upharpoonright \ell) \bigr) \bigl( s_n (\ell) \bigr) = \left( \prod_{\ell < m } \sigma \bigl( (s \upharpoonright \ell) * (t \upharpoonright \ell) \bigr) \bigl( s (\ell) \bigr) \right) \cdot \sigma \bigl(s * (t \upharpoonright m) \bigr) \bigl(s_n (m) \bigr)$. We may assume $s_n (m) =0$ (the case $s_n (m) = 1$ is similar). By the definition of $\sigma \bigl( s * (t \upharpoonright m) \bigr)$, $\mu ([s_n]) = \sup \{ h^{s} (b) (0) \mid b \in Q_{m} \}$. But since $b_{m} \in Q_{m}$, by Claim~\ref{claim:BeckerA}~\ref{BeckerA2}, for all $b \in Q_{m}$, $|h^{s} (b) (0) - h^{s} (b_{m}) (0)| < \min \{ \varepsilon_{k_0} , \varepsilon_{k_1}\}$ for $k_i = k_{s^\frown \langle i \rangle}$. Since $\nu ( [s_n] )  =   \mu_{\tilde{\sigma}, \tau_{(a_0,b_0)}} \bigl( \{ \bigl( (x_0, y_0) , (a'_0, b'_0\bigr) \mid (a'_0 , b'_0) = (a_0, b_0) \text{ and } (\forall \ell < m+1) \ \rho (u_{\ell}) = s_n(\ell)\} = h^{s} (b_{m}) (0)$, we have $|\mu ([s_n]) - \nu ([s_n]) | < \varepsilon_k$ for $k = k_{s_n}$, as desired. 
   

It is easy to see that $\sigma$ is in a countable union of sets in $\bigcup_{n < \omega}\undertilde{\mathbf{\Gamma}}^{\mathrm{II}}_n$. 

This completes the proof of Lemma~\ref{lem A}. 
\end{proof}

\begin{proof}[Proof of Lemma~\ref{small conj}]
Let $\tilde{\tau}$ be an optimal strategy for player II in the game $\tilde{\mathcal{G}}_{\mathrm{Be}} (A)$. 
We will construct a weakly optimal strategy $\tau$ for player II in the game $\mathcal{G}_{\mathrm{Be}} (A)$. 


We first prove the following claim which is corresponding to \cite[Lemma~12.14]{ADplus}: 
\begin{Claim}\label{claim:BeckerB}
There is a function $E$ mapping each $(z_0, \ldots , z_n) \in 2^{n+1}$ for each $n \in \omega$ to a sequence of reals $E(z_0, \ldots z_n) = (x_0, u_0, x_1, u_1, \ldots , x_n , u_n)$ such that for all $m \le n$, $\rho (u_m) = z_m$ and setting $F_m = U^{\mathrm{I}, m+1}_{x_m}$ and $P_{m} = F_m (x_m)$, the following hold: 
\begin{enumerate}[(1)]
\item\label{BeckerB1} For all $m \le n$, $P_{m}$ is a perfect set of reals, and for all $y \in P_m$, there is a unique sequence $(y_0, \ldots , y_{m-1}) \in \R^m$ such that for all $\ell < m-1$, $y_{\ell}$ is identified with $(u_{\ell}, x_{\ell +1} , y_{\ell +1})$ via $P_{\ell}$, and $y_{m-1}$ is identified with $(u_{m-1}, x_m, y)$ via $P_{m-1}$. We call the unique sequence $(y_0, \ldots , y_{m-1})$ for $y$ the {\it induced sequence from $y$ via $(P_{\ell} \mid \ell < m)$}. 

\item\label{BeckerB2} For all $m \le n$ and $\vec{\eta} \in \prod_{\ell \le m} \delta^{\mathrm{II}}_{\ell}$, let $f^{\vec{\eta}} \colon P_m \to [0,1]$ be the following function: For each $y  \in P_m$, $f^{\vec{\eta}} (y) = \mu_{\sigma_{(x_0,y_0)}, \tilde{\tau}} \bigl( \{ \bigl( (x'_0,y'_0) , (a_0, b_0) \bigr) \mid (x'_0, y'_0) = (x_0,y_0) \text{ and } \vec{\eta} \in T_m \}\bigr)$, where $T_m = \pi^{\mathrm{II}}_m (t_m)$, $(y_0, \ldots , y_{m-1})$ is the induced sequence from $y$ via $(P_{\ell} \mid \ell <m)$, and for all $\ell \le m$, setting $H_{\ell} = U^{\mathrm{II}, \ell}_{a_{\ell}}$ and $Q_{\ell} = H_{\ell} (a_{\ell})$, the sequence $(a_0, \ldots , a_{m})$ is the induced sequence from $a_{m+1}$ via $(Q_{\ell} \mid \ell \le m)$, and each $b_{\ell}$ is identified with $(t_{\ell}, a_{\ell+1}, b_{\ell +1})$ via $Q_{\ell}$. 

Then for all $m \le n$ and $\vec{\eta} \in \prod_{\ell \le m} \delta^{\mathrm{II}}_{\ell}$, $f^{\vec{\eta}}$ is continuous on $P_m$ and for all $y , y' \in P_m$, if $(y_0, \ldots , y_{m-1})$ is the induced sequence from $y$ via $(P_{\ell} \mid \ell < m)$ and $(y'_0, \ldots , y'_{m-1})$ is the induced sequence from $y'$ via $(P_{\ell} \mid \ell < m)$, then for all $\ell < m$, $y_{\ell} \upharpoonright m = y'_{\ell} \upharpoonright m$, and $y \upharpoonright m = y' \upharpoonright m$. 
\end{enumerate}
\end{Claim}




\begin{proof}[Proof of Claim~\ref{claim:BeckerB}]
We construct $E(z_0, \ldots z_n) = (x_0, u_0, x_1, u_1, \ldots , x_n , u_n)$ by induction on $n$. 

Let $n =0$ and $z_0 \in 2$. We will construct $E(z_0) = (x_0, u_0)$. 
Fix a real $x$ and an ordinal $\eta \in \delta^{\mathrm{II}}_0$. Let $f^{\eta}_x \colon 2^{\omega} \to [0,1]$ be as follows: Given a real $y$, set $f^{\eta}_x (y) = \mu_{\sigma_{(x,y)}, \tilde{\tau}} \bigl( \{ \bigl( (x'_0,y'_0) , (a_0, b_0) \bigr) \mid (x'_0, y'_0) = (x,y) \text{ and } \eta \in T_0 \}\bigr)$, where $T_0 = \pi^{\mathrm{II}}_0 ( t_0 )$, and setting $H_0 = U^{\mathrm{II}, 0}_{a_0}$ and $Q_0 = H_0 (a_0)$, $b_0$ is identified with $(t_0, a_1, b_1)$ via $Q_0$. 
By $\BlADR$ and Theorem~\ref{Blackwell Baire property}, every set of reals has the Baire property. So by Theorem~\ref{fact:category}, for each $\eta \in \delta^{\mathrm{II}}_0$, the function $f^{\eta}_x$ is continuous on a comeager set. 
By Theorem~\ref{fact:category} again, the meager ideal on the Cantor space is closed under well-ordered unions. Hence for some comeager set of reals $D^0$, for all $\eta \in \delta^{\mathrm{II}}_0$, the function $f^{\eta}_x$ is continuous on $D^0$. 
So for some perfect set $P$, the following holds:
\begin{description}
\item[$(\ast)^0_{x,P}\, $]\label{BeckerBn=0} For all $\eta \in \delta^{\mathrm{II}}_0$, $f^{\eta}_x$ is continuous on $P$. 
\end{description}
Let $R^0 = \{ (x, P)\mid  (\ast)^0_{x,P} \text{ holds} \}$.  
Then $R^0$ is projective in a set in $\undertilde{\mathbf{\Gamma}}^{\mathrm{II}}_0$. So for some function $F_0$ in $\undertilde{\mathbf{\Gamma}}^{\mathrm{I}}_1$, $F_0$ uniformizes $R^0$. 
Fix a real $x_0$ with $U^{\mathrm{I}, 1}_{x_0} = F_0$. 
Let $u_0$ be a real such that $\rho (u_0) = z_0$. 
Set $E(z_0) = (x_0, u_0)$. Then $E(z_0)$ satisfies \ref{BeckerB1} and \ref{BeckerB2} in the claim for $m=n=0$.  

Now let $n > 0$. Take any $(z_0, \ldots , z_{n}) \in 2^{n+1}$. We will construct $E(z_0, \ldots z_n)$. By induction hypothesis, we already obtained $E(z_0, \ldots z_{n-1}) =  (x_0, u_0, x_1, u_1, \ldots , x_{n-1} ,  \\ u_{n-1})$ and the sequence $(P_m \mid m \le n-1)$ that satisfy \ref{BeckerB1} and \ref{BeckerB2} in the claim for all $m \le n-1$. 
We will find a pair of reals $(x_n, u_n)$ such that $\rho (u_n) = z_n$ and that the sequences $(x_0, u_0, x_1, u_1, \ldots , x_{n} , u_{n})$ and $(P_m \mid m \le n)$ satisfy  \ref{BeckerB1} and \ref{BeckerB2} in the claim for $m =n$.

Fix a real $x$ and $\vec{\eta} \in \prod_{m \le n} \delta^{\mathrm{II}}_m$. Define the function $f^{\vec{\eta}}_x \colon 2^{\omega} \to [0,1]$ as follows: For a real $y$, let $(y_m \mid m < n)$ be the induced sequence from $y$ via $(P_m \mid m < n )$.  
Set $f^{\vec{\eta}}_x (y) = \mu_{\sigma_{(x_0, y_0)}, \tilde{\tau}} \bigl( \{ \bigl( (x'_0, y'_0), (a_0, b_0)\bigr) \mid (x'_0, y'_0) = (x_0,y_0) \text{ and } \vec{\eta} \in T_{n}\} \bigr)$, where $T_{n} =   \pi^{\mathrm{II}}_n ( t_{n} )$, and for all $m \le n$, setting $H_m = U^{\mathrm{II}, m}_{a_m}$ and $Q_m = H_m (a_m)$, $b_m$ is identified with $(t_m, a_{m+1}, b_{m+1})$ via $Q_m$. 
By the same argument as in the case $n=0$, for some perfect set $P$, the following holds:
\begin{description}
\item[$(\ast)^{n}_{x,P}\, $] For all $\vec{\eta} \in \prod_{m \le n} \delta^{\mathrm{II}}_m$, $f^{\vec{\eta}}_x$ is continuous on $P$ and for all $y, y' \in P$, if $(y_m \mid m < n)$ is the induced sequence from $y$ via $(P_m \mid m < n)$ and $(y'_m \mid m < n)$ is the induced sequence from $y'$ via $(P_m \mid m < n)$, then for all $m < n$, $y_m \upharpoonright n = y'_m \upharpoonright n$, and $y\upharpoonright n = y' \upharpoonright n$. 
\end{description}
Let $R^{n} = \{ (x, P) \mid (\ast)^{n}_{x,P} \text{ holds}\}$.  
Then $R^{n}$ is projective in a set in $\undertilde{\mathbf{\Gamma}}^{\mathrm{II}}_{n}$. So there is a function $F_{n}$ in $\undertilde{\mathbf{\Gamma}}^{\mathrm{I}}_{n+1}$ such that $F_{n}$ uniformizes $R^{n}$. 
Fix a real $x_{n}$ with $U^{\mathrm{I}, n+1}_{x_{n}} = F_{n}$. 
Let $u_{n}$ be a real such that $\rho (u_{n}) = z_{n}$. Then setting $P_n = F_n (x_n)$, the sequences $(x_0, u_0, x_1, u_1, \ldots , x_{n} , u_{n})$ and $(P_m \mid m \le n)$ satisfy \ref{BeckerB1} and \ref{BeckerB2} in the claim for $m =n$, as desired. 

This finishes the proof of Claim~\ref{claim:BeckerB}. 
\end{proof}

We can now define a mixed strategy $\tau$ for player II in the game $\mathcal{G}_{\mathrm{Be}} (A)$ as follows: Let $n < \omega$, $s = (z_0, \ldots z_n) \in 2^{n+1}$, and $w = (t_0, \ldots , t_{n-1}) \in \R^n$. Let $s * w$ be the concatenation of $s$ and $w$ bit by bit, i.e., $s*w = (z_0, t_0, z_1, t_1, \ldots , z_{n-1}, t_{n-1}, z_n)$. 
Let $E(z_0, \ldots , z_{n}) = (x_0, u_0,  \ldots , x_n, u_n)$, $(P_m \mid m \le n)$, and $f^{\vec{\eta}} \colon P_n \to [0,1]$ be as in Claim~\ref{claim:BeckerB}. 
For each $\vec{\eta} \in \prod_{m \le n} \delta^{\mathrm{II}}_m$, set $p_n (\vec{\eta}) = \sup \{ f^{\vec{\eta}} (y) \mid y\in P_n \}$. 
For each $k < 2^{n+1}$, set 
\[T_n^k = \left\{ \vec{\eta} \in \prod_{m\le n} \delta^{\mathrm{II}}_m\  \middle| \ \frac{k}{2^{n+1}} < p_n (\vec{\eta}) \le 1\right\}.  
\] 
Let $t_n^k$ be a real such that $\pi^{\mathrm{II}}_n (t_n^k) = T_n^k$. 

We now set $\tau (s * w) \in \text{Prob}_{\omega}(\R)$ as follows: when $n=0$, for each real $t$, 
\[\tau (s * w) (t) = \begin{cases}
\ \ \frac{1}{2} & \text{ if $t = t_0^k$ for some $k < 2$}, \\
\ \ 0 & \text{ otherwise}. 
\end{cases}
\]
When $n >0$, for each real $t$,  
\[\tau (s * w) (t) = \begin{cases}
\ \ \frac{1}{2} & \text{ if for some $j < 2^{n}$, $t_{n-1} = t_{n-1}^j$ and (either $t = t_n^{2j}$ or $t= t_n^{2j+1}$)}, \\
\ \ 0 & \text{ otherwise}. 
\end{cases}
\]
Then $\tau (s * w)$ is of finite support, i.e., the set $\{ t \mid \tau (s * w) (t) \neq 0\}$ is finite. 
Moreover, $\tau$ enjoys the following properties: 
\begin{Claim}\label{tau}

${}$

\begin{enumerate}[(1)]
\item\label{plus0} Let $n < \omega$ and $s \in 2^{n+1}$. Then for each $k < 2^{n+1}$, there is a unique $w \in \R^n$ such that $\tau (s * w) (t_n^k) \neq 0$ and for all $m < n $, $\tau \bigl( (s \upharpoonright m+1) * (w \upharpoonright m)\bigr) \bigl( w (m) \bigr) \neq 0$. 

\item\label{plus0.5} For all $n < \omega$, $k < 2^{n+1}$, and $z \in \R$, $\mu_{\sigma_z , \tau} \bigl( \bigl\{ (z' ,\vec{t}) \mid z' = z \text{ and } t_n = t_n^k \bigr\} \bigr) = \frac{1}{2^{n+1}}$.  

\item\label{plus} For all $n < \omega$, $\vec{\eta} \in \prod_{m \le n } \delta^{\mathrm{II}}_m$, and $z \in \R$, $p_n (\vec{\eta}) \le \mu_{\sigma_z, \tau} \bigl( \{ (z',\vec{t}) \mid z' = z \text{ and } \vec{\eta} \in T_n \} \bigr) < p_n(\vec{\eta}) + \frac{1}{2^{n+1}}$. 

\item\label{plus2} For all $n < \omega$, $X \subseteq \prod_{m \le n} \delta^{\mathrm{II}}_m$, and $z \in 2^{\omega}$, 
$\mu_{\sigma_z, \tau} \bigl( \{ (z',\vec{t}) \mid z' = z \text{ and } (\exists \vec{\eta} \in X) \, \vec{\eta} \in T_n \} \bigr) \le \sup \{ p_n (\vec{\eta}) \mid \vec{\eta} \in X \} + \frac{1}{2^{n+1}}$. 

\item\label{plus3} For all reals $z$ and $n < \omega$, $\mu_{\sigma_z , \tau} \bigl( \{ (z', \vec{t}) \mid z'=z \text{ and }T_{n+1} \upharpoonright (n+1) \subseteq T_n \} \bigr) =1$. 
\end{enumerate}
\end{Claim}

\begin{proof}[Proof of Claim~\ref{tau}]
For \ref{plus0}, we argue by induction on $n$. For $n=0$, it is obvious. Assume $n > 0$. Let $w = (t_0, \ldots t_{n-1}) \in \R^n$ be such that $\tau (s * w) (t_n^k) \neq 0$. Then by the definition of $\tau (s * w)$, letting $j$ be a unique natural number with $j < 2^n$ and either $k = 2j$ or $k=2j+1$, we have $t_{n-1} = t_{n-1}^j$. By induction hypothesis, there is a unique $w'\in \R^{n-1}$ such that $\tau\bigl( (s \upharpoonright n * w'\bigr) (t_{n-1}^j) \neq 0$ and for all $m < n-1 $, $\tau \bigl( (s \upharpoonright m+1) * (w' \upharpoonright m)\bigr) \bigl( w' (m) \bigr) \neq 0$. 
Hence a desired  $w$ must exist and be of the form $w'^{\frown} \langle t_{n-1}^j\rangle$, which is unique for $s$ and $k$. 

For \ref{plus0.5}, let $n < \omega$, $k < 2^{n+1}$, and $z \in \R$. By \ref{plus0} for $s = z \upharpoonright (n+1)$,  there is a unique $w_k \in \R^n$ such that $\tau (s * w_k) (t_n^k) \neq 0$ and for all $m < n $, $\tau \bigl( (s \upharpoonright m+1) * (w_k \upharpoonright m)\bigr) \bigl( w_k (m) \bigr) \neq 0$. 
Then 
\begin{align*}
& \ \mu_{ \sigma_z , \tau} \bigl( \bigl\{ (z' ,\vec{t}) \mid z'  = z \text{ and } t_n = t_n^k \bigr\} \bigr) \\
 = &  \sum_{w \in \R^n} \prod_{m < n} \tau \bigl( (z \upharpoonright m+1) * (w \upharpoonright m) \bigr) \bigl( w(m) \bigr) \cdot \tau(z *w) (t_n^k) \\
 = &  \prod_{m < n} \tau \bigl( (z \upharpoonright m+1) * (w_k \upharpoonright m) \bigr) \bigl( w_k(m) \bigr) \cdot \tau(z *w_k) (t_n^k) \\
 = & \ \frac{1}{2^{n+1}}. 
\end{align*}

For \ref{plus}, let $n < \omega$, $\vec{\eta} \in \prod_{m \le n } \delta^{\mathrm{II}}_m$, and $z \in \R$. Let $k_{\vec{\eta}} = \min \{ k \le 2^{n+1} \mid \vec{\eta} \notin T_n^k \text{ or } k = 2^{n+1} \}$. By the definition of $T_n^k$, it is enough to see that $\mu_{\sigma_z, \tau} \bigl( \{ (z',\vec{t}) \mid z' = z \text{ and } \vec{\eta} \in T_n \} \bigr) = \frac{k_{\vec{\eta}}}{2^{n+1}}$. 
Now 
\begin{align*}
\mu_{\sigma_z, \tau} \bigl( \{ (z',\vec{t}) \mid z' = z \text{ and } \vec{\eta} \in T_n \} \bigr)  
& = \mu_{\sigma_z, \tau} \bigl( \{ (z',\vec{t}) \mid z' = z \text{ and } \vec{\eta} \in T_n^k \text{ for some $k$} \} \bigr) \\
& = \sum_{\vec{\eta} \in T_n^k} \mu_{\sigma_z, \tau} \bigl( \{ (z',\vec{t}) \mid z' = z \text{ and } t_n = t_n^k\} \bigr) \\
& = \sum_{\vec{\eta} \in T_n^k} \frac{1}{2^{n+1}} = \frac{k_{\vec{\eta}}}{2^{n+1}},  
\end{align*}
where the third equality above follows from \ref{plus0.5}. 

For \ref{plus2}, let $n < \omega$ and $X \subseteq \prod_{m \le n } \delta^{\mathrm{II}}_m$. 
Set $p_X = \sup \{ p_n (\vec{\eta}) \mid \vec{\eta} \in X\}$ and let $k_X = \min \{ k \le 2^{n+1} \mid p_X <  \frac{k}{2^{n+1}} \text{ or } k = 2^{n+1} \}$. Then it is enough to see that $\mu_{\sigma_z, \tau} \bigl( \{ (z',\vec{t}) \mid z' = z \text{ and } (\exists \vec{\eta} \in X) \, \vec{\eta} \in T_n \} \bigr) \le \frac{k_X}{2^{n+1}}$. 
Now
\begin{align*}
  & \ \mu_{\sigma_z, \tau} \bigl( \{ (z',\vec{t}) \mid z' = z \text{ and } (\exists \vec{\eta} \in X) \, \vec{\eta} \in T_n \} \bigr)  \\
\le & \ \mu_{\sigma_z, \tau} \bigl( \{ (z',\vec{t}) \mid z' = z \text{ and } t_n \neq t_n^{\ell} \text{ for any $\ell \ge k_X$} \} \bigr) \\
=  & \ \mu_{\sigma_z, \tau} \bigl( \{ (z',\vec{t}) \mid z' = z \text{ and } t_n = t_n^{j} \text{ for some $j < k_X$} \} \bigr)\\
= & \sum_{j < k_X} \mu_{\sigma_z, \tau} \bigl( \{ (z',\vec{t}) \mid z' = z \text{ and } t_n = t_n^j\} \bigr) \\
= &  \sum_{j < k_X} \frac{1}{2^{n+1}} = \ \frac{k_X}{2^{n+1}},  
\end{align*}
where the third equality above follows from \ref{plus0.5}. 

For \ref{plus3}, let $z \in \R$ and $n < \omega$. 
We first claim that for all $\vec{\eta} \in \prod_{m \le n+1} \delta^{\mathrm{II}}_m$, $p_{n+1} (\vec{\eta}) \le p_n \bigl(\vec{\eta} \upharpoonright (n+1) \bigr)$. By the definitions of $p_{n+1} (\vec{\eta})$ and $p_n \bigl(\vec{\eta} \upharpoonright (n+1) \bigr)$, the desired inequality follows from that for all $y \in P_{n+1}$, if $(y'_m \mid m \le n)$ is the induced sequence from $y$ via $(P_m \mid m \le n)$, then $y'_n \in P_n$ and $f^{\vec{\eta}} (y) \le f^{\vec{\eta} \upharpoonright (n+1)} (y'_n)$, which in turn follows from that $\tilde{\tau}$ is optimal for player II in the game $\tilde{\mathcal{G}}_{\mathrm{Be}} (A)$. 

To see $\mu_{\sigma_z , \tau} \bigl( \{ (z', \vec{t}) \mid z'=z \text{ and }T_{n+1} \upharpoonright (n+1) \subseteq T_n \} \bigr) =1$, it is enough to verify that for all $k < 2^{n+2}$, $\mu_{\sigma_z , \tau } \bigl( \{ (z', \vec{t}) \mid z' = z, t_{n+1} = t_{n+1}^k, \text{ and } T_{n+1}^k \upharpoonright (n+1) \subseteq T_n \} = \frac{1}{2^{n+2}}$. 
Fix any $k< 2^{n+2}$. Let $w = (t_0, \ldots , t_n) \in \R^{n+1}$ be such that $\tau \bigl( (z\upharpoonright n+2) * w\bigr) (t_{n+1}^k) \neq 0$. Then by the definition of $\tau\bigl( (z\upharpoonright n+2) * w\bigr)$, letting $j$ be a unique natural number with $j < 2^{n+1}$ and either $k=2j$ or $k=2j+1$, we have $t_n = t_n^j$. 
Hence it is enough to show that $\mu_{\sigma_z , \tau } \bigl( \{ (z', \vec{t}) \mid z' = z, t_{n+1} = t_{n+1}^k, \text{ and } T_{n+1}^k \upharpoonright (n+1) \subseteq T_n^j \} \bigr) = \frac{1}{2^{n+2}}$. 
Therefore, it suffices to verify that for all $\vec{\eta} \in T_{n+1}^k$, $\vec{\eta} \upharpoonright (n+1) \in T_n^j$. 
Take any $\vec{\eta} \in T_{n+1}^k$. Then $\frac{k}{2^{n+2}} < p_{n+1} (\vec{\eta}) \le 1$. 
Since  $p_{n+1} (\vec{\eta}) \le p_n \bigl(\vec{\eta} \upharpoonright (n+1) \bigr)$ and either $k =2j$ or $k=2j+1$, we have $\frac{j}{2^{n+1}} < p_n  \bigl(\vec{\eta} \upharpoonright (n+1) \bigr) \le 1$. Hence $\vec{\eta} \upharpoonright (n+1) \in T_n^j$, as desired. 

This completes the proof of Claim~\ref{tau}. 
\end{proof}

We show that $\tau$ is weakly optimal for player II in the game $\mathcal{G}_{\mathrm{Be}} (A)$. Let $z = (z_n \mid n < \omega)$ be a real. 
We will show that $\mu_{\sigma_z , \tau} \bigl( \bigl\{ (z' , \vec{t}) \mid z'=z, (\forall n < \omega) \ T_{n+1} \upharpoonright (n+1) \subseteq T_n, \text{ and }z \in A \iff T\text{ is ill-founded} \bigr\}\bigr) > \frac{1}{2}$, where $T_m = \pi^{\mathrm{II}}_m (t_m)$ and $T =  \{ \emptyset \} \cup \bigcup_{n < \omega} \bigcap_{m > n} T_m \upharpoonright (n+1)$.  
By Claim~\ref{tau}~\ref{plus3} above, for all $n< \omega$, $\mu_{\sigma_z , \tau} \bigl( \{ (z', \vec{t}) \mid z'=z \text{ and }T_{n+1} \upharpoonright (n+1) \subseteq T_n \} \bigr) =1$. 
So it suffices to show that $\mu_{\sigma_z , \tau} \bigl( \bigl\{ (z' , \vec{t}) \mid z'=z \text{ and }z \in A \iff T\text{ is ill-founded} \bigr\}\bigr) > \frac{1}{2}$. 

We assume that $z \in A$ (the case $z \notin A$ is similarly dealt with). Then it is enough to verify that $\mu_{\sigma_z , \tau} \bigl( \bigl\{ (z' , \vec{t}) \mid z'=z \text{ and } T\text{ is ill-founded} \bigr\}\bigr) > \frac{1}{2}$. 

First note that together with $z = (z_n \mid n < \omega)$, $\tau$ produces a Borel probability measure $\mu$ on $\R^{\omega}$ such that for all $n < \omega$ and $w\in \R^n$, $\mu ([w]) = \prod_{m < n } \tau \bigl( (z\upharpoonright (m+1)) * (w \upharpoonright m) \bigr) \bigl( w (m) \bigr)$. Since for all $B \subseteq \R^{\omega}$, $\mu (B) = \mu_{\sigma_z , \tau} \bigl( \{ (z' ,\vec{t}) \mid z' =z \text{ and } \vec{t} \in B \}\bigr)$, it suffices to show that $\mu\bigl(\{ \vec{t} \mid T \text{ is ill-founded}\}\bigr) > 1/2$. 

Let $\bigl( (x_n,u_n) \mid n < \omega \bigr)$ be such that for all $n < \omega$, $(x_0 , u_0, \ldots x_n, u_n) = E(z_0, \ldots , z_{n} )$ as in Claim~\ref{claim:BeckerB}.  
Fix an $n < \omega$. Let $F_n = U^{\mathrm{I}, {n+1}}_{x_n}$ and $P_{n} = F_n (x_n)$. 
Now set $C^n = \bigl\{ y'_0 \in P_{0} \mid \text{for some}$ $y'_{n} \in P_{n}, (y'_m \mid m < n) \text{ is the induced sequence from $y'_n$}  \text{ via } (P_m \mid m < n)\bigr\}$. 
Then $C^n$ is nonempty and each $y'_0 \in C^n$ uniquelly determines a sequence $(y'_m \mid m \le n)$ such that for all $m < n$, $y'_m$ is identified with $(u_m, x_{m+1}, y'_{m+1})$ via $P_m$ and that $y'_n \in P_n$. 
We also note that for all $m , n \in \omega$ with $m < n$, $C^m \supseteq C^n$. Since each $C^n$ is perfect in the Cantor space $2^{\omega}$, each $C^n$ is compact and nonempty. Hence the intersection $\bigcap_{n < \omega} C^n$ is nonempty.

We fix an element $y_0$ of $\bigcap_{n < \omega} C^n$ from now on. Then $y_0$ uniquely determines a sequence $(y_n \mid n < \omega)$  such that for all $n< \omega$, $y_n$ is in $P_n$ and $y_n$ is identified with $(u_n , x_{n+1} , y_{n+1})$ via $P_n$. 
Since $\tilde{\tau}$ is optimal for player II in the game $\tilde{\mathcal{G}}_{\mathrm{Be}} (A)$, we have $\mu_{\sigma_{(x_0, y_0)} , \tilde{\tau}} \bigl( \left\{ \bigl( (x'_0, y'_0) , (a_0, b_0) \bigr) \mid (x'_0, y'_0) = (x_0, y_0) \text{ and } z \in A \iff \text{$T$ is ill-founded}\right\} \bigr) =1$. 
Since $z \in A$, $\mu_{\sigma_{(x_0, y_0)} , \tilde{\tau}} \bigl( \left\{ \bigl( (x'_0, y'_0) , (a_0, b_0) \bigr) \mid (x'_0, y'_0) = (x_0, y_0) \text{ and } \text{$T$ is ill-founded}\right\} \bigr) =1$. 

Also, the measure $\mu_{\sigma_{(x_0, y_0)}, \tilde{\tau}}$ induces a Borel probability measure $\nu$ on $\R^{\omega}$ as follows: For all $n < \omega$ and $w \in \R^n$, set $\nu ([w]) = \mu_{\sigma_{(x_0, y_0)}, \tilde{\tau}} \bigl( \bigl\{ \bigl( (x'_0, y'_0) , (a_0, b_0)\bigr) \mid (x'_0 , y'_0) = (x_0, y_0) \text{ and } (\forall m < n) \ t_m = w(m)\bigr\}$. 
Since $\mu_{\sigma_{(x_0, y_0)} , \tilde{\tau}} \bigl( \bigl\{ \bigl( (x'_0, y'_0) , (a_0, b_0) \bigr) \mid (x'_0, y'_0) = (x_0, y_0) \text{ and } \text{$T$ is ill-founded}\bigr\} \bigr) =1$,  $\nu \bigl(  \{ \vec{t}  \mid \text{$T$ is ill-founded}\}\bigr) = 1$. 

Since $\nu \bigl(  \{ \vec{t}  \mid \text{$T$ is ill-founded}\}\bigr) = 1$, to verify $\mu\bigl(\{ \vec{t} \mid T \text{ is ill-founded}\}\bigr) > 1/2$, it is enough to see that $ \mu\bigl(\{ \vec{t} \mid T \text{ is ill-founded}\}\bigr) = \nu\bigl(\{ \vec{t} \mid T \text{ is ill-founded}\}\bigr)$, which would actually imply $\mu\bigl(\{ \vec{t} \mid T \text{ is ill-founded}\}\bigr) = 1$.

To see $ \mu\bigl(\{ \vec{t} \mid T \text{ is ill-founded}\}\bigr) = \nu\bigl(\{ \vec{t} \mid T \text{ is ill-founded}\}\bigr)$, 
we first show the following: 
\begin{Claim}\label{measure calculations}
${}$

\begin{enumerate}[(1)]
\item\label{ichi} For all $n < \omega$ and $\vec{\eta} \in \prod_{m \le n} \delta^{\mathrm{II}}_m$, $\nu \bigl( \{ \vec{t} \mid \vec{\eta} \in T_n \} \bigr) \le  \mu \bigl( \{ \vec{t} \mid \vec{\eta} \in T_n \} \bigr)$. 

\item\label{ni} For all $n < \omega$ and $\vec{\eta} \in \prod_{m \le n} \delta^{\mathrm{II}}_m$, 

$\mu \left( \left\{ \vec{t} \ \middle| \ \vec{\eta} \in \bigcap_{m > n} T_m \upharpoonright (n+1) \right\}\right) \le \nu \left( \left\{ \vec{t} \mid \vec{\eta} \in T_n \right\}\right).
$
\end{enumerate}
\end{Claim}

\begin{proof}[Proof of Claim~\ref{measure calculations}]

We first verify \ref{ichi}. Since $p_n (\vec{\eta}) = \sup \{ f^{\vec{\eta}} (y) \mid y \in P_{x_n}\}$ and $y_n$ is in $P_n$, $p(\vec{\eta}) \ge f^{ \vec{\eta}} (y_n) = \nu \bigl( \{ \vec{t} \mid \vec{\eta} \in T_n \}\bigr)$. Also by Claim~\ref{tau}~\ref{plus}, $p_n (\vec{\eta}) \le \mu_{\sigma_z, \tau} \bigl( \{ (z',\vec{t}) \mid z' = z \text{ and } \vec{\eta} \in T_n \} \bigr) =  \mu \bigl( \{ \vec{t} \mid \vec{\eta} \in T_n \} \bigr)$. 
Hence $\nu \bigl( \{ \vec{t} \mid \vec{\eta} \in T_n \} \bigr) \le  \mu \bigl( \{ \vec{t} \mid \vec{\eta} \in T_n \} \bigr)$. 

We next verify \ref{ni}. Set $r = \mu \bigl( \bigl\{ \vec{t} \mid \vec{\eta} \in \bigcap_{m > n} T_m \upharpoonright (n+1) \bigr\}\bigr)$ and $s = \nu \bigl( \bigl\{ \vec{t} \mid \vec{\eta} \in T_n \bigr\}\bigr)$. We will argue that $r \le s $. 

Suppose not. Then $r > s$. 
Since $s = f^{\vec{\eta}} (y_n)$ and $f^{\vec{\eta}}$ is continuous on $P_{n}$ by Claim~\ref{claim:BeckerB}, there exists some $\ell \in \omega$ with $\ell > n$ such that for all $y'_n \in P_{n}$, if $y'_n \upharpoonright \ell = y_n \upharpoonright \ell$, then $f^{\vec{\eta}}(y'_n) + \frac{1}{2^{\ell + 1}} < \frac{r+s}{2} < r$. 
However, by Claim~\ref{claim:BeckerB}, for all $y'_{\ell} \in P_{\ell}$, $y'_{n} \upharpoonright \ell = y_{n} \upharpoonright \ell$. 
Hence for all $y'_{\ell} \in P_{\ell}$, we have $f^{\vec{\eta}} (y'_n) + \frac{1}{2^{\ell +1}} < \frac{r+s}{2}$. Now
\begin{align*}
r & = \mu \bigl( \bigl\{ \vec{t} \mid \vec{\eta} \in \bigcap_{m > n} T_m \upharpoonright (n +1) \bigr\}\bigr)\\
  & \le  \mu \bigl( \bigl\{ \vec{t} \mid \vec{\eta} \in T_{\ell} \upharpoonright (n+1) \bigr\}\bigr)\\
  & \le \sup \{ p_{\ell} (\vec{\zeta}) \mid \vec{\zeta} \supseteq \vec{\eta} \} + \frac{1}{2^{\ell +1}}\\
  & = \left( \sup_{\vec{\zeta} \supseteq \vec{\eta}} \ \sup_{ y'_{\ell} \in P_{x_{\ell}}} f^{\vec{\zeta}} (y'_{\ell}) \right)+ \frac{1}{2^{\ell +1}}\\
  & = \left( \sup_{ y'_{\ell} \in P_{x_{\ell}}} \ \sup_{\vec{\zeta} \supseteq \vec{\eta}}  f^{\vec{\zeta}} (y'_{\ell}) \right)+ \frac{1}{2^{\ell +1}}\\
  & =  \left( \sup_{y'_{\ell} \in P_{x_{\ell}}} \ \sup_{\vec{\zeta} \supseteq \vec{\eta}} \mu_{\sigma_{(x_0, y'_0)} , \tilde{\tau}} \bigl(\bigl\{\bigl( (x''_0, y''_0), (a_0, b_0) \bigr) \mid (x''_0, y''_0) = (x_0, y'_0) \text{ and }\vec{\zeta} \in T_{\ell} \bigr\} \bigr) \right)+ \frac{1}{2^{\ell +1}}\\
  & \le  \sup_{y'_{\ell} \in P_{x_{\ell}}} \mu_{ \sigma_{(x_0, y'_0)} , \tilde{\tau}} \bigl(\bigl\{\bigl( (x''_0, y''_0), (a_0, b_0) \bigr) \mid  (x''_0, y''_0) = (x_0, y'_0) \text{ and }(\exists \vec{\zeta} \supseteq \vec{\eta})\ \vec{\zeta} \in T_{\ell}  \bigr\} \bigr) + \frac{1}{2^{\ell +1}}\\
  & = \sup_{y'_{\ell} \in P_{x_{\ell}}}   \mu_{ \sigma_{(x_0, y'_0)} , \tilde{\tau}} \bigl(\bigl\{\bigl( (x''_0, y''_0), (a_0, b_0) \bigr) \mid  (x''_0, y''_0) = (x_0, y'_0) \text{ and }\vec{\eta} \in T_{\ell} \upharpoonright (n+1) \bigr\} \bigr) + \frac{1}{2^{\ell +1}}\\
  & \le \sup_{y'_{\ell} \in P_{x_{\ell}}} \mu_{ \sigma_{(x_0, y'_0)} , \tilde{\tau}} \bigl(\bigl\{\bigl( (x''_0, y''_0), (a_0, b_0) \bigr) \mid  (x''_0, y''_0) = (x_0, y'_0) \text{ and }\vec{\eta} \in T_n \bigr\} \bigr)  + \frac{1}{2^{\ell +1}}
\end{align*}

\begin{flalign*}
 & \ \ = \sup_{y'_{\ell} \in P_{x_{\ell}}} f^{\vec{\eta}} (y'_n) + \frac{1}{2^{\ell +1}}&\\
  & \ \ \le \frac{r+s}{2} &\\
  & \ \ < r , &
\end{flalign*}
where the second inequality above follows from Claim~\ref{tau}~\ref{plus2} above. The fourth inequality above follows from the condition that $\tilde{\tau}$ is optimal in the game $\tilde{\mathcal{G}}_{\mathrm{Be}} (A)$, and that if player II wins in the game $\tilde{\mathcal{G}}_{\mathrm{Be}} (A)$, then for all $m< \omega$, $T_{m+1} \upharpoonright (m+1) \subseteq T_m$. 

The above calculation gives us $r < r$. Contracition! Therefore, $r \le s$, as desired.
This completes the proof of Claim~\ref{measure calculations}. 
\end{proof}

We finally verify that $ \mu\bigl(\{ \vec{t} \mid T \text{ is ill-founded}\}\bigr) = \nu\bigl(\{ \vec{t} \mid T \text{ is ill-founded}\}\bigr)$. 
Let $r' =  \mu\bigl(\{ \vec{t} \mid T \text{ is ill-founded}\}\bigr) $ and $s' = \nu\bigl(\{ \vec{t} \mid T \text{ is ill-founded}\}\bigr)$.  
We will argue that $s' \le r'$ and $s' \ge r'$. 
By Claim~\ref{measure calculations}~\ref{ichi}, it follows that $s' \le r'$. 

Now 
\begin{align*}
s' & =  \nu \bigl( \bigl\{ \vec{t} \mid T = \{ \emptyset \} \cup \bigcup_{n < \omega} \bigcap_{m > n} T_m \upharpoonright (n +1) \text{ is ill-founded} \bigr\}\bigr)\\
& = \nu \bigl( \bigl\{ \vec{t} \mid  (\exists g \in \prod_{n < \omega} \delta^{\mathrm{II}}_n ) \, ( \forall n < \omega) \, g\upharpoonright (n+1) \in  \bigcap_{m > n} T_m \upharpoonright (n+1) \bigr\}\bigr)\\
 & = \nu \bigl( \bigl\{ \vec{t} \mid  (\exists g \in \prod_{n < \omega} \delta^{\mathrm{II}}_n ) \, ( \forall n < \omega) \, g\upharpoonright (n+1) \in  T_n \bigr\}\bigr)\\
& \ge \mu \bigl( \bigl\{  \vec{t} \mid (\exists g \in \prod_{n < \omega} \delta^{\mathrm{II}}_n ) \, ( \forall n < \omega) \, g\upharpoonright (n+1) \in \bigcap_{m > n} T_m \upharpoonright (n+1) \bigr\}\bigr)\\
 & = \mu \bigl( \bigl\{ \vec{t} \mid  T = \{ \emptyset \} \cup \bigcup_{n < \omega} \bigcap_{m > n} T_m \upharpoonright n \text{ is ill-founded} \bigr\}\bigr)\\
& = r'.
\end{align*}
The inequality above follows from Claim~\ref{measure calculations}~\ref{ni}. 

The above calculation shows $s' \ge r'$. Hence $r' = s'$, as desired.  

This completes the proof of Lemma~\ref{small conj}. 
\end{proof}

\section{Equiconsistency between $\ADR$ and $\BlADR$}\label{sec:Con}

In the last section, we proved the equivalence between $\ADR$ and $\BlADR$ in $\ZF$+$\DC$. However, assuming $\DC$ to see the equivalence between $\ADR$ and $\BlADR$ is not optimal by the following theorem: 
\begin{Thm}[Solovay]
The theory $\ZF$+$\ADR$+$\DC$ proves the consistency of $\ZF$+$\ADR$. Hence the consistency of $\ZF$+$\ADR$+$\DC$ is strictly stronger than that of $\ZF$+$\ADR$. 
\end{Thm}

\begin{proof}
See \cite[Corollary~5.8]{Solovay_AD_R}.
\end{proof}

One can ask whether $\ADR$ and $\BlADR$ are equivalent in $\ZF$+$\AC_{\omega}(\R)$ without assuming $\DC$. 
While we do not know the answer to this question, we prove the following:
\begin{Thm}[Theorem~\ref{equiconsistency}]\label{equicon}
The theories $\ZF$+$\ADR$ and $\ZF$+$\AC_{\omega} (\R)$+$\BlADR$ are equiconsistent.
\end{Thm}

\begin{proof}
We use the following result of Trang and Wilson~\cite[Theorem~1.6]{Nam-Trevor}:
\begin{Thm}[Trang and Wilson~\cite{Nam-Trevor}]\label{nam-trevor}
The following theories are equiconsistent:
\begin{enumerate}
\item $\ZF$+$\ADR$.

\item $\ZF$+$\DC_{\wp (\omega_1)}$+\lq\lq There is a fine measure on $\wp_{\omega_1} (\R)$ and $\Theta$ is singular''. 
\end{enumerate}
\end{Thm}

By Theorem~\ref{nam-trevor}, it is enough to show that assuming $\ZF$+$\AC_{\omega} (\R)$+$\BlADR$, there is a model of $\ZF$+$\DC_{\wp (\omega_1)}$+\lq\lq There is a fine measure on $\wp_{\omega_1} (\R)$ and $\Theta$ is singular''. 
By Theorem~\ref{R sharp}, $\BlADR$ implies that there is a fine measure on $\wp_{\omega_1}(\R)$. 
We also note that $\BlADR$ implies $\DC_{\wp (\omega_1)}$:  
By Theorem~\ref{uniformization}, $\BlADR$ implies $\DC_{\R}$. 
By the weak version of Moschovakis' Lemma obtained in Theorem~\ref{weak coding}, $\BlAD$ implies that there is a surjection from $\R$ to $\wp (\omega_1)$ and so $\DC_{\R}$ implies $\DC_{\wp (\omega_1)}$ while by Proposition~\ref{easy_Blackwell}, $\BlADR$ implies $\BlAD$. All in all, $\BlADR$ implies $\DC_{\wp (\omega_1)}$. 
Therefore, it is enough to verify that assuming $\ZF$+$\AC_{\omega} (\R)$+$\BlADR$, there is a model of $\ZF$+$\BlADR$+\lq\lq $\Theta$ is singular''. 

We will prove something stronger than needed: For a set of reals $A$, by $\Theta_A$ we mean the supremum of ordinals $\delta$ such that there is a surjection $\pi \colon \R \to \delta$ so that $\pi$ is ordinal definable with the parameters $A$ and a real. 
\begin{Thm}\label{BLADR-cofomega}
Assume $\BlADR$. Then there is an ordinal definable function $f \colon \Theta \to \Theta$ such that if $\gamma < \Theta$ is closed under $f$, then letting $\Gamma_{\gamma} = \{ A \subseteq \R \mid \Theta_A < \gamma \}$, $M = \text{HOD}_{\R \cup \Gamma_{\gamma}} $ is a model of $\BlADR$ and $\gamma = \Theta^{M}$. In particular, $\BlADR$ implies there is an inner model of $\ZF$+$\BlADR$+\lq\lq $\Theta$ is of countable cofinality''.  
\end{Thm}

\begin{proof}[Proof of Theorem~\ref{BLADR-cofomega}]
The arguments are similar to those by Solovay~\cite[Theorem~2.5]{Solovay_AD_R}. We simulate Solovay's arguments using $\BlADR$ instead of $\ADR$. 

We first argue that for all sets of reals $A$, $\Theta_A < \Theta$. 
Take any set of reals $A$. Let $C = \{ (x,y) \in \R \times \R \mid $ $y$ is not ordinal definable from $A$ and $x\}$. By $\BlAD$ and Corollary~\ref{cor:BlAD-measure}, 
there is no injection from $\omega_1$ to the reals. Hence for any real $x$, there is a real $y$ with $(x,y) \in C$. 
By $\BlADR$ and Theorem~\ref{uniformization}, there is a function $F$ which uniformizes $C$. 
By the argument of Solovay~\cite[Lemma~2.2]{Solovay_AD_R}, $F$ is not ordinal definable from $A$ and a real. 
By $\BlAD$ and the weak version of Wadge Lemma (Lemma~\ref{weak Wadge}), for all sets of reals $B$ which are ordinal definable from $A$ and a real, $B$ is $\undertilde{\mathbf{\Sigma}}^1_2$ in $F$. Hence, using $F$, one can construct a surjection from $\R$ to the collection of all sets of reals which are ordinal definable from $A$ and a real, which gives us a surjection from $\R$ to $\Theta_A$, because for each ordinal $\alpha < \Theta_A$, there is a prewellordering $R$ on the reals which is ordinal definable from $A$ and a real such that the length of $R$ is $\alpha$. Therefore, $\Theta_A < \Theta$. 

We next claim that for all $\alpha < \Theta$, there is a set of reals $A$ such that $\alpha < \Theta_A$. This is easy because given any ordinal $\alpha < \Theta$, one can take $A$ as a prewellordering on $\R$ whose length is greater than $\alpha$. 

From the last two paragraphs, one can define the following function $g \colon \Theta \to \Theta$: For any $\alpha < \Theta$, let $g(\alpha) = \min \{ \Theta_A \mid \alpha < \Theta_A \}$. Then $g$ is ordinal definable. 

The following is the key point: We identify sets of reals with subsets of $\R^{\omega}$ via a simple bijection between $\R$ and $\R^{\omega}$. 
\begin{Claim}\label{coding optimal strategies}
Let $A$ be a set of reals. Then there is a set of reals $C$ such that for all $B \subseteq \R^{\omega}$ which is $\undertilde{\mathbf{\Sigma}}^1_2$ in $A$, one of the players has an optimal strategy $\sigma$ in $B$ with $\sigma \in \mathrm{L}(\R, C)$. 
\end{Claim}

\begin{proof}[Proof of Claim~\ref{coding optimal strategies}]

Take any set of reals $A$. 
Consider the following game $\tilde{\mathcal{G}}_{\mathrm{Bl}}$: There are two players. Player I first chooses a pair $(x, \phi)$ where $x$ is a real and $\phi$ is a $\Sigma^1_2$-formula such that $A$ and $x$ define a subset $B$ of $\R^{\omega}$  via the formula $\phi$. Next Player II chooses \lq I' or \lq II'. Then they play the Blackwell game $\mathcal{G}_{\mathrm{Bl}}(B)$ with the payoff set $B$ where Player II becomes the first player of $\mathcal{G}_{\mathrm{Bl}}(B)$ if she chose \lq I', otherwise Player II becomes of the second player of $\mathcal{G}_{\mathrm{Bl}} (B)$. 

By $\BlADR$, for all $B \subseteq \R^{\omega}$, one of the players has an optimal strategy in the game $\mathcal{G}_{\mathrm{Bl}}(B)$. So Player I cannot have an optimal strategy in the game $\tilde{\mathcal{G}}_{\mathrm{Bl}}$. Hence by $\BlADR$, Player II has an optimal strategy $\tau$ in the game $\tilde{\mathcal{G}}_{\mathrm{Bl}}$. 
Also, notice that for each pair $(x, \phi)$ such that $A$ and $x$ define a subset $B$ of $\R^{\omega}$ via the formula $\phi$, one can compute from $\tau$ and $(x, \phi)$ an optimal strategy in the Blackwell game $\mathcal{G}_{\mathrm{Bl}}(B)$ in a simple manner. Hence any set of reals $C$ coding $\tau$ in a simple manner is the desired one.
\renewcommand{\qedsymbol}{$\square \ (\text{Claim~\ref{coding optimal strategies}})$}
\end{proof}

By Claim~\ref{coding optimal strategies}, one can define the following function $h \colon \Theta \to \Theta$: Given an $\alpha < \Theta$, let $h(\alpha) $ be the least ordinal of the form $\Theta_C$ such that for all $B\subseteq \R^{\omega}$ with $\Theta_B < \alpha$, one of the players has an optimal strategy in $\mathrm{L}(\R , C)$. Such a $C$ exists for any $\alpha < \Theta$ because given a prewellordering $A$ on the reals of length $\alpha$, $A$ is not ordinal definable from any $B \subseteq \R^{\omega}$ with $\Theta_B < \alpha$ and a real, so by the weak version of Wadge Lemma (Lemma~\ref{weak Wadge}), any such $B$ is $\undertilde{\mathbf{\Sigma}}^1_2$ from $A$. Now one can apply Claim~\ref{coding optimal strategies} to $A$ to obtain the desired $C$. 
Notice that $h$ is ordinal definable. 

Let $f \colon \Theta \to \Theta$ be such that $f(\alpha) = \sup \{ g(\alpha) , h (\alpha)\}$. 
Then since $g$ and $h$ are ordinal definable, so is $f$. 
Let $\gamma$ be any ordinal closed under $f$. We show that $M =\text{HOD}_{\R \cup \Gamma_{\gamma}}$ is a model of $\BlADR$ and $\gamma = \Theta^{M}$.  

We first verify that $\Gamma_{\gamma}  = \wp (\R) \cap \text{HOD}_{\R \cup \Gamma_{\gamma}}$. 
The inclusion $\Gamma_{\gamma}  \subseteq \wp (\R) \cap \text{HOD}_{\R \cup \Gamma_{\gamma}}$ is obvious. 
We argue for the reverse inclusion. Let $A$ be a set of reals in $\text{HOD}_{\R \cup \Gamma_{\gamma}}$. We see that $A$ is in $\Gamma_{\gamma}$. Since $A$ is in $\text{HOD}_{\R \cup \Gamma_{\gamma}}$, there is a set of reals $B$ such that $B$ is in $\Gamma_{\gamma}$ and $A$ is ordinal definable from $B$ and a real. Since $B$ is in $\Gamma_{\gamma}$, $\Theta_B < \gamma$. But since $A$ is ordinal definable from $B$ and a real, $\Theta_A \le \Theta_B$. So $\Theta_A < \gamma$ and $A$ is in $\Gamma_{\gamma}$, as desired. 

We next verify that $\gamma = \Theta^{M}$, where $M = \text{HOD}_{\R \cup \Gamma_{\gamma}}$. 
We first see $\gamma \le \Theta^M$. 
Let $\alpha < \gamma$. We show that $\alpha < \Theta^{M}$.  
Since $\alpha < g(\alpha) \le f(\alpha) < \gamma$, there is a set of reals $A$ such that $\alpha < \Theta_A < \gamma$. So $A$ is in $\Gamma_{\gamma}$. Since $\alpha < \Theta_A$, there is a surjection $\pi \colon \R \to \alpha$ which is ordinal definable with parameters $A$ and a real. So there is a prewellordering $B$ on the reals of length $\alpha$ which is ordinal definable with parameters $A$ and a real. Since $\Theta_B \le \Theta_A <\gamma$, $B$ is in $\Gamma_{\gamma}$, and hence $\alpha < \Theta^{M}$, as desired. 
We next see $\Theta^M \le \gamma$. 
Let $\beta < \Theta^{M}$. We show that $\beta < \gamma$. Since $\beta < \Theta^{M}$, there is a prewellordering $C$ on the reals of length $\beta$ such that $C$ is in $M = \text{HOD}_{\R \cup \Gamma_{\gamma}}$. Since $\Gamma_{\gamma}  = \wp (\R) \cap \text{HOD}_{\R \cup \Gamma_{\gamma}}$, $C$ is in $\Gamma_{\gamma}$. So $\Theta_C < \gamma$ while $\beta < \Theta_C$ because $\beta$ is the length of $C$. Hence $\beta < \gamma$, as desired. 

We now verify that $M = \text{HOD}_{\R \cup \Gamma_{\gamma}}$ is a model of $\BlADR$. Let $A$ be any subset of $\R^{\omega}$ in $\text{HOD}_{\R \cup \Gamma_{\gamma}}$. We show that one of the players has an optimal strategy in the Blackwell game $\mathcal{G}_{\mathrm{Bl}} (A)$ in the model $\text{HOD}_{\R \cup \Gamma_{\gamma}}$. Since $\Gamma_{\gamma}  = \wp (\R) \cap \text{HOD}_{\R \cup \Gamma_{\gamma}}$, $A$ is in $\Gamma_{\gamma}$. So $\Theta_A < \gamma$. Since $h (\Theta_A) \le f(\Theta_A) < \gamma$, by the definition of $h(\Theta_A)$, there is a set of reals $C$ with $\Theta_C < \gamma$ such that one of the players has an optimal strategy $\sigma$ in the Blackwell game $\mathcal{G}_{\mathrm{Bl}} (A)$ with $\sigma \in \mathrm{L}(\R ,C)$. Since $\Theta_C < \gamma$, $C$ is in $\Gamma_{\gamma}$ and hence $\sigma \in \mathrm{L}(\R , C) \subseteq \text{HOD}_{\R \cup \Gamma_{\gamma}}$, as desired.

All in all, we have verified that $f$ is the desired function. 

Lastly, we verify that $\BlADR$ implies that there is an inner model of $\ZF$+$\BlADR$+ \lq\lq $\Theta$ is of countable cofinality''. 
Let $f \colon \Theta \to \Theta$ be the ordinal definable function we defined above. 
Let $(\alpha_n \mid n < \omega)$ be such that $\alpha_0 = 0$ and for all $n < \omega$, $\alpha_{n+1} = f( \alpha_n)$. Set $\gamma_0 = \sup_{n < \omega} \alpha_n$. Then $\gamma_0$ is closed under $f$. Let $N = \text{HOD}_{\R \cup \Gamma_{\gamma_0}}$. 
Then by the property of $f$, $N$ is an inner model of $\ZF$+$\BlADR$ such that $\gamma_0 = \Theta^N$. 
We argue that $\gamma_0 = \Theta^N$ is of countable cofinality in $N$.  
Since $f$ is ordinal definable, so is the sequence $(\alpha \mid n < \omega)$, and hence $\gamma_0 = \sup_{n < \omega} \alpha_n$ is of countable cofinality in $\text{HOD}$. Since $\text{HOD} \subseteq \text{HOD}_{\R \cup \Gamma_{\gamma_0}} = N$, $\gamma_0 = \Theta^N$ is of countable cofinality in $N$, as desired. 
\renewcommand{\qedsymbol}{$\square \ (\text{Theorem~\ref{BLADR-cofomega}})$}
\end{proof}
This completes the proof of Theorem~\ref{equicon}. 
\end{proof}

\section{Questions}\label{sec:Q}

We close this paper by raising some open questions:
\begin{Q}\label{Q1}
Are $\AD$ and $\BlAD$ equivalent?
\end{Q}

Martin conjectured that the answer to Question~\ref{Q1} is \lq Yes' (Conjecture~\ref{Martin-conjecture}). 

\begin{Q}
Are $\ADR$ and $\BlADR$ equivalent without assuming $\DC$?
\end{Q}

As in Theorem~\ref{main theorem}, $\ADR$ and $\BlADR$ are equivalent in $\ZF$+$\DC$. 



\bibliographystyle{plain}
\bibliography{myreference}

\end{document}